\documentclass{amsart}
\usepackage{amscd}
\usepackage{graphicx}
\usepackage{epstopdf}
\numberwithin{equation}{section}

\let\cal\mathcal

\def\id{\text{id}}

\def\Sym{\mathop{\text{\upshape{Sym}}}\nolimits}

\def\id{\operatorname {id}}

\def\r{\rightarrow}

 \DeclareMathOperator{\Aut}{Aut}

\newtheorem{lemma}{Lemma}[section]
\newtheorem{proposition}[lemma]{Proposition}
\newtheorem{theorem}[lemma]{Theorem}

\newtheorem{openquestion}[lemma]{Open Question}

\newtheorem{corollary}[lemma]{Corollary}

\theoremstyle{definition}

\newtheorem{example}[lemma]{Example}
\newtheorem{definition}[lemma]{Definition}

\newtheorem{notation}[lemma]{Notation}

{

}

\theoremstyle{remark}

\newtheorem{remark}[lemma]{Remark}

\newcommand{\Acal}{\mbox{$\cal A$}}

\newcommand{\Dcal}{\mbox{$\cal D$}}

\newcommand{\Lcal}{\mbox{$\cal L$}}

\newcommand{\Ocal}{\mbox{$\cal O$}}

\newcommand{\Rcal}{\mbox{$\cal R$}}

\newcommand{\mll}{{\rm\bf ml2}}
\newcommand{\mrr}{{\rm\bf mr2}}
\newcommand{\mlw}{{\rm\bf ml2w}}
\newcommand{\mrw}{{\rm\bf mr2w}}
\newcommand{\mla}{{\rm\bf ml1a}}
\newcommand{\mra}{{\rm\bf mr1a}}

\newcommand{\U}{{U}}

\newcommand{\la}{{\triangleright}}
\newcommand{\ra}{{\triangleleft}}

\newdimen\uboxsep \uboxsep=1ex
\def\uboxn#1{\vtop to 0pt{\hrule height 0pt depth 0pt\vskip\uboxsep
\hbox to 0pt{\hss #1\hss}\vss}}

\def\uboxs#1{\vbox to 0pt{\vss\hbox to 0pt{\hss #1\hss}
\vskip\uboxsep\hrule height 0pt depth 0pt}}

\title[Matched pairs approach to YBE]{Matched pairs approach to set theoretic solutions of the Yang-Baxter equation} 
\thanks{The first author was partially supported by the Royal Society, UK, by 
Grant MI 1503/2005 of the Bulgarian National Science Fund
  of the Ministry of Education and Science, and
  by The Abdus Salam International Centre
   for Theoretical Physics (ICTP). Both authors were supported by the Isaac Newton Institute during completion of the work.}
\author{Tatiana Gateva-Ivanova and Shahn Majid}
\address{TGI: Institute of Mathematics and Informatics\\
Bulgarian Academy of Sciences\\
Sofia 1113, Bulgaria\\
S.M:  Queen Mary, University of London\\
School of Mathematics, Mile End Rd, London E1 4NS, UK}

\date{21.11.2006}

\email{tatianagateva@yahoo.com, tatyana@aubg.bg,
s.majid@qmul.ac.uk}

\begin{document}

\begin{abstract} We study set-theoretic solutions $(X,r)$ of the Yang-Baxter
equations on a set $X$ in terms of the induced left and right
actions of $X$ on itself. We give a characterization of involutive
square-free solutions in terms of cyclicity conditions. We 
characterise general solutions in terms of abstract matched pair properties of 
the associated monoid  $S(X,r)$ and we show that $r$ extends as a solution $(S(X,r),r_S)$.
 Finally, we study extensions of solutions both directly and in terms
of matched pairs of their associated monoids. We also prove
several general results about matched pairs of monoids $S$ 
of the required type,  including iterated products $S\bowtie
S\bowtie S$ equivalent to $r_S$ a solution, and
extensions $(S\bowtie T,r_{S\bowtie T})$. Examples  include a general
 `double' construction $(S\bowtie S,r_{S\bowtie S})$ 
 and some concrete extensions, their actions and graphs based on small sets.
\end{abstract}

\maketitle

\section{Introduction }

It was established in the last two decades that solutions of the  linear braid or `Yang-Baxter equations' (YBE) 
\[ r^{12}r^{23}r^{12}=r^{23}r^{12}r^{23}\]
on a vector space of the form $V^{\otimes3}$  lead to very remarkable algebraic structures. Here $r:V\otimes V\to V\otimes V$,  $r^{12}=r\otimes\id, r^{23}=\id\otimes r$ is a notation and structures include coquasitriangular bialgebras $A(r)$, their quantum group (Hopf algebra) quotients, quantum planes and associated objects at least in the case of specific standard solutions, 
see \cite{R,Ma:book}. On the other hand, the variety of all solutions on vector spaces of a given dimension has remained rather elusive in any degree of generality. Some years ago it was proposed by V.G. Drinfeld\cite{D} to consider the same equations in the category of sets, and in this setting several general results are now known. It is clear that a set-theoretic solution extends to a linear one, but more important than this is that set-theoretic solutions lead to their own remarkable algebraic and combinatoric structures, only somewhat  analogous to quantum group constructions.
 In the present paper we develop a systematic approach to set-theoretic solutions based on the monoid that they generate and on a theory of matched pairs or  factorisation of monoids. We present results both at the level of $(X,r)$, at the level of unital semigroups (monoids) and on the `exponentiation problem' from one to the other. 

 Let $X$ be a
nonempty set and let $r: X\times X \longrightarrow X\times X$ be a
bijective map. In this case  we shall use notation $(X,r)$ and
refer to it as a \emph{quadratic set} or set with quadratic map $r$. We present the
image of $(x,y)$ under $r$ as
\begin{equation}
\label{r} r(x,y)=({}^xy,x^{y}).
\end{equation}
The formula (\ref{r}) defines a ``left action" $\Lcal: X\times X
\longrightarrow X,$ and a ``right action" $\Rcal: X\times X
\longrightarrow X,$ on $X$ as:
\begin{equation}
\label{LcalRcal} \Lcal_x(y)={}^xy, \quad \Rcal_y(x)= x^{y},
\end{equation}
for all $ x, y \in X.$ The map  $r$ is \emph{non-degenerate}, if the maps $\Rcal_x$ and $\Lcal_x$ are
bijective for each $x\in X$.   In this paper we shall  often assume that $r$ is non-degenerate,  as
will be indicated. Also, as a notational tool, we  shall often
identify the sets $X\times X$ and $X^2,$ the set of
all monomials of length two in the free monoid $\langle
X\rangle.$

 As in \cite{T06} to each quadratic map $r: X^2 \r X^2$ we
associate canonically algebraic objects (see
Definition~\ref{associatedobjects}) generated by $X$ and with
quadratic defining relations $\Re =\Re(r)$ naturally determined as
\begin{equation}
\label{defrelations} xy=zt \in \Re(r),\quad
  \text{whenever}\quad r(x,y) = (z,t).
\end{equation}
Note that in the case when $X$ is finite, the set $\Re(r)$  of
defining
 relations is also finite, therefore the associated algebraic
 objects are finitely presented. Furthermore in many cases they
 will be standard finitely presented with respect to the
 degree-lexicographic ordering.  It is known in particular that the algebra generated by the monoid
 $S(X,r)$ defined in this way  has remarkable homological properties when $r$ obeys the braid or Yang-Baxter equations
 in $X\times X\times X$ and other restrictions such as square-free, involutive. Set-theoretic solutions were introduced in  \cite{D,W} and studied in \cite{T1,T2,T3,TM,H,ESS,Lu,JO,Laffaille} as well as more recently in \cite{TJO,TLovetch,E1,T06,T07,Rump} and other works. We say that a set $X$ is `braided' when it is equipped with an invertible solution.

 In this paper we study the  close relations between the combinatorial properties
of the defining relations, i.e. of the map $r$, and the structural
 properties of the associated algebraic objects,  particularly through the `actions'
 above. Section 2 contains preliminary material and some elementary results based on direct calculation with $\Rcal_x,\Lcal_x$. These include cancellation properties in $S(X,r)$ (Proposition~2.13 and related results) that
 will be needed later and a full characterisation of when $(X,r)$ obeys the Yang-Baxter equations
 under the assumption of square free, involutive and nondegenerate. In particular we find that
 they are equivalent to a cyclicity condition
 \[  {}^{({}^yx)}({}^yz)={}^{({}^xy)}({}^xz)\]
 for all $x,y,z\in X$ (see Theorem~2.36). 
Cyclicity conditions, originally called ``cyclic conditions"
 were  discovered in 1990 while the first author
 was looking for new classes of Noetherian Artin-Shelter regular rings, \cite{T1,T2,T3}, and have
 already played an important role in the theory, see \cite{T1,T2,T3,T06,T07,TM,TJO,JO,Rump}.
 The cyclic conditions are used in  the proof of the nice algebraic and homological
 properties of the binomial skew polynomial rings,
 such as being Noetherian, \cite{T1,T2}, Gorenstein (Artin-Schelter regular), \cite{T3,T07}, and for their close relation
 with the set-theoretic solutions of YBE \cite{TM, T07}. The cycle sets, as Rump calls them have essential role
 in his decomposition theorem,  \cite{Rump}.

 In Section 3 we show that when $(X,r)$ is a braided set the `actions' $\Rcal_x,\Lcal_x$ indeed extend to actions of the  monoid $S=S(X,r)$ on itself to form a matched pair $(S,S)$ of monoids in the sense of \cite{Ma:phy,Ma:book}. This means that $S$ acts on itself from the left and the right in a compatible way to construct a new monoid $S\bowtie S$. We also have an induced map $r_S(u,v)=({}^uv,u^v)$ and moreover  \[ uv={}^uv u^v\]
holds in $S$ for all $u,v\in S$ (in other words, $r_S$-commutative). This result is in Theorem~\ref{theoremA} and leads us also to introduce a general theory of monoids $S$ with these properties, which we call an `{\bf M3}-monoid'. We remark that matched pairs of groups are a notion known in group theory since the 1910s and also known \cite{Lu}  to be connected with set-theoretic   solutions of the YBE. There is also a nice review of the main results in \cite{Lu,ESS} from a matched pair of groups point of view in \cite{Takeuchi}, among them a good way to think about the properties of the group $G(X,r)$ universally generated from $(X,r)$. In particular it is known that the group $G(X,r)$ is itself a braided set  in an induced manner and this is the starting point of the these works. Our own goals and methods are certainly very different and in particular for an {\bf M3}-monoid one cannot automatically deduce that $r_S$ is invertible and obeys the YBE as would be true in the group case. When it does, we say that $(S,r_S)$ is a `braided monoid' and  we eventually show in Theorem~\ref{theoremAextra} that this is nevertheless the case for $S=S(X,r)$. Most of this section is concerned with these two theorems. Theorem~\ref{S,S matchedth} at the end of the section provides the converse under a  mild 2-cancellative assumption on the quadratic set $(X,r)$; this is braided if and only if $(S(X,r),r_S)$ is a braided monoid.

Among results at the level of a general {\bf M3}-monoid $S$ we give some explicit criteria for $r_S$ to obey the YBE  and a complete characterisation for this in Theorem~\ref{tricrossedproduct} as $(S, S\bowtie S)$ forming a matched pair by extending the actions. In this case we find iterated double crossproduct monoids $S\bowtie S$, $S\bowtie S\bowtie S$ etc., as the analogue of the construction for bialgebras $A(r)$ of iterated double cross products  $A(r)\bowtie A(r), A(r)\bowtie A(r)\bowtie A(r)$, etc. in \cite{Ma:mor}.  We also provide an interpretation of the YBE as zero curvature for  `surface-transport' around a cube and this is related to the `labelled square' construction for the bicrossproduct Hopf algebras associated to (finite) group matched pairs introduced in \cite{Ma:book}.

 In Section 4 we use our matched pair characterisation to study a natural notion of `regular extension'  $Z=X\sqcup Y$ of set theoretic solutions. Our first main result, Theorems~\ref{propositionBZ} and~\ref{theoremBcancellative} provide explicit if and only if conditions for such an extension to obey the YBE.  Theorem~\ref{theoremB} finds similarly explicit but slightly weaker if and only if conditions for a matched pair $(S,T)$ and monoid $S\bowtie T$ where $S,T$ are the associated braided monoids built from $X,Y$ using the results of Section~3. Our inductive construction from the $(X,r_X)$ and $(Y,r_Y)$ data here (and also the construction of $(S,S)$ in Section~3) is  somewhat analogous to the construction of a Lie group matched pair from a Lie algebra one by integrating vector fields and connections \cite{Ma:mat}.

Finally in Section 4 we bring these results together with a complete if and only if characterisation of when  is a regular extension obeys the YBE as equivalent to $U=S(Z,r_Z)$ being a double cross product $U=S\bowtie T$ and forming matched pairs with $S,T$. This is Theorem~\ref{bellaT2} and is the main result  of the section. We also provide a corresponding  theory of regular extensions at the level of abstract {\bf M3}-monoids in Theorem~\ref{M3exttheorem} and Corollary~\ref{braextcor}. A further immediate corollary is that every braided monoid $S$ has a `double' braided monoid $S\bowtie S$, which is an analogue of the induced coquasitriangular structure on the bialgebra $A(R)\bowtie A(R)$ in \cite{Ma:book}.  
 
 Section~5 looks at explicit constructions for extensions using our methods for a pair of initial $(X,r_X)$, $(Y,r_Y)$. We consider particularly a special case which we call `strong twisted union' and which is a nevertheless a far-reaching generalisaton of  the `generalised twisted union' in \cite{ESS}. We give characterizations of when a strong twisted union obeys the YBE linking back to our cyclicity results of Section~2. We conclude with some  concrete nontrivial examples of extensions. The solutions here are  involutive, nondegenerate, square-free with the condition `lri'  also from Section~2. We provide  graphs of the initial and extended solutions  and illustrate how the methods of the paper translate into very practical tools for their construction. These graphs can  also be viewed as a natural source of discrete noncommutative geometries, a point of view to be developed elsewhere. 

 \section{Cancellation and cyclicity properties for quadratic sets}

The new results in this section relate to 2-cancellation conditions needed later on
and to cyclicity conditions of interest in the theory of quadratic sets. We show in particular that
a number of possible such conditions are all equivalent under mild hypotheses. 

\subsection{Preliminaries on quadratic and braided sets}

 For a non-empty set $X$, as usual, we denote
by $\langle X \rangle$,  and ${}_{gr}\langle X \rangle$,
respectively the free unital semigroup (i.e., free monoid), and the free group
generated by $X,$ and by $k\langle X \rangle$-
 the free associative $k$-algebra
generated by $X$, where $k$ is an arbitrary field. For a set $F
\subseteq k\langle X \rangle$, $(F)$ denotes the two sided ideal
of $k\langle X \rangle$, generated by $F$.

\begin{definition} \cite{T06}
\label{associatedobjects} Assume that $r:X^2 \longrightarrow X^2$
is a bijective map.

(i) The monoid
\[
S =S(X, r) = \langle X; \Re(r) \rangle,
\]
 with a
set of generators $X$ and a set of defining relations $ \Re(r),$
is called \emph{the monoid associated with $(X, r)$}.

(ii) The \emph{group $G=G(X, r)$ associated with} $(X, r)$ is
defined as
\[
G=G(X, r)={}_{gr} \langle X; \Re (r) \rangle.
\]

(iii) For arbitrary fixed field $k$, \emph{the $k$-algebra
associated with} $(X ,r)$ is defined as
\begin{equation}
\label{Adef} \Acal = \Acal(k,X,r) = k\langle X \rangle/(\Re(r)).
\end{equation}
\end{definition}
Clearly $\Acal$ is a quadratic algebra, generated by $X$ and
 with defining relations  $\Re(r).$
Furthermore, $\Acal$  is isomorphic to the monoid algebra
$kS(X, r).$

\begin{remark}
\label{orbitsinG}
 \cite{T07}
When we study the monoid $S=S(X,r)$, the group $G=G(X,r)$,  or
the algebra $A = A(k, X,r)$ associated with $(X,r)$,
it is  convenient to use the action of the infinite group,
$\Dcal(r)$, generated by maps associated with the quadratic
relations, as follows. As usual, we consider the two bijective
maps $r^{ii+1}: X^3 \longrightarrow X^3$, $ 1 \leq i \leq 2$ ,
where $r^{12} = r\times Id_X$, and $r^{23}=Id_X\times r.$ Then
\[
\Dcal= \Dcal(r)=\: _{\rm{gr}} \langle r^{12}, r^{23} \rangle
\]
acts on $X^3.$ If  $r$ is involutive, the bijective maps $r^{12}$
and $r^{23}$ are involutive as well, so in this case $\Dcal(r)$ is
\emph{the infinite dihedral group},
\[
\Dcal= \Dcal(r)=\: _{\rm{gr}} \langle r^{12}, r^{23}: (r^{12})^2= e,\quad  (r^{23})^2=e \rangle
\]
Note that all monomials $\omega^{\prime}\in X^3$  which belong to
the same orbit $\Ocal_{\Dcal}(\omega)$ satisfy $\omega^{\prime} =\omega$ as elements of $G$ (respectively, $S$, $A$).
\end{remark}
\begin{definition}
\label{defvariousr}
\begin{enumerate}
 \item
$r$ is \emph{square-free} if $r(x,x)=(x,x)$ for all $x\in X.$
\item A non-degenerate involutive square-free map $(X,r)$ will be
called \emph{(set-theoretic) quantum binomial map}. \item
\label{YBE} $r$ is \emph{a set-theoretic solution of the
Yang-Baxter equation}
 (YBE) if  the braid relation
\[r^{12}r^{23}r^{12} = r^{23}r^{12}r^{23}\]
holds in $X\times X\times X.$  In this case $(X,r)$ is also called
\emph{a braided set}. If in addition $r$ is involutive $(X,r)$ is
called \emph{a symmetric set}.
\end{enumerate}
\end{definition}
The following lemma is  straightforward.
\begin{lemma}
\label{someproperties} Suppose  $(X,r)$ is given, and let
${}^x\bullet$, and ${\bullet}^x$ be the associated left and right
actions. Then
\begin{enumerate}
\item \label{involutive} $r$ is involutive if and only if
\[
{}^{{}^xy}{(x^y)} = x, \; \text{and}\;
({}^yx)^{y^x} = x, \; \text{for all}\; x,y \in X.\]

 \item \label{SF} $r$ is square-free if and only if
\[
{}^xx=x,\; \text{and} \; x^x = x \; \text{for all}\; x \in X.
\]
\item \label{SFnondeg} If $r$ is non-degenerate and square-free,
then
\[
{}^xy=x\ \Longleftrightarrow\ x^y = y  \Longleftrightarrow\  y=x\ \Longleftrightarrow\ r(x,y)=(x,y).\]
\end{enumerate}
\end{lemma}

It is also straightforward to write out the Yang-Baxter equations
for $r$ in terms of the actions. This is in  \cite{ESS}  but we
recall it here in our notations for convenience.
\begin{lemma}\label{ybe} Let $(X,r)$ be given in the notations above. Then $r$ obeys the YBE 
(or $(X,r)$ is a braided set) {\em iff} the following conditions
hold
\[
\begin{array}{lclc}
 {\bf l1:}\quad& {}^x{({}^yz)}={}^{{}^xy}{({}^{x^y}{z})},
 \quad\quad\quad
 & {\bf r1:}\quad&
{(x^y)}^z=(x^{{}^yz})^{y^z},
\end{array}\]
 \[ {\rm\bf lr3:} \quad
{({}^xy)}^{({}^{x^y}{z})} \ = \ {}^{(x^{{}^yz})}{(y^z)},\]
 for all $x,y,z \in X$.
 \end{lemma}
  \begin{proof} We refer to the diagram
(\ref{ybediagram}). This diagram contains elements of the orbit of
arbitrary monomial $xyz \in X^3$, under the action of the group
$\Dcal(r).$ Note that each two monomials in this orbit are equal
as elements of $S$
\begin{equation}
\label{ybediagram}
\begin{CD}
 xyz @>r_{12}>> r(xy)z=({{}^xy}x^y)z\\
@V  r_{23} VV @VV r_{23} V\\
  xr(yz)=x({{}^yz}y^z)@. ({}^xy)r((x^y)z)=({{}^xy}) ({}^{x^y}z)(x^y)^z\\
@V r_{12} VV @VV r_{12} V\\
\kern -20pt r(x({{}^yz}))(y^z)={}^{x}{({}^yz)}(x^{{}^yz})(y^z)@.
\kern-80pt r(({}^xy)({}^{x^y}z))(x^y)^z=[{}^{{}^xy}{({}^{x^y}z)}][({}^xy)^{({}^{x^y}z)}][(x^y)^z]=w_1\\
  @V r_{23} VV  \\
 {}^{x}{({}^yz)}r((x^{{}^yz})(y^z))=[{}^{x}{({}^yz)}][{}^{(x^{{}^yz})}{(y^z)}]
[(x^{{}^yz})^{y^z}]=w_2@.
 \end{CD}
\end{equation}

from which we read off {\bf l1,lr3,r1} for equality of the words
in $X^3$.
 \end{proof}

 The following proposition was also proved in \cite{ESS}, see 2.1.

\begin{proposition}
\label{ESS1} \cite{ESS}.  (a) Suppose that $(X,r)$ is involutive,
$r$ is right non-degenerate and the assignment $x \rightarrow
\Rcal_x$ is a right action of $G(X,r)$ on $X$. Define the map $T:
X \longrightarrow X$ by the formula $T(y) = \Rcal_y^{-1}(y)$, then
one has $\Rcal_x^{-1}\circ T= T\circ \Lcal_x$.

Suppose in addition that  $\Lcal_x$ are invertible , i.e. $(X,r)$
is involutive and non-degenerate. Then

(b) $T$ is invertible, and the two left actions of $G(X,r)$ given
by $x \rightarrow \Rcal_x^{-1}$ and $x \rightarrow \Lcal_x$ are
isomorphic to each other.

(c) Condition {\bf r1} implies {\bf l1} and {\bf lr3}.
Thus $(X,r)$ is symmetric if and only if {\bf r1} is satisfied.
\end{proposition}
\begin{remark}
 It was shown by a different argument in  \cite{T06} that if (X,r)
is square-free symmetric set, then $\Rcal_x = \Lcal_x^{-1}$ for
all $x\in X$
\end{remark}
We want to find out now what is the relation between  the actions,
and the above conditions,  if we assume only that $r$ in
non-degenerate, without any further restrictions. Compatible left and right actions appear in the  notion of matched
pairs of groups, see \cite{Ma:book}, and motivated by this  theory we also define
their natural extensions:

\begin{definition}\label{extendedleftaction} Given $(X,r)$ we extend the  actions ${}^x\bullet$ and  $\bullet
^x$ on $X$ to left and right actions on $X\times X$ as follows.
For $x,y,z \in X$ we define:
\[ {}^x{(y,z)}:=({}^xy,{}^{x^y}z), \quad
\text{and} \quad (x,y)^z:= (x^{{}^yz}, y^z)\]
We say that $r$ is respectively left and right invariant if
\[{\bf l2}: \; \; \;  \; \; \;
r({}^x{(y,z)})={}^x{(r(y,z))},\quad {\bf r2}: \; \; \; \; \; \;
r((x,y)^z)={(r(x,y))}^z\]
hold for all $x,y,z\in Z$.
\end{definition}

\begin{lemma}
\label{leml1r2} Let $(X,r)$ be given in the notations above. Then the following are equivalent
\begin{enumerate}
\item  $(X,r)$ is a braided set.
\item {\bf l1}, {\bf r2}.
\item  {\bf r1}, {\bf l2}.
\end{enumerate}
\end{lemma}
\begin{proof}  From the  formula (\ref{r}) and Definition~\ref{extendedleftaction}: 
\begin{equation} 
\label{xr}
r({}^x{(y,z)})=r({}^xy, {}^{(x^y)}z) = ({}^{{}^xy}{({}^{(x^y)}z)}, ({}^xy)^{({}^{(x^y)}z)}).
 \end{equation}
\begin{equation} 
\label{rx}
 {}^x{(r(y,z))}={}^{x}{({}^yz,y^z)} = ({}^{x}{({}^yz)},{}^{(x^{({}^yz}))}{(y^z)}).
 \end{equation}
Hence condition \textbf{l2} is an equality of pairs in $X\times
 X$:
\[
({}^{{}^xy}{({}^{(x^y)}z)}, ({}^xy)^{({}^{(x^y)}z)})=({}^{x}{({}^yz)},{}^{(x^{({}^yz)})}{(y^z)})
\]
as required.  The case of \textbf{r2} is similar. We have proved that for any quadratic set
\begin{equation}\label{l2l1}
{\bf l2} \Longleftrightarrow {\bf l1,lr3};\quad 
{\bf r2} \Longleftrightarrow {\bf r1,lr3}
\end{equation}
and we then use Lemma~\ref{ybe}.  It is also easy enough to see the result directly from the diagram (\ref{ybediagram}). Thus 
\begin{equation}\label{w1w2} w_1={}^{{}^xy}({}^{x^y}z)r(x,y)^z,\quad w_2={}^{x}({}^yz) r((xy)^z)\end{equation}
on computing these expressions as in the proof of Lemma~\ref{leml1r2}. Equality of the first factors for all $x,y,z$ is {\bf l1} and of the 2nd factors is {\bf r2} for $x,y,z$, viewed in $X^2$. Similarly for {\bf r1,l2}. 
 \end{proof}

These observations will play a role later on. Clearly, all of these conditions and {\bf lr3} hold in the case of a braided set.

\subsection{Cancellation conditions}

To proceed further we require and investigate next some
cancellation conditions. A sufficient but not necessary condition
for them is if $X\subset G(X,r)$ is an inclusion.
\begin{definition}
\label{cancellsoldef} Let $(X,r)$ be a quadratic set. We say that $r$ is \emph{2-cancellative} if for
every positive integer $k$, less than the order of $r$, the
following two condition holds:
\[
r^k(x,y)= (x,z) \Longrightarrow z=y \;\;\text{(left
2-cancellative )}
\]
\[
r^k(x,y)= (t,y) \Longrightarrow x=t \;\;\text{(right
2-cancellative)}
\]
\end{definition}
It follows from Corollary~\ref{involcancelcorol}
 that every non-degenerate involutive
quadratic map  $(X,r)$ is 2-cancellative.
\begin{proposition}
\label{cancellsolprop} Let $(X,r)$ be a set with quadratic map,
with associated monoid $S=S(X,r).$ Then:

(1)  $S$  has (respectively left, right) cancellation on monomials of length 2 {\em  if and
only if}  $r$ is (respectively left, right)  2-cancellative;

(2) Suppose that $(X,r)$ is 2-cancellative  solution.
 Then $S$ has  cancellation on monomials of length 3.
\end{proposition}
\begin{proof}
We recall that the defining relations of $S$ come from the map $r.$
We do not assume $r$ is necessarily of finite order,  so, in
general,
\[
xy=xz\; \text{in}\;\; S \Longleftrightarrow xz=r^k(xy),\;\text{for\
some\ positive\  integer}\  k,\]
\[\quad\quad\quad\quad \quad\quad\quad\quad\quad \text{or} \;\;xy=r^k(xz)\;\text{for\
some\ positive\ integer}\ k .
\]
Assume $xy=xz$ in $S$. Without loss of generality, we may assume
$xz=r^k(xy).$ But $r$ is 2-cancellative, so $z=y$

Clearly if $r$ is not left 2-cancellative, then there is an
equality $xz=r^k(xy)$ for some integer $k \geq 1$, and $z\neq y$.
This gives $xz=xy$ in $S,$ therefore $S$ is not left cancellative.
This proves the left case of (1). The proof of the right case is
analogous.

 We shall prove (2). Suppose  $(X,r)$ is a 2-cancellative solution. Let $xyz=xpq$ be an
equality in $S$. Then the monomial $xpq$, considered as an element
of $X^3,$ is in the orbit $\Ocal_{\Dcal}(xyz)$ of $xyz$ under the
action of the group $\Dcal(r)$ on $X^3$, and therefore occurs in
the YB-diagram (\ref{ybediagram}). We study the six possible
cases. In the first four  cases we follow the left vertical branch
of the diagram.
\[
\text(a) \;\;  (x,p,q)=(x,y, z) \;\text{in} \; X^{\times 3}
\;\;\Longrightarrow (p,q)=(y,z) \;\text{in} \; X\times X
\;\Longrightarrow pq=yz \;\text{in} \; S.
\]
\[
\text(b) \;\;  (x,p,q)=(x,r(y,z)) \;\text{in} \; X^{\times 3}
\;\Longrightarrow (p,q) =r(y,z) \;\text{in} \; X\times X
\;\Longrightarrow pq=yz \;\text{in} \; S.
\]
\[
\text(c) \;\;  (x,p,q)=(r(x,{}^yz), y^z) \;\text{in} \; X^3
\;\;\Longrightarrow  q=y^z \;\; \text{and}\;\; (x,p)=r(x,{}^yz)
\]
\[
\Longrightarrow^{\text{by $r$ 2-cancellative}}\; p= {}^yz
\Longrightarrow (x,p, q) = (x, {}^yz,y^z )
\]
\[\Longrightarrow (p,q)= r(y,z) \Longrightarrow pq=yz \;\;\text{in}\; S.
\]
\[
\text(d) \;\;  (x,p,q)=({}^{x}{({}^yz)},
r((x^{{}^yz}),y^z)\Longrightarrow
x={}^{x}{({}^yz)}\Longrightarrow^{\text{by $r$ 2-cancellative}}\;
x^{{}^yz}= {}^yz
\]
\[
\Longrightarrow (x,p,q)=(x, r({}^yz,y^z)= (x,
r^2(y,z))\Longrightarrow (p,q) = r^2(y,z)\Longrightarrow pq=yz
\;\;\text{in}\; S.
\]
\[
\text(e) \;\; (x,p,q)= (r(x,y), z) \;\text{in} \; X^{\times 3}
\;\Longrightarrow q = z\; \text{and}\; (x,p) = r(x,y)
\]
\[
\Longrightarrow^{\text{by $r$ 2-cancellative}} p = y\Longrightarrow
(x,p,q) = (x,y,z) \;\text{in} \; X^{\times 3}
\]
\[
\text(f) \;\; (x,p,q)= ({}^xy, r(x^y, z)) \;\text{in} \; X^{\times
3} \;\Longrightarrow {}^xy = x \Longrightarrow^{\text{by $r$
2-cancellative}}\;
\]
\[
x^y=y\Longrightarrow (x,p,q)= (x, r(y, z))\Longrightarrow (p,q)= r(y, z)\Longrightarrow pq=yz \;\text{in} \; S.
\]
We have shown that
\[
xyz=xpq \;\text{is an equality  in}\; S^3\Longrightarrow yz=pq \;
\text{in } \; S^2
\]
Assume now that $xyz=xyt$ is an equality in $S^3.$ It follows from
the previous that $yz=yt,$ in $S^2,$ therefore by the
2-cancellativeness of $S$ $z=t.$

Analogous argument verifies the right cancellation in $S^3$.
\end{proof}

The following Lemma~\ref{Zrlemma3} and Corollary
\ref{involcancelcorol} show that each  non-degenerate involutive
quadratic map is 2-cancellative:

\begin{lemma}
\label{Zrlemma3} Suppose $(X,r)$ is  non-degenerate and involutive
quadratic map (not necessarily a solution of YBE),  $x,y \in X.$
The following conditions are equivalent:
\begin{enumerate}
\item \label{1}
 $r(x,y) = (x,t),$ for some $t \in X$;
 \item
 \label{2} $r(x,y)= (s,y),$ for some $s \in X$;
\item \label{3} ${}^xy = x;$ \item \label{4} $x^y = y;$ \item
\label{5} $r(x,y) = (x,y).$
\end{enumerate}
\end{lemma}
\begin{proof}
Suppose $r(x,y) = (x,t),$ for some $t \in X$. Clearly it follows
from the standard equality $r(x,y)= ({}^xy, y^x)$ that
$(\ref{5})\Longrightarrow (\ref{1})\Longleftrightarrow (\ref{3}),$ and
$(\ref{5})\Longrightarrow (\ref{2})\Longleftrightarrow (\ref{4}).$

We will show $(\ref{3})\Longrightarrow (\ref{4}).$  Assume ${}^xy= x.$
Then the equalities
\[
{}^xy= x=^{r\; \text{involutive}}
\;{}^{{}^xy}{(x^y)}={}^x{(x^y)}.
\]
give ${}^xy={}^x{(x^y)},$ which by the non-degeneracy of $r$
implies $y = x^y.$ One analogously proves the implication $(\ref{4})\Longrightarrow
(\ref{3}).$

We have shown that each of the first four conditions implies the
remaining three and hence, clearly (\ref{5}) also holds.  This
completes the proof of the lemma.
\end{proof}
\begin{corollary}
\label{involcancelcorol} Let $(X,r)$ be  non-degenerate and
involutive quadratic set. Then $r$ is 2-cancellative and
$S=S(X,r)$ has cancellation on monomials of length 2.

Furthermore, if $(X,r)$ is a solution of YBE, then   $S$ has
cancellation on monomials of length 3.
\end{corollary}

The following example gives a non-degenerate bijective solution
$(X,r)$ which is not 2-cancellative,
\begin{example}
Let $X=\{x,y,z\},$ $\rho = (x y z)$ , be a cycle of length three
in $Sym(X).$ Define $r(x,y):=(\rho(y), x).$ Then $r:X\times
X\longrightarrow X\times X$ is  a non-degenerate bijection of
order $6$.
\[
(x, x)\longrightarrow^r (y, x)\longrightarrow^r (y,
y)\longrightarrow^r (z, y)\longrightarrow^r (z,
z)\longrightarrow^r (x, z)\longrightarrow^r (x, x),
\]
\[
(x, y)\longrightarrow^r (z, x)\longrightarrow^r (y,
z)\longrightarrow^r (x, y).
\]
It is easy to check that $r$ is a solution of YBE, (this is a
permutation solution). The two actions satisfy: $\Lcal_x\Lcal_y=\Lcal_z=(x y z);  \Rcal_x=\Rcal_y=\Rcal_z=e,$ so $r$ is
nondegenerate. Note that $r^2$ fails to be nondegenerate, since
$r^2(y,x)=(z, y),$ $r^2(y,y)=(z, z).$ Moreover, $r$ is not
2-cancellative, since $xx=yx$ is an equality in $S$. Note also
that in the group $G(X,r)$ all generators are equal: $x=y=z.$
\end{example}
\begin{lemma}\label{cancellativelemma} Let $(X,r)$ be a 2-cancellative quadratic set. 
 If  the monoid $S(X,r)$
has cancellation on monomials of length 3 then
\begin{enumerate}
\item [\text{(i)}] ${\bf l2}\ \Longleftrightarrow\ {\bf r2}\
\Longleftrightarrow\ {\bf l1}\ \&\ {\bf r1}\ \Longleftrightarrow\
(X,r)\ is\ braided\ set.$

Furthermore, if $r$ is nondegenerate and involutive then 
\item [\text{(ii)}] 
${\bf l1} \Longleftrightarrow {\bf r1}\; \Longleftrightarrow
\;(X,r)$ is a symmetric set.
\end{enumerate}
\end{lemma}
\begin{proof}
 We look again at the requirements of the YBE in the
proof recalled above. The monoid $S(X,r)$ has cancellation on
monomials of length 3. So if we look at this diagram
(\ref{ybediagram}) in the monoid $S(X,r)$ then as words in
$S(X,r)$ we already have equality $w_1=w_2$. If we assume {\bf l1,
lr3} we can cancel the first two factors and deduce {\bf r1} holding in
$S(X,r)$. But $X\subset S(X,r),$ so we can conclude {\bf r1} in
$X$. Similarly any two of ${\bf l1,lr3,r1}$ allows us to conclude
the third. It follows from Lemma~\ref{ybe} that any two is
equivalent in this case to $(X,r)$ braided. This proves (i). Also,
if we assume now that $r$ is involutive and nondegenerate 
then by Proposition~\ref{ESS1}, 
one has $\textbf{r1}\;
\Longrightarrow \; \textbf{l1, lr3}$, similarly $\textbf{l1}\;
\Longrightarrow \; \textbf{r1, lr3},$ thus we have
$\textbf{l1}\;
\Longleftrightarrow \; \textbf{r1},$ which by (i) gives  the last
part.
\end{proof}
Clearly for the proof of the lemma it is enough to   assume
$X\subset G(X,r)$ is an inclusion. It is important to note that
the condition $X\subset G(X,r)$ or an equivalent one is not empty
as shows the following example.
\begin{example}
Consider $X=\{x,y,z\}$ and $r$ on monomials:
\[ xy\to xz\to yz\to yy\to xy, \quad xx\to zz\to yx\to zy\to xx,\quad zx\to
zx.\] Clearly, $r$ is non-degenerate and as we show below obeys
{\bf l1, r1}. On the other hand it does not obey the YBE. Note
that $x=y=z$ in $G(X,r)$ by cancelling in the group, so $X$ is not
contained in $G(X,r)$. Here we give some details. It is not
difficult to see that for each pair $\xi,\eta \in X$ there is an
equality
\[
r(\xi \eta)= (\Lcal(\eta), \Rcal(\xi)), \;\text{where}\; \Lcal =(x \;z\; y)\in S_3,\quad \Rcal= (x\; z) \in S_3.
\]
In other words $r$ is \emph{a permutational map}, (see Definition
\ref{Permutationalsolution}), and in the notation of (\ref{LcalRcal}) we have $\Lcal _x=\Lcal _y=\Lcal _z=\Lcal = (x\;
z\; y),$ a cycle of length three in $S_3$, and $\Rcal _x=\Rcal _y
=\Rcal _z=\Rcal = (x \;z),$ a transposition in $S_3.$ It is clear
then that the left action satisfy
\[
{}^{a}{({}^b{\xi})}={\Lcal}^2(\xi)= {}^{{}^ab}{({}^{a^b}{\xi})},
\;\;\text{for all}\;\; \xi, a, b \in X,
\]
therefore \textbf{l1 } holds.
 Similarly for the right action one has:
\[
({\xi}^a)^b={\Rcal}^2(\xi)= \xi = ({\xi}^{{}^ab})^{a^b},
\;\text{for all}\; \xi, a, b \in X,
\]
which verifies  \textbf{r1 }. Furthermore, as a permutational map
$r$ obeys all cyclic conditions, see Definition
\ref{cyclicconditionsall}. Direct computation shows that
$r^{12}r^{23}r^{12}(xyz) \neq r^{23}r^{12}r^{23}(xyz),$ so $(X,r)$
does not obey the YBE.
\end{example}

\subsection{Cyclic conditions}

 In various cases in a nodegenerate $(X,r)$ the left and the
right actions are inverses, i.e. $\Rcal_x= \Lcal_x^{-1}$ and
$\Lcal_x = \Rcal_x^{-1}$ for all $x\in X$. For example
this is true for every square-free symmetric set $(X,r)$. It
is natural then to single out a class of non-degenerate sets
$(X,r)$ by a stronger condition \textbf{lri} defined below, and study the
the relation of this property to the other conditions on the left
(resp. right) action.

\begin{definition}
\label{lri} Let  $(X,r)$ be a quadratic set. We define the condition
\[ \textbf{lri:} \;\;\;\;\;\;\;\;\;\; ({}^xy)^x= y={}^x{(y^x)} \;\text{for all} \; x,y \in X.\]
In other words \textbf{lri} holds if and only if  $\Rcal_x=\Lcal_x^{-1}$ and $\Lcal_x = \Rcal_x^{-1}$.
\end{definition}

\begin{lemma}
\label{(yx)y}  Let $(X,r)$ be a quadratic set. Then the following are equivalent:

\begin{enumerate} \item {\bf lri} holds
\item  $ ({}^xy)^x= y$ for all $x,y \in X$
\item ${}^x(y^x)= y$ for all $x,y\in X$.
\end{enumerate}

In this case $(X,r)$ is nondegenerate.
\end{lemma}
\begin{proof}
Assume (2).  We will show first that $(X,r)$ is left
non-degenerate. Assume ${}^xy={}^xz.$ Then 
\[
y=({}^xy)^x=({}^xz)^x=z,
\]
which proves the left non-degeneracy. Assume now that
$y^x=z^x.$ By the left non-degeneracy, there exist unique $s,t \in X,$
with $y = {}^xs, z ={}^xt.$ Then
\[
y^x=z^x \Longrightarrow \; s={({}^xs)}^x={({}^xt)
}^x=t \; \Longrightarrow \; y=z,
\]
so $(X,r)$ is  right non-degenerate, and therefore non-degenerate. 

Next we show that (2) $\Rightarrow$ (3). By the nondegeneracy of $(X,r)$ for any $x,y \in X$
there exists unique
$z \in X,$ with $y={}^xz.$ Then 
\[{}^x{(y^x)}={}^x{(({}^xz)^x)}=^{\text{by
assumption}}\;{}^xz=y.\]
The proof of the implication in the opposite direction
is analogous. Thus either part implies {\bf lri}.
\end{proof} 

We will show later that under some  cyclicity
restrictions, the involutiveness of $r$   is equivalent to
condition \textbf{lri}. The following lemma is straight forward.
\begin{lemma}
\label{LequivalentR} Assume $(X,r)$ satisfies {\bf
lri}. Then
${\bf l1}\Longleftrightarrow {\bf r1}$, and
${\bf l2}\Longleftrightarrow {\bf r2}.$
\end{lemma}
More generally, clearly, \textbf{lri} implies that whatever property is
satisfied by the left action, an analogous property is
valid for the right action and vice versa. In particular, this is
valid for the left and right `cyclic conditions'. Such conditions
were discovered, see \cite{T1,T2,T06} when the
 first author studied binomial rings with skew polynomial relation and square-free solutions of YBE.
It is interesting to know that the proofs of the good
algebraic and homological properties of these algebras and monoids
use in explicit or implicit form the existence of the full
cyclic condition in the form below (Definition~\ref{cyclicconditionsall}). This includes
the properties of being Noetherian, Gorenstein, therefore
Artin-Shelter regular, and "producing" solutions of YBE, see
\cite{T1,T2,T3,T06,T07,TM,TJO,JO}. Compared with these
works, here we do not assume that $X$ is finite, and initially the
only restriction on the map $r$ we impose is "$r$ is
non-degenerate". We recall first the notion of "cyclic conditions" in
terms of the left and right actions and study the implication of
the cyclic conditions on the properties of the actions, in
particular  how are they related with involutiveness of $r$ and and
condition
\textbf{lri}.
\begin{definition}
\label{cyclicconditionsall} Let $(X,r)$ be a quadratic
set. We define the conditions
\[\begin{array}{lclc}
 {\rm\bf cl1:}\quad&  {}^{y^x}x= {}^yx \quad\text{for all}\; x,y \in X;
 \quad&{\rm\bf cr1:}\quad &x^{{}^xy}= x^y, \quad\text{for all}\; x,y \in
X;\\
 {\rm\bf cl2:}\quad &{}^{{}^xy}x= {}^yx, \quad\text{for all}\; x,y \in X; \quad & {\rm\bf
cr2:}\quad  &x^{y^x}= x^y \quad\text{for all}\; x,y \in X.
\end{array}\]
We say $(X,r)$ is {\em weak cyclic} if {\bf cl1,cr1}
hold and is {\em cyclic} if all four of the above hold.
\end{definition}

One can also define left-cyclic as {\bf cl1,cl2} and similarly right-cyclic,  which is related to cycle
sets to be considered later.

\begin{example}
\label{Permutationalsolution} (\emph{Permutational solution},
Lyubashenko, \cite{D}) Let $X$ be nonempty set, let
$f, g$ be bijective maps $X\longrightarrow X$, and let  $r(x,y)=
(g(y),f(x))$. We call $r$  \emph{a permutational map}. Then a)
$(X,r)$ is braided if and only if $fg=gf.$ Clearly, in this case
${}^xy={}^zy = g(y),$ and $y^x=y^z = f(y)$ for all $x,z, y \in X$,
so $\Lcal_x=g,$ and $\Rcal_x=f$ for all $x \in X,$
hence $(X,r)$ is cyclic. b) $(X,r)$ is symmetric if and
only if $f=g^{-1}.$ Note that if $n \geq 2,$ and $f \neq \id_X$, or
$g \neq \id_X$, the solution $(X,r)$ is never square-free.
\end{example}
The following examples come from Lemma~\ref{equivalenceofallcc}.

\begin{example}\ 
\begin{enumerate}
\item Every square-free non-degenerate braided set
$(X,r)$ is weak cyclic.
 \item Every square-free left cycle set (see Definition~\ref{csl&csr}) satisfies \textbf{cl2}.
 \item Every square-free right cycle set (see Def.~\ref{csl&csr})  satisfies
\textbf{cr2}.
\end{enumerate}
\end{example}
\begin{lemma}
\label{equivalenceofallcc} Assume $(X,r)$ is  with
{\bf lri}. Then the following conditions are equivalent
\begin{enumerate}
\item $(X,r)$ satisfies {\bf cl1}; \item $(X,r)$
satisfies {\bf cl2}; \item $(X,r)$ satisfies {\bf cr1}; 
\item $(X,r)$ satisfies {\bf cr2}; \item $(X,r)$ is cyclic.
\end{enumerate}
\end{lemma}
\begin{proof} Note first that by  Lemma~\ref{(yx)y}   $(X,r)$ is nondegenerate.
 Suppose first that {\bf cl1} holds, so   $ {}^{y^x}x={}^yx $ for all $x,y \in X$.
 In this equality we set $z=y^x,$ and (by \textbf{lri}) $y={}^xz,$
 and obtain by the nondegeneracy  $ {}^zx={}^{{}^xz}x$ for all $x,z \in X$, i.e. {\bf cl2}.
 The implication
\[ \textbf{cl2} \Longrightarrow \textbf{cl1}\]
is analogous. It follows then that \textbf{cl1} and \textbf{cl2}
are equivalent. But under the assumption of  \textbf{lri} one has
\[
 \textbf{cl1}\Longleftrightarrow \textbf{cr1} \; \;\text{and} \; \;  \textbf{cl2}\Longleftrightarrow \textbf{cr2}.
 \]
We have shown that the conditions \textbf{cl1}, \textbf{cl2},
\textbf{cr1} \textbf{cr2}, are equivalent. Clearly, each of them
implies the remaining three conditions, and therefore (see 
Definition~\ref{cyclicconditionsall}) $(X,r)$ is cyclic.  The
inverse implications follow straightforwardly from Definition~\ref{cyclicconditionsall}.
\end{proof}
\begin{proposition}
\label{lri&INVOLUTIVE}
Let  $(X,r)$ be a quadratic set. Then any two of the following conditions imply the remaining third
condition.
\begin{enumerate}
\item \label{cond1} $(X,r)$ is involutive
\item \label{cond2} $(X,r)$ is nondegenerate and cyclic.
\item \label{cond3} {\bf lri} holds.
\end{enumerate}
\end{proposition}
\begin{proof} 
(\ref{cond1}), (\ref{cond2}) $\Longrightarrow$ (\ref{cond3}). 
By assumption $r$ is involutive, hence $({}^yx)^{y^x}
= x$, which together with \textbf{cl1} in the form ${}^yx={}^{y^x}x,$  implies
\[
({}^{y^x}x)^{y^x} = x.
\]
By the non-degeneracy of $r$ it follows that given
$x\in X$, every $z \in X$ can be presented as $z=y^x,$ for appropriate
uniquely determined $y$. We have shown
\[
({}^zx)^z = x \;\;\text{for every} \;\; x,z \in X.
\]
It follows from Lemma~\ref{(yx)y} that ${}^z(x^z)= x$
is also in force for all $x,z \in X,$ so {\bf lri} is in force.

(\ref{cond2}), (\ref{cond3}) $\Longrightarrow$ (\ref{cond1}).
For the involutiveness of $r$ it will be enough to show $({}^yx)^{y^x}=x$ for
all $x,y \in X.$  We set
\begin{equation}
\label{z} z={}^xy, \;\;  \text{so}\;\;  y=^{\textbf{lri}}z^x.
\end{equation}
The equalities
\[ x=^\textbf{lri}({}^yx)^y  =^\textbf{cl2}({}^{{}^xy}x)^y  =({}^zx)^{z^x}\]
by (\ref{z}), imply that $({}^zx)^{z^x}=x$, which, by the non-degeneracy of
 $r$ is valid for all $z, x \in X.$ Therefore $r$ is involutive
and the implication (\ref{cond2}), (\ref{cond3}) $\Longrightarrow$ (\ref{cond1}) is verified.

(\ref{cond1}), (\ref{cond3}) $\Longrightarrow$ (\ref{cond2}). 
By  Lemma \ref{(yx)y}  $(X,r)$ is nondegenerate. Let $x,y \in X.$ We will show {\bf cl1}. Consider the
equalities
\[
({}^{y^x}x)^{y^x}=^{{\bf lri}}\; x=^{\text{$r$
involutive}}\; ({}^yx)^{y^x}.
\] 
Thus by the nondegeneracy, ${}^{y^x}x={}^yx,$ which
verifies {\bf cl1}. The hypothesis of Lemma \ref{equivalenceofallcc} is satisfied, therefore
{\bf cl1} implies $(X,r)$ cyclic. The proposition has been proved.
\end{proof}
\begin{remark}
In Proposition~\ref{lri&INVOLUTIVE} one can replace (2) by the weaker condition\\
\noindent (2')\quad $(X,r)$ nondegenerate and  {\bf cl1} holds.
\end{remark}

We give now the definition of \emph{(left) cycle set}.
\begin{remark}
\label{remarkcycleset} The notion of cycle set was
introduced by Rump, see \cite{Rump} and was used in the proof of the
decomposition theorem. In his definition Rump assumes
that the left and the right actions on $X$ are inverses (or in
our language, condition \textbf{lri} holds), and that $r$
is involutive. We keep the name ``cycle set" but we
suggest a bit more general definition here. We do not assume that $r$ is
involutive, neither that \textbf{lri} holds. Therefore we have to
distinguish left and right cycle sets. Furthermore, Corollary
\ref{squarefreelri&rinvol}  shows that for a square-free left
cycle set conditions $r$ is involutive and
\textbf{lri} are equivalent.
\end{remark}

\begin{definition}
\label{csl&csr} Let  $(X,r)$ be a quadratic set.

a) \cite{Rump} $(X,r)$ is called \emph{a left cycle
set} if  
\[ \textbf{csl} \; \; \;\; \; \; {}{}^{({}^yt)}{({}^yx)}={}{}^{({}^ty)}{({}^tx)} \; \;
\text{for all}\; x,y,t \in X\]

b) Analogously we define \emph{a right cycle set} by
the condition
\[ \textbf{csr} \; \; \;\; \; \; (x^y)^{t^y}=(x^t)^{y^t}
\; \; \text{for all}\; x,y,t \in X.\]
\end{definition}

\begin{proposition}\label{csl&l1} Let $(X,r)$ be
involutive with {\bf lri}. Then 
\begin{enumerate}
\item [\text{(i)}] $(X,r)$ is non-degenerate and cyclic. 
\item [\text{(ii)}] 
${\bf csl} \Longleftrightarrow {\bf l1}$.
\item [\text{(iii)}] ${\bf csl}$ $\Rightarrow$ $({}^xz)^y={}^{x^y}(z^{y^x})$ and $ {}^x(z^y)=({}^{{}^yx}z)^{{}^xy}$ for all $x,y,z \in X.$
\item [\text{(iv)}] Suppose that $(X,r)$ is
 2-cancellative and the monoid $S(X,r)$
has cancellation on monomials of length 3 then
${\bf csl} \Longleftrightarrow {\bf l1}\; \Longleftrightarrow
\;(X,r)$ is a symmetric set.
\end{enumerate}
\end{proposition}
\begin{proof}
Proposition \ref{lri&INVOLUTIVE} implies (i). 
To prove (ii) we first show  {\bf l1} $\; \Longrightarrow {\bf
csl}$.  Let $x,y,t \in X.$ Consider the equalities
\[
{}{}^{({}^yt)}{({}^yx)}=^{{\bf l1}} {}^{{}^{{}^yt}y}{({}^{({}^yt)^y}x)}=^{{\bf cl2, lri}} \; {}{}^{({}^ty)}{({}^tx)}.
\]
This verifies {\bf csl}. Next we assume  {\bf csl} holds in $(X,r)$ and shall
verify {\bf l1}. Let $x,y,z \in X.$  Set $t=z^y.$ Then, by {\bf lri} $z={}^yt,$ and the
following equalities verify {\bf l1}
\[
{}^{z}{({}^yx)}={}^{{}^yt}{({}^yx)}=^{{\bf csl}} \;{}{}^{({}^ty)}{({}^tx)}={}^{({}^{z^y}y)}{({}^{z^y}x)}=^{{\bf
cl1}} {}^{({}^zy)}{({}^{z^y}x)}.
\]
We have shown ${\bf csl} \Longleftrightarrow {\bf l1}$.

We shall prove the first equality displayed in (iii). 
Assume  \textbf{l1}. Note that \textbf{lri} and $(X,r)$ cyclic imply the
following equalities in $S(X,r):$
\begin{equation}
\label{CS3} x.{y^x}={}^x{(y^x)}.x^{y^x}=y.x^y,
\end{equation}
and
\begin{equation}
\label{CS4} {}^yx.y={}^{{}^yx}y.({}^yx)^y={}^xy.x
\end{equation}
 Consider the equalities:
\[
{}^xz=^{\bf lri} {}^x{({}^{y^x}{(z^{y^x})})}=^{\bf l1}{}^{(xy^x)}{(z^{y^x})} 
={}^{(yx^y)}{(z^{y^x})}={}^{y}{({}^{x^y}{(z^{y^x})})}\]
using \textbf{l1} and (\ref{CS3}) followed by {\bf l1}. We apply $\bullet^y$ to both sides to get
\[
({}^xz)^y={}^{x^y}{(z^{y^x})}.
\]
For the second equality displayed in (iii), we
consider
\[
z^y=^{\bf lri}(({}^{{}^yx}z)^{{}^yx})^y =^{\bf r1}({}^{{}^yx}z)^{({}^yx.y)}=(({}^{{}^yx}z)^{{}^xy.x}=(({}^{{}^yx}z)^{{}^xy})^x
\]
using \textbf{r1} and  (\ref{CS4}) followed by {\bf r1}. Then by 
 \textbf{lri}, we obtain
\[
{}^x{(z^y)}=({}^{{}^yx}z)^{{}^xy}.
\]
This verifies (iii).

Assume now $(X,r)$ is 2-cancellative and the monoid $S(X,r)$
has cancellation on monomials of length 3. Then by  Lemma~\ref{cancellativelemma} one has 
\begin{equation}
\label{l1r1sym1} {\bf l1, r1} \; \Longleftrightarrow \;(X,r) \;\; \text{is a symmetric set}.
\end{equation} 
Clearly, {\bf lri} implies {\bf l1}
$\Longleftrightarrow$ {\bf r1}, which together with (\ref{l1r1sym1}) and (ii) gives
\[
{\bf csl} \Longleftrightarrow {\bf l1} \; \Longleftrightarrow \;(X,r) \;\; \text{is a symmetric set}.
\]
The proof of the proposition is now complete. 
\end{proof}

We now look at the more general non-degenerate case.

\begin{lemma}
\label{squarefreecl1cl2} Assume $(X,r)$ is non-degenerate and
square-free. (We do not assume involutiveness.) Then
\[\begin{array}{lclc}
 {\rm\bf l1}\ \Longrightarrow\ {\rm\bf cl1};\quad&{\rm\bf r1}\ \Longrightarrow\ {\rm\bf
 cr1};\quad&{\rm\bf lr3}\ \Longrightarrow\ {\rm\bf cl1,cr1};\\
{\rm\bf csl}\  \Longrightarrow\ {\rm\bf cl2}; \quad& {\rm\bf csr}
\ \Longrightarrow\ {\rm\bf cr2}.
\end{array}\]
\end{lemma}
\begin{proof}
We prove first ${\rm\bf l1} \Longrightarrow\ {\rm\bf cl1}.$
 Suppose $(X,r)$ satisfies {\bf l1}. Consider the equalities
\begin{equation}
\label{e4} {}^yx={}^{y}{({}^xx)} =^{\bf l1}
{}^{{}^yx}{({}^{y^x}x)}.
\end{equation}
By Lemma~\ref{someproperties}.\ref{SFnondeg}, ${}^tx=t$ implies
$x=t$, which together with (\ref{e4}) gives  ${}^{y^x}x = {}^yx.$
We have shown $\textbf{l1} \Longrightarrow \textbf{cl1}$. The
implication ${\rm\bf r1} \Longrightarrow\ {\rm\bf cr1}$  is analogous.

Suppose {\bf lr3} holds. We set $z=y$ in the  equality
\begin{equation}
 \label{lr3again}
 {({}^xy)}^{({}^{x^y}{(z)})}={}^{(x^{{}^yz})}{(y^z)}, \quad \text{for all} \; x,y,z, \in X,
 \end{equation}
and obtain:
\begin{equation}
 \label{wcc1}
 {({}^xy)}^{({}^{x^y}{(y)})} ={}^{(x^y)}{y}.
 \end{equation}
 By hypothesis $(X,r)$ is square-free, therefore (\ref{wcc1}) and Lemma~\ref{someproperties}
imply
 \[
 {}^xy
 ={}^{(x^y)}{y},
 \]
 which proves \textbf{cl1}. If we set $x=y$ in (\ref{lr3again}), similar argument proves
\textbf{cr1}. Therefore \textbf{lr3} implies that $(X,r)$ is weakly
cyclic as stated.

Assume \textbf{csl} is satisfied. So, for all $x,y,t\in X$ one
has $ {}^{{}^yt}{({}^yx)}= {}^{{}^ty}{({}^tx)}$ in
which  we substitute in $y=x$  and obtain
\[
{}^{{}^xt}{({}^xx)}= {}^{{}^tx}{({}^tx)}.
\]
This, since $r$ is square-free, yields ${}^{({}^xt)}x={}^tx$,
which verifies $\textbf{csl} \Longrightarrow\textbf{cl2}$. The
proof of $\textbf{csr} \Longrightarrow \textbf{cr2}$ is analogous.
\end{proof}

Lemma~\ref{squarefreecl1cl2} and Proposition~\ref{lri&INVOLUTIVE}
imply the following corollary.
\begin{corollary}
\label{squarefreelri&rinvol} Suppose $(X,r)$ is
non-degenerate and square free. If any of the conditions {\bf l1}, {\bf r1},
{\bf lr3}, {\bf csl}, {\bf csr} hold then
\[
r \;\; \text{is involutive} \Longleftrightarrow {\bf lri}
\]
In this case $(X,r)$ is cyclic.
\end{corollary}
\begin{lemma}
\label{squarefreelr3impliesYBE} Suppose $(X,r)$ is
non-degenerate, square-free and involutive (i.e. a quantum binomial set). Then
condition {\bf lr3} implies that $(X,r)$ is symmetric set with {\bf lri}.
\end{lemma}
\begin{proof}
Note first that by Corollary~\ref{squarefreelri&rinvol}
\textbf{lri} is in force.  By Lemma~\ref{squarefreecl1cl2}
\textbf{lr3} implies the weak cyclic conditions
\textbf{cl1, cl2}.  We shall prove the implication
 \[
\textbf{lr3} \Longrightarrow \textbf{l1}.
\]
Let $x,u, t \in X.$ We have to show 
\begin{equation}
\label{equ0} {}^x{({}^ut)}= {}^{{}^xu}{({}^{x^u}t)}.
 \end{equation}
By assumption the equality (\ref{lr3again}) holds. We
apply the left action ${}^{{}^{x^y}z}{\bullet}$ , to each side of
(\ref{lr3again}), and by \textbf{lri} we obtain
\begin{equation}
\label{equ1} 
{}^xy={}^{({}^{x^y}{z})}{({}^{(x^{{}^yz})}{(y^z)})}
\end{equation}
for all  $x,y,z \in X$. First we set
\begin{equation}
\label{equ2} z = u^t,\quad y = {}^zt, 
\end{equation}
hence by \textbf{lri},
\begin{equation}
\label{equ3} y^z=t,\quad {}^tz=u; 
\end{equation}
\begin{equation}
\label{equ4}
{}^yz={}^{y^z}z=^{\textbf{cl1}}\;{}^tz=u, \quad \; x^{{}^yz}=x^u.
\end{equation}
Next
\begin{equation}
\label{equ5}
{}^zt={}^{u^t}t =^{\textbf{cl1}}\; {}^ut; \quad
\Rightarrow\quad\; x^y=x^{{}^zt}=x^{{}^ut},
\end{equation}
which together with (\ref{equ2}) and \textbf{lr3} yields
\begin{equation}
\label{equ6}
{}^{x^y}z={}^{x^{{}^ut}}{(u^t)}={({}^xu)^{{}^{x^u}t}}.
\end{equation}
We now use (\ref{equ2}), (\ref{equ3}), and
(\ref{equ5}),  to obtain for the left-handside of (\ref{equ1})
\begin{equation}
\label{equ7} {}^xy={}^x({}^zt)={}^x({}^ut).
\end{equation}
For right-handside of (\ref{equ1}),  we use (\ref{equ3}),(\ref{equ4}), and (\ref{equ6}) to find
\begin{equation}
\label{equ8}
{}^{({}^{x^y}{z})}{({}^{(x^{{}^yz})}{(y^z)})}={}^{{({}^xu)^{{}^{x^u}t}}}{({}^{x^u}t)}=^{\textbf{cl1}}
\; {}^{({}^xu)}{({}^{x^u}t)}.
\end{equation}
Now (\ref{equ7}), (\ref{equ8}), and (\ref{equ1}) imply (\ref{equ0}).
\end{proof}

\begin{corollary}
\label{csl&symmetrics} Let $(X, r)$ be a
non-degenerate, square-free involutive set. Then
\[
(X,r) \;\;\text{is symmetric} \;\Longleftrightarrow
{\bf csl}.
\]
In this case $(X,r)$ is cyclic and satisfies  {\bf lri}.
\end{corollary}
\begin{proof}
Assume \textbf{csl} holds. Note first that both \textbf{lri} and
the  conditions for $(X,r)$ cyclic are satisfied. Indeed,  by
hypothesis $r$ is involutive, which  by Corollary~\ref{squarefreelri&rinvol}  
implies  \textbf{lri}. It follows from
Lemma~\ref{squarefreecl1cl2} that \textbf{cl2} is satisfied, so
\textbf{lri} and Lemma~\ref{equivalenceofallcc} imply
all cyclic conditions for $(X,r)$ cyclic. Now the hypothesis of 
Proposition~\ref{csl&l1} is satisfied, therefore $\textbf{csl}
\Longrightarrow \textbf{l1}$. Clearly \textbf{l1} and \textbf{lri}
imply \textbf{r1} (see Remark~\ref{LequivalentR}). It
follows by Proposition~\ref{ESS1}  that $(X,r)$ is braided and
therefore a symmetric set.

Assume next $(X,r)$ is symmetric set. Clearly then
\textbf{l1} holds, and since $(X,r)$ is square-free and
involutive,  Corollary~\ref{squarefreelri&rinvol}  implies \textbf{lri}.
Furthermore, by Lemma~\ref{squarefreecl1cl2} \textbf{l1} implies
\textbf{cl1}, and therefore by  \textbf{lri} all the conditions for
$(X,r)$ cyclic are satisfied. We use again Proposition~\ref{csl&l1} to deduce \textbf{csl}.
\end{proof}
\begin{remark}
Note that in \cite{Rump} it is shown that under the
assumptions that  $(X,r)$ is square-free, non-degenerate,
involutive and satisfies \textbf{lri},   then $(X,r)$ is braided if
and only if it is a cycle set. Here we show that a weaker
hypothesis is enough: Assuming that $(X,r)$  is square-free,
non-degenerate, and involutive, we show that condition \textbf{csl} is
equivalent to $(X,r)$ braided, and each of them implies
\textbf{lri}.
\end{remark}
\goodbreak

\begin{theorem}
Suppose $(X,r)$ is a quantum binomial set (i.e. non-degenerate, involutive and square-free). Then the
following conditions are equivalent:
\begin{enumerate}
\item \label{YBESquarefree} $(X,r)$ is a set-theoretic
solution of the Yang-Baxter equation. \item \label{YBEl1} $(X,r)$
satisfies {\bf l1}. \item \label{YBEl2} $(X,r)$ satisfies {\bf l2}.
\item \label{YBEr1} $(X,r)$ satisfies {\bf r1}. 
\item \label{YBEr2} $(X,r)$ satisfies {\bf r2}.
 \item\label{YBElr3}
$(X,r)$ satisfies {\bf lr3} \item \label{cycleset} $(X,r)$
satisfies {\bf csl}.
\end{enumerate}
In this case $(X,r)$ is cyclic and satisfies  {\bf lri}.
\end{theorem}
\begin{proof}
By Lemma~\ref{ybe} and Lemma~\ref{leml1r2}
\[(X,r) \; \text{symmetric} \Longrightarrow \textbf{l1},\textbf{l2}, \textbf{r1},
\textbf{r2},\textbf{lr3}
\]
By Proposition~\ref{ESS1} each of the conditions
\textbf{l1}, \textbf{r1} implies "$(X,r)$ is symmetric". Note that \textbf{l2}
is just \textbf{l1} and \textbf{lr3},  \textbf{r2} is just
\textbf{r1} and \textbf{lr3}. Lemma~\ref{squarefreelr3impliesYBE}
gives $\textbf{lr3}\Longrightarrow (X,r)$ is symmetric. This
verifies the equivalence of the first six conditions. Finally,
Corollary~\ref{csl&symmetrics} implies the equivalence
$\textbf{csl}\Longleftrightarrow (X,r)$ is symmetric.
\end{proof}

 \section{Matched pairs of monoids and braided monoids}

In this section we shall study  the `exponentiation' of a braided set $(X,r)$ to a matched
pair $(S,S)$ of monoids with property {\bf M3}. We shall also study the converse.

\subsection{Strong matched pairs and monoid factorisation}

The notion of matched pair of groups in relation to group
factorisation has a classical origin. For finite groups it was
used in the 1960s by Kac and Paljutkin in the construction (in
some form) of certain Hopf algebras \cite{Kac}. More recently,
such  `bismash product' or `bicrossproduct' Hopf algebras were
rediscovered respectively by Takeuchi \cite{Tak1} and Majid
\cite{Ma:phy}. The latter also extended the theory to Lie group
matched pairs \cite{Ma:mat}  and  to general Hopf algebra matched
pairs \cite{Ma:phy}. By now there have been  many works on matched
pairs in different contexts and we refer to the text
\cite{Ma:book} and references therein. In particular, this notion
was used by Lu, Yan and Zhu  to study the set-theoretic solution
of YBE and the associated `braided group', see \cite{Lu} and the
excellent review \cite{Takeuchi}.

In fact we will  need the notion of  a matched pair of monoids \cite{Ma:phy},
to which we shall add some refinements that disappear in the group case.

 \begin{definition}\label{MLaxioms} $(S,T)$ is a matched pair of monoids if $T$
 acts from the left on $S$ by ${}^{(\ )}\bullet$ and $S$ acts $T$
 from the right by $\bullet^{(\ )}$ and these two actions obey
\[\begin{array}{lclc}
{\bf ML0}:\quad & {}^a1=1,\quad  {}^1u=u;\quad &{\bf MR0:} \quad &1^u=1,\quad a^1=a \\
 {\rm\bf ML1:}\quad& {}^{(ab)}u={}^a{({}^bu)},\quad& {\rm\bf MR1:}\quad  & a^{(uv)}=(a^u)^v \\
{\rm\bf ML2:}\quad & {}^a{(u.v)}=({}^au)({}^{a^u}v),\quad &{\rm\bf
MR2:}\quad & (a.b)^u=(a^{{}^bu})(b^u),
\end{array}\]
for all $a, b\in T, u, v \in S$.
\end{definition}

The following
proposition is well-known, though more often for groups rather
than monoids. For
 completeness only, we recall briefly  the proof as well.

\begin{proposition}  A matched pair $(S,T)$ of monoids implies a monoid $S\bowtie T$
(called the {\em double cross product}) built on $S\times T$ with product and unit
\[ (u,a)(v,b)=(u.{}^av,a^v.b),\quad 1=(1,1),\quad \forall u,v\in S,\ a,b\in T\]
and containing $S,T$ as submonoids. Conversely, if   there exists
a monoid $R$ factorising into $S,T$ in the sense that the product
$\mu:S\times T\to R$ is bijective then $(S,T)$ are a matched pair
and $R\cong S\bowtie T$ by this identification $\mu$.
\end{proposition}
\begin{proof} A direct proof immediate, see for example \cite{Ma:phy} in the monoid case:
 \[
W_1=[(u,a).(v,b)].(w,c)=(u.({}^{a}{v}),
({a}^{v}).b).(w,c)=(a({}^{a}{v})({}^{(({a}^{v})b)}{w}),
((({a}^{v})b)^{w}) c )
\]
\[
=(u[{}^{a}{(v({}^{b}{w}))}],[(({a}^{v})b)^{w}] c )
\]
\[
W_2=(u,a).[(v,b).(w,c)]= (u,a).(v({}^{b}{w}),(b^{w})c)(u.({}^{a}{(v({}^{b}{w}))}), {a}^{(v({}^{b}{w})}.(b^{w})c))
\]
\[
=(u[{}^{a}{(v({}^{b}{w}))}], [(({a}^{v})b)^{w}]c)
\]
Clearly, $W_1=W_2,$ on using the matched pair properties, hence
$S\bowtie T$ is a monoid. It is also clear from the construction
that the converse is true: for $S\bowtie T$ with product in the
form stated  to be a monoid we need the matched pair conditions.
For example,
\[
(u,a). (1,1) = (u.{}^a1, a^1.1) = (u,a)
\]
requires $u.{}^a1=u$ and $a^1=a,$ for all $u \in S, a\in T.$ The
first equality implies ${}^a1=1$, similarly the other cases of
{\bf ML0, MR0}. Likewise, consider the equalities:
\[
W_1=[(1,a)(u,1)](v,1)= ({}^a{u}, a^{u}).(v, 1)=({}^a{u}.({}^{a^{u}}{v}), (a^{u})^{v})
\]
\[
W_2=(1,a)[(u,1)(v,1)]= (1,a)(u.{}^1v,1^{v})=(1,a)(uv,1)=({}^{a}{(uv)},a^{(uv)}).
\]
To be associative we need  $W_1=W_2,$ therefore {\bf MR1, ML2}.
Analogously,  we obtain the other requirements for a matched pair
from $S\bowtie T$ a monoid.

Now suppose $R$  factors into $S,T$. Consider $a\in T,u\in S$ and
$au\in R$. By the bijectivity it must be the product of some
unique elements $u'a'$ for $u'\in S$ and $v'\in T$. We define
${}^{(\ )}\bullet:T\times S\to S$ and $\bullet^{(\ )}:T\times S\to
T$ by $u'={}^au$ and $a'=a^u$. It is then easy to see that these
are actions and form a matched pair. Indeed, the product map
allows us to identify $R\cong S\times T$ by $u.a=(u,a)$.   In this
case associativity in $R$ implies
\[
(u.a).(v.b)=u(a.v)b =u.({}^{a}{v}.{a}^{v})b= (u.{}^{a}{v}).({a}^{v}.b)
\]
 i.e. the product of $R$ has the double cross product form when referred to $S\times T$ for the maps ${}^{(\ )}\bullet, \bullet^{(\ )}$ defined as above.  \end{proof}

 Next we introduce the following natural notion in this context. It is automatically satisfied in the group case but is useful in the monoid case:
 \begin{definition} A {\em strong monoid factorisation} is a factorisation  in submonoids
 $S,T$ as above such that $R$ also factorises into $T,S$.
 We say that a matched pair is {\em strong} if it corresponds to a strong factorisation.
 \end{definition}

In this case we have two bijections
\[ \mu_1:S\times T\to R,\quad \mu_2:T\times S\to R\]
and hence an invertible map
 \[ r_{T,S}=\mu_1^{-1}\mu_2:T\times S\to S\times T,\quad  r_{T,S}(a,u)=\mu_1^{-1}({}^au.a^u)=({}^au,a^u).\]
We also have two double crossed products $S\bowtie T$ and
$T\bowtie S$ and two underlying matched pairs. If $(S,T)$ is a strong matched pair we shall denote the actions for the accompanying matched pair $(T,S)$ by different
 notations $\la:S\times T\to T$ and $\ra:S\times T\to S$ to keep them distinct from the previous ones. They are defined by 
 \[ r_{T,S}^{-1}(u,a)=(u\la a,u\ra a)\]
and are characterised with respect to the original actions by
 \begin{equation}\label{lara} {}^au\la a^u=a,\quad {}^a u\ra a^u=u,\quad {}^{u\la a}(u\ra a)=u,\quad (u\la a)^{u\ra a}=a,\quad \forall u\in S, a\in T.\end{equation}
 
\begin{figure}
\[ \includegraphics{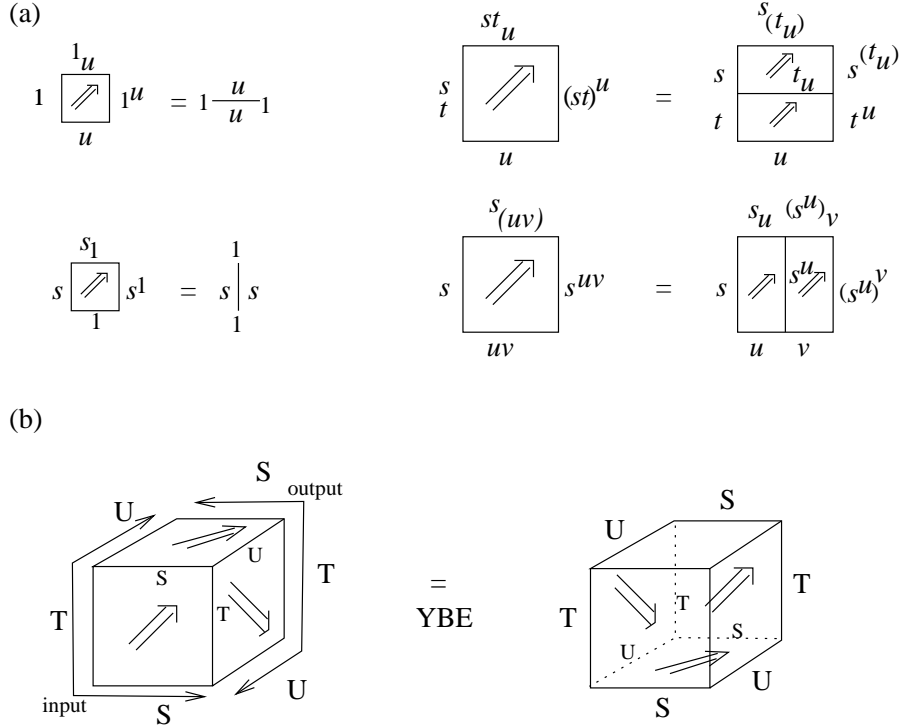}\]
\caption{Notation (a) and (b) subdivision property   encoding the
axioms of a matched pair with $\Rightarrow$ the map $r_{T,S}$, and
(c) Yang-Baxter equation as surface transport around a cube}
\end{figure}

We note also that the axioms Definition~\ref{MLaxioms} of a matched
pair $(S,T)$ have the nice interpretation and representation as a
subdivision property, see \cite{Mac:dou, Ma:book}. Thus is
recalled in Figure~1 where a box is labelled on the left and lower
edge by $T,S$ and the other two edges are determined by the
actions. This operation $\Rightarrow$ is exactly the map $r_{T,S}$
above. If one writes  out the subdivision property as a
composition of maps, it says
\[ r_{T,S}(ab,u)=(\id\times \cdot)r_{12}r_{23}(a,b,u),\quad r_{T,S}(a,uv)=(\cdot\times\id)r_{23}r_{12}(a,u,v)\]
where $r=r_{T,S}$ and the numerical suffices denote which factors
it acts on. This is just the subdivision property written out
under a different notation, but is suggestive of the axioms of a
coquasitriangular structure in the case when $T=S$, a point of
view used in \cite{Lu}.   Finally, returning to Figure~1, if the
matched pair is strong it means precisely that $\Rightarrow$ is
reversible. The reversed map evidently obeys the same subdivision
property and therefore corresponds exactly to another matched
pair, namely $(T,S)$. This proves that a matched pair is strong in the sense $r_{T,S}$ invertible if an only if $(T,S)$ is another matched pair with actions related by (\ref{lara}).

 We also define a braided monoid analogously to the term `braided group' in the
 sense of   \cite{Takeuchi}, \cite{Lu}. (This term should not be confused with
 its use to describe Hopf algebras in braided categories in another context).

 \begin{definition}
 \label{braidedmonoiddef}
 An \emph{{\bf M3}-monoid} is a monoid $S$ forming part of a matched pair $(S,S)$
 for which the actions are such that 
 \[ {\rm\bf M3}:\quad {}^uvu^v=uv\]
holds in $S$ for all $u,v\in S$. We define the \emph{associated map}   $r_S: S\times S\to S\times S$  by  $r_S(u,v)= ({}^uv,u^v).$
 A \emph{braided monoid}  is an {\bf M3}-monoid $S$ where
 $r_S$ is bijective and obeys the YBE.
 \end{definition}
 
Similarly a strong {\bf M3}-monoid is one where $(S,S)$ is a strong matched pair and this clearly happens precisely when $r_S$ is invertible. 

\begin{proposition}\label{M3cancellative} If an {\bf M3}-monoid $S$ has left cancellation then $r_S$ obeys the YBE. Hence a strong {\bf M3}-monoid with left cancellation is necessarily a braided monoid. 
\end{proposition}\begin{proof} Clearly {\bf l1,r1} on the set of $S$ hold in view of {\bf ML1,MR1,M3}. Next we  consider the identities
\[ {}^u(vw)=^{\bf M3}{}^u({}^vw v^w)=^{\bf ML2,ML1}{}^{uv}w.{}^{(u^{{}^vw})}(v^w)\]
\[ {}^uv.{}^{u^v}w=^{\bf M3,ML1}{}^{{}^uv u^v}w.({}^uv)^{({}^{u^v}w)}=^{\bf M3}{}^{uv}w.({}^uv)^{({}^{u^v}w)}.\]
The two left hand sides are equal by {\bf ML2} and hence if we can cancel on the left, the second factors on the right are equal, which is {\bf lr3} on the set $S$. We then apply  Lemma~\ref{ybe} to the set $S$. \end{proof}

This covers the group case where  $r_S$ does automatically obey
the YBE and have an inverse. Also note that it is perfectly
equivalent to think of an {\bf M3}-monoid as a pair $(S,r_S)$
where $S$ is a monoid and $r_S$ is a set map such that ${}^uv$ and
$u^v$ defined by $r_S$ via the formula above makes $(S,S)$ a
matched pair and such that the product of $S$ is $r_S$-commutative
in the sense $\mu r_S(u,v)=uv$, where $\mu$ is the product in $S$.

There is also a `surface transport' description of the YBE. Here, it was already observed in \cite{Ma:book} that the subdivision property in Figure~1(b) implies a notion of surface transport defined by the matched pair.   Fitting in with this now is a pictorial description of the generalised Yang-Baxter equation shown in part (c) of the figure, which holds in any category for a collection of such maps $r_{T,S}$; we require three objects (in our case monoids) $S,T,U$ and such 'exchange' maps which we represent as $\Rightarrow$ as above. The Yang-Baxter equation then has the interpretation that if we 'surface transport' elements along the edges $U,T,S$ shown as `input' around the top and front of the cube to the output using $\Rightarrow$, and keeping the plane of the surface with respect to which the notation in Figure~1(a) is interpreted with normal outward, we get the same answer as the going around the bottom and back of the cube with normal pointing inwards. In other words the net 'surface-holonomy' with normal always outwards right round the cube should be the identity operation. In this way the YBE has the interpretation of zero `higher curvature'. This works in any category (but not a point of view that we have seen before), but in our case fits in with the subdivision property to imply a genuine surface-transport gauge theory. We note that it appears to be somewhat different from  ideas of 2-group gauge theory currently being proposed in the physics literature. Theorem~\ref{tricrossedproduct} will use this to give a matched pair point of view on when an {\bf M3}-monoid  is braided.

Assume $(X,r)$ is a braided
set. We shall mostly be interested in $S=S(X,r)= \langle X| \Re(r) \rangle$ the associated
Yang-Baxter monoid. Clearly, $S$ is graded by length:
\[
 S = \bigcup_{n \geq 0} S^{n}; \quad
S^0 = 1,\quad S^1 = X,\quad S^n = \{u \in S ;  |u|= n \} \] and
\[
S^n.S^m \subseteq S^{n+m}.
\]
In the sequel whenever we speak of a graded monoid $S$ we shall
mean that it is generated  by $S^1=X$ and graded by length. An {\bf M3} or braided monoid $(S, r_S)$ is \emph{graded} if $(S,S)$ is a graded  matched pair.
For a  graded {\bf M3}-monoid $(S,r_S)$, we define the
\emph{restriction} $r= r_X$ of $r_S$ on $X\times X$, where $X$ is the degree 1 generating set. The left and
the right actions respect the grading of $S,$ and condition {\bf
M3} implies
 \[
 r(x,y) = ({}^xy,x^y), \quad r: X\times X \to X\times X.
 \]
 as required for consistency with our notations in Section~2.

\subsection{Matched pair construction of the braided monoid $(S,r_S)$ from $(X,r)$.}

Here, to each braided set $(X,r)$ we will associate a
matched pair $(S,S)$ ($S=S(X,r)$) with left and right actions
uniquely determined by $r,$ which defines a unique `braided monoid'
$(S, r_S)$ associated to $(X,r).$  This is not a surprise given the analogous results for
 the group $G(X,r)$ \cite{Lu} but our approach is necessarily  different. In fact we first construct 
the matched pair or monoids (Theorem~\ref{theoremA}) which is a self-contained result and then consider $r_S$ (Theorem~\ref{theoremAextra}). The reader should  be aware that due to the possible
lack of cancellation in $S$ (we want the statements to be as
general as possible) the proofs of our results for monoids   are
difficult and necessarily involve different computations and combinatorial
arguments. Surely, the results can not be extracted from the
already known results from the group case. Nevertheless, the monoid case is the one naturally arising in this context. Both the monoid  $S(X,r)$ and the
quadratic algebra  $\Acal = \Acal(k,X,r)$ over a field $k$, see
\ref{Adef} are of particular interest. It is known for example that in some special cases of finite
solutions $(X,r),$ $S,$ and $A$ have remarkable algebraic and
homological properties, see \cite{TM, T06, T07}.

We shall extend the left and right actions ${}^x\bullet$ and
$\bullet^x$ on $X$ defined via $r$ , see (\ref{r}),  to a left
action
\[
{}^{(\; )}{\bullet} :S \times S \longrightarrow S
\]
and a right action
\[
{\bullet}^{(\; )} :S \times S \longrightarrow S.
\]
By construction, these actions agree with the grading of $S$, i.e.
$|{}^au|= |u|= |u^a|,$ for all $a,u \in S$.

\begin{theorem}
\label{theoremA} Let  $(X,r)$ be a braided set and
$S=S(X,r)$ the associated monoid. Then the left and the
right actions
\[{}^{(\;\;)}{\bullet}: X\times X  \longrightarrow
 X  , \; \text{and}\;\; \bullet^{(\;\;)}: X
\times X \longrightarrow  X
\]
defined via $r$ can be extended in a unique way to a left and a
right action
\[{}^{(\;\;)}{\bullet}: S\times S  \longrightarrow
 S  , \; \text{and}\;\; \bullet^{(\;\;)}: S
\times S \longrightarrow  S.
\]
which make $S$ a strong graded {\bf M3}-monoid. The associated
bijective map $r_S$ restricts to $r$.
\end{theorem}
\begin{proof}
The proof of the theorem will be made in several steps and we will
need several intermediate results under the hypothesis of the
theorem.

 We will define first left and right actions of $S$ on the free
monoid $\langle X \rangle.$
\[{}^{(\;\;)}{\bullet}: S\times \langle X \rangle \longrightarrow
\langle X \rangle , \; \text{and}\;\; \bullet^{(\;\;)}: \langle X
\rangle \times S \longrightarrow \langle X \rangle
\]
We set
\[
{}^u1:=1, \;1^u:=1 \;\text{for all} \; u \in S.
\]
Clearly, the free monoid $\langle X \rangle$ is graded by
length
\[
\langle X \rangle = \bigcup_{n \geq 0} X^{n}, \; \text{where} \;
X^0 = 1, X^1 = X, X^n = \{u \in \langle X \rangle ;  |u|= n \}
\]
and
\[
X^n.X^m \subseteq X^{n+m}.
\]
 \textbf{Step 1.} We
define the \emph{``actions"}

\begin{equation}
\label{preaction} {}^{(\;\;)}{\bullet}: X\times \langle X \rangle
\longrightarrow \langle X \rangle , \;
\text{and}\;\bullet^{(\;\;)}: \langle X \rangle \times X
\longrightarrow \langle X \rangle
\end{equation} recursively as
follows.  For $n=1$ we have $X^1=X$, and the actions \[
{}^{(\;\;)}{\bullet}: X\times  X^1 \longrightarrow X^1 , \;
\text{and}\;\bullet^{(\;\;)} X^1\times X \longrightarrow
 X^1
\]
are well defined via $r$, see  (\ref{r}).

 Assuming that the
``actions" $X\times X^n \longrightarrow X^n$, and $X^n\times X
\longrightarrow X^n$, are defined for $n$, we will define them for
$n+1$. Let $u \in X^{n+1}.$ Then $u=y.a= b.z,$ where $y,z \in X,
a,b \in X^n .$ We set
\begin{equation}
\label{defleftactiononS} {}^xu = {}^x(y.a): =({}^xy)({}^{x^y}a),
\end{equation}
\begin{equation}
\label{defrightactiononS} u^x=(bz)^x :=(b^{{}^zx})(z^x).
\end{equation}
This way we have defined the ``actions" $X\times X^n\longrightarrow
X^n$, and $X^n\times X\longrightarrow X^n$, for all $n\geq 1$. The next lemma shows now that these definitions are consistent with the multiplication
 in $\langle X\rangle,$ therefore these are well-defined  ``actions" (\ref{preaction}) as required to complete step 1.
 
\begin{lemma}
\label{lemML1} \

a) {\bf ML2} holds for $X$ ``acting" (on the left) on $\langle
X\rangle,$ that is:
\[ {}^x{(ab)}= ({}^xa).({}^{x^a}b)\;\; \text{for
all}\;\; x \in X, a,b \in \langle X\rangle,\; \; |a|,|b| \geq 1.\]
b) {\bf MR2} holds for $X$ ``acting" (on the right) on $\langle
X\rangle,$ that is:
\[ (ab)^x= (a^{{}^bx})b^x\;\; \text{for all}\;\; x \in
X, \; a,b \in \langle X\rangle,\; \; |a|,|b| \geq 1.\]
\end{lemma}
Note that these are equalities of monomials in the free monoid
$\langle X \rangle.$
\begin{proof}
We prove (a) by induction on the length $|a|=n$. Definition
\ref{defleftactiononS} gives the base for the induction. Assume the statement part a) is true for all $a,b \in \langle X\rangle,$
with $|a|\leq n.$ Let $|a|=n+1.$ Then $a=ya_1$, with $ y \in X,
|a_1|=n$. Consider the equalities:
\begin{eqnarray*}
{}^x{(ab)}&=& {}^x{(y(a_1b))}= ({}^xy).({}^{x^y}{(a_1.b)}) \; :
\text{by}\; (\ref{defleftactiononS})
\\
&=& ({}^xy).(({}^{x^y}{a_1}).({}^{((x^y)^{a_1})}b) \; : \text{by the
inductive assumption}\\
&=&({}^x{(ya_1)}).({}^{x^{ya_1}}b) \; : \text{by}\;
(\ref{defleftactiononS})\\
&=&{}^xa.({}^{x^a}b).\end{eqnarray*}
This verifies (a). The proof of (b) is analogous.
\end{proof}

\textbf{Step 2.}  Extend the ``actions" (\ref{preaction}) to a left
and a right actions ${}^{(\;\;)}{\bullet}: S\times\langle X\rangle
\longrightarrow \langle X\rangle , \; \text{and} \;
{\bullet}^{(\;\;)}:\langle X\rangle \times S \longrightarrow
\langle X\rangle $ of $S$ onto $\langle X\rangle$.

Note here that conditions \textbf{l1} (respectively
\textbf{r1}) on $X$ imply that there is a well defined left action
$S\times X \longrightarrow X$, and  a right action $X\times S
\longrightarrow X$ given by the equalities:
\[
{}^{(x_1...x_k)}y : {}^{x_1}{(...({}^{x_{k-1}} {({}^{x_k}y)} )..)}
\]
and
\[
y^{{(x_1...x_k)}} :=(...((y^{x_1})^{x_2})...)^{x_k}.
\]
Clearly ${}^{(ab)}x= {}^a{({}^bx)},$ and $x^{(ab)}= (x^a)^b,$ for
all $x \in X, a,b \in S.$

\begin{proposition}
\label{ML1&MR1} The actions ${}^x\bullet,$ and $\bullet^x$ on
$\langle X\rangle$ extend to ${}^a\bullet,$ and $\bullet^a$ for
arbitrary monomials $a\in S.$ That is to left and right actions
$
 S\times \langle
X\rangle \longrightarrow \langle X\rangle,  \;\; \text{and}\;\;
\langle X\rangle \times S \longrightarrow \langle X\rangle
$
obeying
\[
{\bf ML1}\; \text{for}\; \;{}^S{\langle X\rangle}:\quad {}^{ab}u=
{}^a{({}^bu)},\quad 
{\bf MR1}\; \text{for} \;{\langle X\rangle}^S: \quad u^{ab}=
(u^a)^b
\]
for all $a, b\in S$, $u \in \langle X\rangle$.
\end{proposition}
\begin{proof}
We have to show that the following equalities hold for all  $x,y
\in X, u \in \langle X\rangle$
\begin{equation}
\label{L11} \textbf{L1} \; \text{for} \; {}^S{\langle X\rangle} :\quad 
{}^x{({}^yu)} = {}^{{}^xy}{({}^{x^y}u)}\end{equation}
\begin{equation}
\label{R11} \textbf{R1} \; \text{for} \; {\langle X\rangle}^S :\quad 
(u^x)^y = (u^{{}^xy})^{x^y}.
\end{equation}
We denote these conditions in upper case to distinguish them from {\bf l1,r1} on $X$ itself. We prove (\ref{L11}) by induction on $|u|=n.$ Clearly, when
$|u|=1$, (\ref{L11})
 is simply condition \textbf{l1} on $X$.
Assume  (\ref{L11}) is true for all $u \in X^n.$ Let $u \in
X^{n+1}.$ We can write $u = tv,$ with $t \in X$, $v \in X^n.$ Then
\begin{equation}
\label{eqs3}  {}^x{({}^yu)}={}^x{({}^y(tv))}
=^{Lemma~\ref{lemML1}a)}\;{}^x{(({}^yt).({}^{y^t}v))}=^{Lemma~\ref{lemML1}a)}
({}^x{({}^yt)})({}^{x^{{}^yt}}{({}^{y^t}v)})
\end{equation}
\begin{equation}
\label{eqs3'} =^{\text{inductive
ass.}}({}^x{({}^yt)})({}^{((x^{{}^yt}).(y^t))}v)=^{Lemma~\ref{lemML1}b)}
({}^x{{}^yt})({}^{(xy)^t}v).
\end{equation}
We have shown
\begin{equation}
\label{eqs4} {}^x{({}^yu)= ({}^{xy}t})({}^{(xy)^t}v)= w.
\end{equation}
 Similarly, as in (\ref{eqs3}), (\ref{eqs3'}) we obtain:
\begin{equation}
\label{eqs5} {}^{{}^xy}{({}^{x^y}u)}={}^{{}^xy}{({}^{x^y}tv)} ({}^{{}^xy}{({}^{x^y}t)})({}^{({}^xy.x^y)^t}v)=w_1.
\end{equation}
We use conditions \textbf{L1} and \textbf{R2} on ${}^SX$ to
simplify the right-hand side monomial $w_1$ in this equality:
\begin{equation}
\label{eqs6} w_1= ({}^{{}^xy}{({}^{x^y}t)})({}^{({}^xy.x^y)^t}v)
=^{\textbf{L1},\textbf{R2}}({}^{xy}t)({}^{(xy)^t}v) =w.
\end{equation}
Now  (\ref{eqs4}), (\ref{eqs5}), (\ref{eqs6}) imply
\[
{}^{{}^xy}{({}^{x^y}u)}= w = {}^x{({}^yu)}
\]
We have shown that (\ref{L11}) holds for all $u \in \langle X
\rangle,$ and all $x,y \in X.$ Therefore the equality
\[
{}^{(x_1...x_k)}u : {}^{x_1}{(...({}^{x_{k-1}} {({}^{x_k}u)} )..)}
\]
gives a well defined left action $S\times \langle X \rangle
\longrightarrow \langle X \rangle$ , which satisfies \textbf{ML1}.

The proof of (\ref{R11}) is analogous. So the right action
$\langle X \rangle \times S  \longrightarrow \langle X \rangle$ is
also well defined, and satisfies \textbf{MR1}.
\end{proof}

\begin{proposition} \ 
\label{propML2}
\[
\begin{array}{lclc}
 {\rm\bf ML2\ for\ }{}^S\langle X\rangle:\quad& {}^a{(uv)}=({}^au).({}^{(a^u)}v),
\quad& {\rm\bf MR2\ for\ }\langle X\rangle^S:\quad  &(vu)^a=(v^{{}^ua}).(u^a),
\end{array}\]
for all $a \in S,  u, v\in  \langle X \rangle$.
\end{proposition}
\begin{proof}
We prove first the  following equalities:
\begin{equation}
\label{ML3}
\begin{array}{lclc}
{}^a{(yv)}= ({}^ay)({}^{a^y}v), \quad & (vy)^a=(v^{{}^ya})(y^a),
\text{for all}\; \; a\in S,  y\in X, v\in \langle X \rangle.
\end{array}\end{equation}
The proofs of the two identities in (\ref{ML3}) are analogous. In
both cases one uses induction on $|a|$. We shall prove the
left-hand side equality. By the definition of the left action one
has ${}^y{xv}=({}^yx)({}^{y^x}v).$ This gives the base for the
induction. Assume (\ref{ML3}) holds for all $a,$ with $|a|\leq n$.
Let $a_1 \in S, $ $|a_1|= n+1$, $y \in X,$ $v \in \langle X
\rangle.$ Clearly, $a_1=ax,$ where $x \in X,$ $a\in S, |a|=n.$
Consider the equalities
\[
{}^{ax}{(yv)}=^{\textbf{ML1}}
{}^a{({}^x{(yv)})}=^{(\ref{defleftactiononS})}\;{}^a{[({}^xy)({}^{x^y}v)]}
=^{\text{inductive ass.}} \;
[{}^a{({}^xy)}][{}^{(a^{{}^xy})}{({}^{x^y}v)}]
\]
\[
 =^{\textbf{ML1}}[{}^{(ax)}y][{}^{((a^{{}^xy}).(x^y\emph{}))}v] =^{(\ref{defrightactiononS})}
[{}^{(ax)}y][{}^{(ax)^y}v].
\]
We have verified the left-hand side of (\ref{ML3}).

Now we verify Proposition~\ref{propML2}) using again induction on $|a|$.
\textbf{Step 1.} Base for the induction.  The equality
${}^x{(uv)}=({}^xu).({}^{x^u}v)$, for all $x \in X$, $u, v \in
\langle X \rangle,$ is verified by Lemma~\ref{lemML1}. \textbf{Step
2.} Assume the statement holds for all monomials $a\in S$, with
$|a|\leq n$. Let $a_1 = ax \in S$, where $|a|=n$, $x \in X.$ The
following equalities hold in $\langle X \rangle$:
\[
{}^{ax}{(uv)}=^{\textbf{ML1}}{}^a{({}^x{(uv)})}={}^a{({}^xu).({}^{x^u}v)}
 =^{\text{inductive ass.}}
[{}^a{({}^xu)}][{}^{(a^{{}^xu})}{({}^{x^u}v)}]
\]
\[
=^{\textbf{ML1}} [{}^{ax}u][{}^{((a^{{}^xu})(x^u))}v]
=^{(\ref{ML3})} [{}^{ax}u][{}^{(ax)^u}v].
\]
This verifies the left-hand half for the stated proposition.
 Analogous argument proves the right-hand half. 
\end{proof}
We will show next that the left and right actions of $S$ on
$\langle X\rangle$ defined and studied above induce naturally left
and right actions  ${}^{(\;)}{\bullet}: S\times S \longrightarrow
S $, and ${\bullet}^{(\;)}: S\times S \longrightarrow S $.
 We need to verify  that the actions
 agree with the relations in $S$. We start with the following
lemma, which gives an analogue of \textbf{L2}, but for $S$ acting
on monomials of length 2.
\begin{lemma}
\label{beautifullemma} The following equalities hold in $X^2$.
\[\begin{array}{lclc}
 \quad& {}^a{(yz)}={}^a{({}^yz.y^z)},
\quad& {(yz)}^a={({}^yz.y^z)}^a ,\quad  &\text{for all} \; a\in S;
y,z \in X.
\end{array}\]
\end{lemma}
\begin{proof} As usual we use induction on the
length $|a|=n.$

\textbf{Step 1.} $|a|=1$, so $ a =x \in X.$ Then by \textbf{l2} on
$X$ we have
\[
{}^x{(yz)}={}^x{({}^yz.y^z)},\quad  \text{for all} \; x, y,z \in X.
\]
\textbf{Step 2.} Assume the statement of the lemma holds for all $y,z \in X,$ and
all $a \in S$, with $|a|=n.$ Let $a_1 \in S, |a_1|=n+1,$ so
$a_1=ax, x \in X, a\in S, |a|=n.$ Now
\begin{equation}
\label{eqs11}
{}^{(ax)}{(yz)}=^{\textbf{ML1}}{}^a{[{}^x{(yz)}]}={}^a{[({}^xy)({}^{x^y}z)]}
 =^{\text{induct.\ 
ass.}}{}^a{[({}^{({}^xy)}{({}^{x^y}z)}).({}^xy)^{({}^{x^y}z)}]}={}^a{(uv)},
\end{equation}
where for convenience we denote
\[
u = {}^{({}^xy)}{({}^{x^y}z)}=^{ \textbf{l1}} {}^x{({}^yz)},\quad v= ({}^xy)^{({}^{x^y}z)}=^{\textbf{lr3}}= {}^{x^{{}^yz}}{(y^z)}.\]
So
\begin{equation}
\label{eqs12} uv= ({}^x{({}^yz))({}^{x^{{}^yz}}{(y^z)})}=^{Lemma~\ref{beautifullemma}} {}^x{({}^yz.y^z)}.
\end{equation}
is an equality in $X^2.$
 It follows from (\ref{eqs11}) and (\ref{eqs12}) that
\[
{}^{(ax)}{(yz)}= {}^a{(uv)}=^{(\ref{eqs12})}
{}^a{({}^x{({}^yz.y^z)})}=^{\textbf{ML1}}{}^{ax}{({}^yz.y^z)}.
\]
This proves the first equality stated in the lemma. Analogous argument
verifies the second.
\end{proof}

The following statement shows that the left action agrees with all
replacements coming from the defining relations of $S$, and
therefore agrees with equalities of words in $S$
\begin{proposition}
\label{agreementwithreductions} The left and right actions
${}^{(\;)}{\bullet}: S\times S\longrightarrow S$, and
${\bullet}^{(\;)}: S\times S\longrightarrow S$ are well defined.
They make $(S,S)$ a graded matched pair of monoids.
\end{proposition}
\begin{proof}
We need to verify the  following equalities  in $S$:
\begin{equation}
\label{reduction} {}^a(u.yz.v)= {}^a(u.({}^yz.y^z).v).
\end{equation}
It easily follows from (\ref{reduction}) that
\[
w_1 = w_2 \;\;\text{is an equality in}\; S \; \Longrightarrow
{}^a{w_1} ={}^a{w_2} \;\;\text{is an equality in} \; S.
\]
Indeed,  $w_1=w_2$ is an equality in $S$ if and only if $w_2$ can
be obtained from $w_1$ after a finite number of replacements
coming from the relations $\Re(r)$. Note that each relation in
$\Re(r)$ has the shape $yz = {}^yz.y^z,$ where $y,z \in X.$

We prove now (\ref{reduction}). Let $a,u,v \in S, y,z \in X.$ Then
\[
{}^a(u.yz.v)=[{}^au].[{}^{a^u}{(yz)}].[{}^{a^{u(yz)}}v]
=^{Lemma~\ref{beautifullemma}} [{}^au].[{}^{a^u}{({}^yz.y^z)}].[{}^{a^{u(yz)}}v]
\]
\[ = {}^a{[u.({}^yz.y^z)]}.[{}^{a^{u(yz)}}v]=^{(\ref{eqs15})}
{}^a{[u.({}^yz.y^z)]}.[{}^{(a^{(u({}^yz.y^z))})}v]=^{\textbf{ML2}}
{}^a[u.({}^yz.y^z).v].
\]
We used above the following equality implied by condition {\bf
M3}:
\begin{equation}
\label{eqs15} a^{u(yz)}=((a^u)^y)^z)=((a^u)^{{}^yz})^{y^z}.
\end{equation}
\end{proof}
\begin{proposition}
 \label{compatibleactions} $S$ in
Proposition~\ref{agreementwithreductions} is an  \textbf{M3}-
monoid (i.e. for every $u,v \in S$, the equality $uv={}^uv.u^v$
holds in $S.$)
\end{proposition}
\begin{proof}
Using induction on $|w|$ we first show that there is an equality
in $S$
\begin{equation}
\label{eqs10} xw={}^xw.x^w \; \text{for all} \; w \in S, x \in X.
\end{equation}
When $|w|=1$, one has $r(xw)={}^xw.x^w,$ therefore (\ref{eqs10})
either   belongs to the set of defining relations $\Re(r)$
(\ref{defrelations}) for S, or is an equality in $X^2$. Assume
(\ref{eqs10}) is true for all $w$ with $|w|\leq n.$ Let $v \in S,
|v|= n+1,$ $x \in X.$ Present $v=yw,$ where $y \in X, w\in S, |w|n.$ The following equalities follow from the associativity of the
multiplication in $S$, the inductive assumption, and \textbf{ML2}.
\[
x.v = x.(yw) = (xy)w = ({}^xy.x^y)w = ({}^xy)(x^y.w)({}^xy)[({}^{x^y}w){(x^y)}^w]
\]
\[
=^{\textbf{MR1}} [({}^xy)({}^{x^y})w][x^{(yw)}] =^{\textbf{ML2}}\;
[{}^x{(yw)}].[x^{(yw)}].
\]
This proves (\ref{eqs10}) for all $x, w,$ $x\in X, w \in S.$

Next we use induction on $|u|$ to prove the statement of the
proposition.

 Assume
\[
uv={}^uv.u^v\quad \text{for all}\quad u,v \in S,
\quad\text{with}\quad |u|\leq n.
\]
 Let
$u_1 = ux,$ where $u\in S, |u|= n,$ $x \in X.$ Then the
associative law in $S$, the inductive assumption,
\textbf{ML1}, and \textbf{MR2} imply the following equalities
\[
(ux)v=u(xv)=^{(\ref{eqs10})}u({}^xv.x^v) = [u.{}^xv ].x^v\quad\quad\quad\quad\quad
\]
\[
= [{}^u{({}^xv)}].[u^{({}^xv)}).(x^v)]=^{\textbf{ML1},
\textbf{MR2}}\; [{}^{(ux)}v].[(ux)^v]
\]
\end{proof}
We have shown that the left and the right actions constructed
above make $(S,S)$ a graded matched pair of monoids with condition
\textbf{M3}. Moreover the matched pair axioms imply that these
actions are uniquely determined by $r$. So we  define the
associated braided monoid $(S,r_S),$ with $r_S(u,v)=({}^uv, u^v)$,
which is also uniquely determined by $r$.

It remains to show that $(S,S)$ is a strong matched pair. By
Theorem~\ref{theoremAextra}.1 $(S,r_S)$ is a braided set, in
particular $r_S$ is a bijective map, therefore by  Proposition
\ref{M3impliesbraidedth}.2 $(S,S)$ is a strong matched pair.
 This completes the proof of Theorem 
\ref{theoremA} \end{proof}

In the case when $(S,r_S)$ is a
graded {\bf M3}-monoid generated by $X$, $(X=S^1)$, by definition,
the set $X$ is invariant under the left and the right actions. The next proposition
will be used in this case, but later on we will also need it more generally and 
therefore we prove it now in this greater generality. Thus, suppose $(S,r_S)$ is an {\bf
M3}-monoid (not necessarily graded). Suppose $S$ is generated  by
a set $X$ which is invariant under the left and the right actions
of the matched pair. We define the restriction $r=r_X$ of
$r_S$ to $X\times X$, so $r(x,y) = ({}^xy,x^y)$ as
usual. Then $(X,r)$ is a set with a quadratic map and condition
{\bf M3} implies that  the set of defining relations $\Re(S)$
of $S$ satisfy
\[
  \Re(S) \supseteq \Re(r),
\]
where in general  there might be a strict inequality ($\Re(r)$ are the
quadratic relations coming from $r$, see (\ref{defrelations})).

The monoid $S$ can be presented as $S=\langle X; \Re(S) \rangle$
and, it is, in general, homomorphic image of the monoid
$S(X,r)=\langle X; \Re(r) \rangle.$ Clearly, $S$ is graded
\emph{iff} $\Re(S)$ consists of homogeneous relations. All
relations in $S$ are results of replacements coming from the set
of defining relations $\Re(S)$. In particular, a relation of the
shape $x_{i_1}\cdots x_{i_m}=1$ can hold in $S$ \emph{iff } the set
$\Re(S)$ explicitly contains relations of the form $u=1$ where $u\in\langle X\rangle$ has length $|u|\le m$.   Every monomial $u \in S$
can be presented as a product
\begin{equation}
\label{eq_u} u=x_{i_1}\cdots x_{i_m}
\end{equation}
 of elements of $X$, (since $X$ is a generating set) but in general, there
might be presentations of $u$ as words of different length.
However we can always consider a presentation (\ref{eq_u}) in
which the length $m=m(u)$ is minimal. Clearly, for each $u\in S$
the minimal length $m(u)$ is uniquely determined. Furthermore,
$m(u.v) \leq m(u)+m(v).$ Clearly,   every sub-word $a$ of a
"minimal" presentation of $w$ is also a presentation of minimal
length for $a,$ i.e. if $w=^{in S}w_1w_2,$  $w_1,w_2\in
\langle X \rangle $ is a presentation of minimal length, then
$m(w)= m(w_1) +m(w_2).$ In particular,  if $w=xu, x \in X, u\in S$
is a presentation of $w$ of minimal length then $m(u)= m(w)-1$.

\begin{proposition} \label{LR3onM3monoids}
Let $(S, r_S)$ be an {\bf M3}-monoid with associated map $r_S$.
Suppose $S$ is generated by a set $X$ which is invariant under the
left and the right actions of the matched pair. Let $r: X\times
X\to X\times X$ be the restriction of $r_S$ on $X\times X.$
Suppose condition {\bf lr3} holds on the quadratic set $(X,r).$ Then
\[
\begin{array}{lclc}
 {\rm\bf LR3:}\quad& {({}^aw)}^{({}^{a^w}{b})}={}^{(a^{{}^wb})}{(w^b)},\quad&\text{for
 all}\; a,b, w \in S.
\end{array}\]
\end{proposition}
\begin{proof}
Using induction on $m(b)=n$ we  prove first
\begin{equation}
\label{LR3Sright}
 {({}^zt)}^{({}^{z^t}{b})}={}^{(z^{{}^tb})}{(t^b)}, \quad\text{for
 all}\; b \in S, t,z  \in  X.
 \end{equation}
Clearly, when $m(b)=1,$ condition (\ref{LR3Sright}) is simply
\textbf{lr3} on $(X,r)$, which gives the base for induction.
Assume (\ref{LR3Sright}) is true for all $b \in S, m(b)\leq n,$
and all $t, z \in X.$ Let $b \in S, m(b)= n+1,$ so $b=xu,$ $x \in
X, u \in S, m(u)= n.$ Consider the equalities:
\begin{eqnarray*}
{}^{(z^{{}^t{(xu)}})}{[t^{xu}]}&=&
{}^{(z^{({}^tx)})^{({}^{t^x}u)})}{[(t^x)^u]} =^{\text{inductive
ass.}}[{}^{(z^{{}^tx})} {(t^x)}]^{{}^{[(z^{{}^tx})^{t^x}]}u}\\
&&=^{\textbf{lr3}, \textbf{ML2}}\
[({}^zt)^{{}^{z^t}x}]^{{}^{(z^{(tx)})}u}= ({}^zt)^{[({}^{z^t}x).({}^{(z^t)^x)}u]}=({}^zt)^{({}^{z^t}{(xu)})}
\end{eqnarray*}
This proves (\ref{LR3Sright}).

Analogous argument with induction on $m(a)=n$ verifies
\begin{equation}
\label{LR3Srightleft}
  {({}^at)}^{({}^{a^t}{b})}={}^{(a^{{}^tb})}{(t^b)},\quad \text{for
 all}\; a,b \in S, t  \in  X .
 \end{equation}
Note that  (\ref{LR3Sright}) gives the base for the induction.

Finally  we prove  {\rm\bf LR3} as stated. We use induction on $m(w)=n.$  In the case $m(w)=1,$ {\rm\bf LR3} is exactly  condition (\ref{LR3Srightleft}). This is the base for
the induction. Assume {\rm\bf LR3} holds for all $a,b, w \in S$,
where $m(w)\leq n$. Let  $w \in S, m(w)= n+1,$ then $w=xv,$ $x\in
X,$ $m(v)= n.$ We have to show
\begin{equation}\label{LR3induction}
{}^{(a^{{}^{(xv)}b})}{[(xv)^b]} = [{}^a{(xv)}]^{({}^{a^{(xv)}}b)}
\end{equation}
We compute
\begin{eqnarray*}
W_1 &\equiv & {}^{(a^{{}^{(xv)}b})}{[(xv)^b]}
=^{\textbf{MR2}}{}^{(a^{{}^{(xv)}b})}{[(x^{{}^vb}).(v^b)]}=^{\textbf{ML1}}{}^{a^{{}^x{({}^vb)}} }{[(x^{{}^vb}).(v^b)]}
\\
&&=^{\textbf{ML2}}[{}^{  a^{{}^x{({}^vb)}} }{(x^{{}^vb})}].[{}^ {
(a^{ {}^x{({}^vb)} }  )^{(x^{{}^vb}) }} {(v^b)}]
=^{\textbf{MR1}}[{}^{  a^{{}^x{({}^vb)}} }{(x^{{}^vb})}].
[{}^{(a^{(x.{}^vb)})}{(v^b)}]\\
&&=^{\text{inductive ass.}}[({}^ax)^{({}^{a^x}{({}^vb))}}]. [{}^{(
(a^x)^{{}^vb} )
} {(v^b)}]\\
W_2&\equiv & [{}^a{(xv)}]^{{}^{(a^{(xv)})}b}=[({}^ax).({}^{a^x}v)]^{{}^{((a^x)^v)}b}
= [({}^ax)^{{}^{({}^{a^x}v)}{{}^{((a^x)^v)}b}}]
.[({}^{a^x}v)^{{}^{((a^x)^v)}b}]\\
&&=^{\textbf{ML1}} [({}^ax)^{({}^{a^x}{({}^vb))}}].
[({}^{a^x}v)^{{}^{((a^x)^v)}b}]=^{\text{inductive ass.}}[({}^ax)^{({}^{a^x}{({}^vb))}}]. [{}^{(
(a^x)^{{}^vb} ) } {(v^b)}].
\end{eqnarray*}
We have shown $W_1=W_2$  which proves (\ref{LR3induction}). Hence {\bf LR3} is proven on $S$.
\end{proof}

We have completed all the necessary parts for the construction of
a solution $r_S$ of the YBE on $S$.
 \begin{theorem}
\label{theoremAextra} Let $(X,r)$ be a braided set and
$(S,r_S)$ the induced graded {\bf M3}-monoid with associated map $r_S$ by
Theorem~\ref{theoremA}. Then
\begin{enumerate}
\item \label{theoremAextra1}$(S,r_S)$ is a graded braided
monoid.
 \item \label{theoremAextranondeg}
$(S, r_S)$ is non-degenerate iff $(X,r)$ is non-degenerate.
\item\label{theoremAextrainvol} $(S, r_S)$ is involutive iff
$(X,r)$ is involutive.
\end{enumerate}
\end{theorem}
\begin{proof}
We  have already in Theorem~\ref{theoremA} extended the left and the right actions on $X$ induced by
$r$ to a left and a right actions of $S$ onto itself which make
$S$ a graded   a graded \textbf{M3}-monoid with associated map $r
_S$. To show that $r_S$ obeys the YBE we use Lemma~\ref{ybe} applied to the 
set $S$. We already have {\bf l1,r1} on $S$ from {\bf ML1,MR1,M3} by the same argument as in the proof of Proposition~\ref{M3cancellative}. $(S,r_S)$ satisfies the
hypothesis of Proposition~\ref{LR3onM3monoids}, therefore
we also have \textbf{lr3}  on $S$ as required.

For bijectivity of $r_S$ we consider $(X,r^{-1})$ (that $r$ is
bijective is our convention throughout the paper) and consruct a
matched pair $(T,S)$ (where $T=S$) using our previous results
applied to $(X,r^{-1})$. Thus, we define $\la,\ra$ by
$r^{-1}(x,y)=(x\la y,x\ra y)$. One may then prove inductively with
respect to the grading that these obey (\ref{lara}) and hence
provide the inverse of $r_S$. Thus, suppose these equations for
all $a,u$ of given degrees and also for $x$ of degree 1 in the
role of $a$. Then
\begin{eqnarray*}{}^{xa}u\la (xa)^u&=&{}^x({}^au)\la (x^{{}^au}a^u)={}^x({}^au)\la x^{({}^au)}.
({}^x({}^au)\ra x^{({}^au)})\la a^u=x.({}^au\la a^u)=xa\\
{}^{xa}u\ra (xa)^u&=&({}^x({}^au)\ra x^{({}^au)})\ra a^u={}^au\ra
a^u=u\end{eqnarray*} using that both sets of actions form matched
pairs and the assumptions. Similarly for the other cases. At the
lowest level (\ref{lara}) holds as  $r^{-1}$ is inverse to $r$. We
have shown that $(S,r_S)$ is a braided monoid.

 Next we observe one direction of part
(\ref{theoremAextranondeg}). By definition the solution
$r_S$ respects the grading by length of $S,$ hence clearly, the
non-degeneracy of $r_S$ implies that $r$ is non-degenerate. For
the converse we will need some additional statements. We recall a
standard notation.

\begin{notation}
For $1\leq i,$  the replacements $r_i, 1 \leq i\leq m-1$ are
defined as
\[
r_i(y_1y_2 \cdots y_m)= y_1\cdots y_{i-1}r(y_iy_{i+1})y_{i+2}
\cdots y_m
\]
\[ =y_1\cdots y_{i-1}[{}^{y_i}{y_{i+1}}{y_i}^{y_{i+1}}]y_{i+2}
\cdots y_m,
\]
for all $m > i,$ and $y_1, \cdots, y_m \in X.$
\end{notation}
 Note that, clearly, each replacement $r_i$ agrees
with the defining relations of $S,$ so for all $u \in \langle
X\rangle,$ and all $i < |u|$ the equality $r_i(u)=u$ holds
in $S.$

\begin{lemma}
\label{beautifullemma2} For any $a \in S,  u \in \langle
X\rangle,$ and any sequence $r_{i_1}\cdots r_{i_s},$ with $\;1
\leq i_j < |u|,$ for $ 1 \leq j \leq s,$ the following is
an equality of words in $\langle X\rangle$
\[{}^a(r_{i_1}\circ \cdots \circ r_{i_s}(u)) = r_{i_1}\circ\cdots \circ r_{i_s}({}^au)
\]
\end{lemma}
\begin{proof}
We show first that ${}^{a}{r_i(u)}= r_i({}^au),$ for all $a \in
S,$ $u \in \langle X,\rangle, 1 \leq i < | u|.$ For our
convenience we present $u$ as $u =vxyw,$ where $|v| = i-1,
x,y \in X.$ The following are equalities  in $\langle X\rangle$:
\begin{eqnarray*}
r_i({}^au) &=& r_i({}^a{vxyw})=^{\textbf{ML2}} r_i({}^av.
({}^{a^v}x.{}^{a^{vx}}y).{}^{a^{vxy}}w)\\
&&=({}^av).[r({}^{a^v}x.{}^{a^{vx}}y)].({}^{a^{vxy}}w) =({}^av).
[{}^{{}^{a^v}x}{({}^{(a^v)^x}y)}.({}^{a^v}x)^{{}^{a^{vx}}y}].{}^{a^{vxy}}w\\
&&=^{\textbf{ML1}} ({}^av).
({}^{(a^v.x)}y).({}^{a^v}x)^{({}^{{(a^v)}^x}y)}].{}^{a^{vxy}}w\\
&&=^{\textbf{ML1,
LR3}}({}^av).[({}^{a^v}{({}^xy)}).{}^{(a^v)^{{}^xy}}{(x^y)}].{}^{a^{vxy}}w\\
&&=^{\textbf{ML2}} {}^av.[{}^{a^v}{({}^xy.x^y)}].{}^{a^{vxy}}w
=^{\textbf{M3}}{}^av.[{}^{a^v}{({}^xy.x^y)}].{}^{a^{v({}^xy.x^y)}}w\\
&&=^{\textbf{ML2}} {}^a{(v.({}^xy.x^y).w)}=
{}^a{r_i(v.xy.w)}={}^a{r_i(u)}
\end{eqnarray*}
The statement of the lemma then follows easily and we leave it to the
reader.
\end{proof}
\begin{proposition}
\label{rsnondegenerateprop} Suppose $(X,r)$ is non-degenerate,
then for any $a, b, u,v \in S$
\begin{enumerate}
\item \label{rsleftnondeg} ${}^au={}^av \quad \text{in $S$}
\Longrightarrow u=v \quad \text{in $S$}$
\item\label{rsrightnondeg}
$a^u=b^u \quad \text{in $S$} \Longrightarrow
a=b \quad \text{in $S$}$
\item \label{rsnondeg}$(S, r_S)$ is a non-degenerate solution.
\end{enumerate}
\end{proposition}
\begin{proof}
We will show first (\ref{rsleftnondeg}). From the definition of
the left action of $S$ on $\langle X \rangle$ and the
non-degeneracy of $r$ we deduce
\begin{equation}\label{eN1}
{}^a(y_1y_2\cdots y_m)={}^a(z_1z_2\cdots z_m) \; \text{in} \;
\langle X \rangle \Longrightarrow y_1y_2\cdots y_m =z_1z_2\cdots
z_m \; \text{in} \; \langle X \rangle.
\end{equation}
 Suppose ${}^au={}^av$
holds in $S$, where $a, u, v \in S.$ Clearly then  the monomial
${}^av,$ considered as an element of $\langle X \rangle$ is
obtained from ${}^au,$ by applying finite sequence of replacements
which come from the defining relations. So there exist  $r_{i_1},
r_{i_2}, \cdots, r_{i_s},$ such that there is an equality in
$\langle X \rangle$:
\[r_{i_1}\circ\cdots \circ r_{i_s}({}^au) = {}^av.
\]
It follows from Lemma~\ref{beautifullemma2} that
\[ {}^a(r_{i_1}\circ \cdots \circ r_{i_s}(u))= r_{i_1}\circ\cdots \circ
r_{i_s}({}^au) ={}^av \quad\text{in}\;\; \langle X \rangle.
\]
Hence, by (\ref{eN1}), the equality  \[ r_{i_1}\circ \cdots \circ
r_{i_s}(u) = v\]
 holds in $\langle X \rangle.$  Then, clearly, $u=v \in S.$ This verifies the left non-degeneracy of $r_S,$
(\ref{rsleftnondeg}). The proof of right non-degeneracy
(\ref{rsrightnondeg}) is analogous. (\ref{rsleftnondeg}) and
(\ref{rsrightnondeg}) imply (\ref{rsnondeg}).
\end{proof}
This completes the part (2) of Theorem~\ref{theoremAextra}. Now
we will prove part (\ref{theoremAextrainvol}). Clearly, the
involutiveness of $(S, r_S)$ implies $(X,r)$ involutive. The
following lemma gives the opposite implication.
\begin{lemma}
\label{semigroupsolutioninvollemma} Under the hypothesis of the
theorem, assume $(X,r)$ is an involutive solution. Then $(S, r_S)$
is involutive.
\end{lemma}
\begin{proof}
We have to show that
equalities in $S:$
\begin{equation}
\label{eqinv} {}^{{}^uv}{(u^v)}=u, \quad ({}^uv)^{u^v}=v, \quad
\text{for all}\quad u,v \in S. \end{equation} We show first the two equations
\begin{equation}
\label{eqinv1} {}^{{}^{xy}z}{(xy^z)}=xy
 \quad \text{for all}\quad x,y,z \in X, \end{equation}
\begin{equation}
\label{eqinv2} ({}^{xy}z)^{({xy}^z)}=z
 \quad \text{for all}\quad x,y,z \in X. \end{equation}
 Indeed,
\begin{eqnarray*}
{}^{{}^{xy}z}{(xy^z)}&=^{\textbf{MR2}}&{}^{{}^{xy}z}{(x^{{}^yz}.y^z)}=^{\textbf{MR1, ML2}} [{}^{{}^{x}{({}^yz)}}{(x^{{}^yz})}].
{}^{{({}^{x}{({}^yz)})}^{(x^{{}^yz})}}{(y^z)}\\
&&=^{r\;\text{involutive}} x.{}^{{}^yz}{(y^z)}
=^{r\;\text{involutive}} xy.\end{eqnarray*}
This proves (\ref{eqinv1}). For the second equality, we have
\begin{eqnarray*}
({}^{xy}z)^{({xy}^z)}&=^{\textbf{ML1, MR2}}&
[{}^{x}{({}^yz)}]^{(x^{{}^yz}.y^z)}\\
&=^{\textbf{ML1}}&
[[{}^{x}{({}^yz)}]^{(x^{{}^yz})}]^{y^z}=^{r\;\text{involutive}}({}^yz)^{y^z}=z,
\end{eqnarray*}
hence (\ref{eqinv2}) also holds. We give a sketch of the proof of
${}^{{}^uv}{(u^v)}=u,$ for all $u,v \in S,$ and leave the details
for the reader. One uses double induction on $| u|=m$ and
$|v|=n.$  \textbf{Step 1. } By induction on $| v|=n,$
one proves the equality ${}^{{}^zv}{(z^v)}=z, z\in X, v \in S$. To
do this one uses a tecnique similar to the one in the proof of
(\ref{eqinv2}). \textbf{Step 2.} Assuming ${}^{{}^uv}{(u^v)}=u,$ for
all $u, v \in S, |u| \leq m,$ one uses argument similar to
the proof of (\ref{eqinv1}) to show that ${}^{{}^{xu}v}{(xu^v)}=xu$,
for all $x \in X, |u| = m, u,v, \in S.$ This proves the
first equality in (\ref{eqinv}). The second is proven by analogous
arguments.\end{proof} The proof of Theorem~\ref{theoremAextra} is complete now.
\end{proof}

\subsection{{\bf M3}-monoids  and braided monoids $(S,r_S)$}

In this section we will study general \textbf{M3}-monoids $S$ with associated
map $r_S$. We recall that \emph{a braided group}, \cite{Takeuchi}, is a pair
$(G,\sigma),$ where $G$ is a group and $\sigma : G\times
G\rightarrow G\times G$ is a map such that the left and the right
actions induced by $\sigma$ make $(G,G)$ a matched pair of groups
with condition \textbf{M3}. Note that in the case of matched pair
of groups the notions of an {\bf M3}-group $(G, r_G)$ and a
braided group are equivallent;  it is
shown in \cite{Lu}, Theorem 1 that a braided group $(G,\sigma)$
forms a non-degenerate braided set. The proof follows
straightforwardly from the definition of a braided group and the
cancellative law in $G$. Analogous argument can not be applied in
the general case of \textbf{M3}-monoid $(S, r_S)$. We show that an
\textbf{M3}-monoid $(S, r_S)$ is braided in cases when some
natural additional conditions on $S$ are imposed.

In various cases, when necessary,  we shall impose the condition
of 2-cancellativity.
\begin{definition} We
say that a monoid $S$ is \emph{2-cancellative} with respect to a  generating set $X$  when it has cancellation on monomials of length 2 in the generators in the sense:
\[xy=xz \Longrightarrow y=z; \quad xz=yz \Longrightarrow x=y,
\quad \text{for all}\quad x,y,z \in X.\]
\end{definition}

 Clearly  for a graded {\bf M3}-monoid $(S,r_S)$,
 the  2-cancellativity of $S$  implies that $r$  is
 2-cancellative.
 Recall that  in the particular case when $S = S(X,r)$,
 $S$ is 2-cancellative \emph{iff} $r$  is 2-cancellative,  see
 Definition~\ref{cancellsoldef} and Proposition~\ref{cancellsolprop}. Furthermore, by Corollary~\ref{involcancelcorol} every nondegenerate
 involutive quadratic set   $(X,r)$
is 2-cancellative, so the
 condition of 2-cancellativeness is a natural restriction.
\begin{lemma}
\label{S,Smatchedlemma} Let $S$ be an {\bf M3}-monoid with
associated map $r_S$. Suppose $S$ is 2-cancellative with respect to a generating set $X$ which
is invariant under the left and the right actions of the matched
pair. Let $r: X\times X\to X\times X$ be the restriction of $r_S$ on
$X\times X.$ Then $(X,r)$ is a set-theoretic solution of the
Yang-Baxter equation. Furthermore, if  the map $r$ is a bijection
then $(X,r)$ is a braided set.
\end{lemma}
\begin{proof}
We have  to show that $(X,r)$ satisfies the relations
\textbf{lr3}. The following equalities hold in $S$ :
\[{}^{x}(yz) = ^{\textbf{ML2}} {}^{x}y. {}^{x^y}z =^{\textbf{M3}} {}^{{}^{x}y}({}^{x^y}z)).({}^{x}y)^ {{}^{x^y}z}\]
and
\[
{}^{x}{(yz)}=^{\textbf{M3}} {}^{x}{({}^yz.y^z)}=^{\textbf{ML2}}
{}^{x}{({}^yz)}.{}^{x^{{}^yz}}{(y^z)}.
\] Comparing the right hand sides of these equalities we obtain
\begin{equation}
\label{LR3forgradedmp1} {}^{{}^{x}y}{({}^{x^y}z)}.({}^{x}y)^
{{}^{x^y}z}= {}^{x}({}^yz).{}^{x^{({}^yz)}}{(y^z)}. \end{equation}
By {\textbf{ML1, M3}} one has
${}^{{}^{x}y}{({}^{x^y}z)}={}^{x}{({}^yz)},$ which together with 
(\ref{LR3forgradedmp1}) and the 2-cancellativity yields
\[
({}^{x}y)^{{}^{x^y}z}={}^{{x}^{{}^yz}}{(y^z)}.
\]
This verifies \textbf{lr3} on $(X,r)$. Clearly, the conditions
\textbf{ML1, MR1 } give  \textbf{l1, r1 } on $(X,r),$ hence
$(X,r)$ is a solution.
\end{proof}

\begin{proposition}
 \label{M3impliesbraidedth}
Let  $S$ be an \textbf{M3}-monoid with associated map $r_S$.
Suppose  one of the following  conditions is satisfied:

(i) $S$ has a generating set $X$ invariant under the left and the
right actions of the matched pair, the restriction $r$ of $r_S$ on
$X\times X$ is bijective, and {\bf lr3} holds on the quadratic set
$(X,r)$;

(ii) $S$ has a generating set $X$ invariant under the left and the
right actions of the matched pair, $S$ is 2-cancellative with respect to it, and the
restriction $r$ of $r_S$ on $X\times X$ is bijective;

(iii) $(S,r_S)$ is graded, $S$ is 2-cancellative,  and the
restriction  $r: X\times X\to X\times X$ on degree 1 is a
bijection;

 (iv)
$S$ is a monoid (not necessarily graded) with left cancellation.

  Then $(S,r_S)$ is a set-theoretic solution of the YBE.
 \end{proposition}
 \begin{proof} As in the proof of Proposition~\ref{M3cancellative} we only need to have {\bf lr3} on the set of $S$, which condition we denote in upper case as {\bf LR3} to avoid confusion with the condition on $X$. 
 
 If (i) holds then  $(S,r_S)$ satisfies the
hypothesis of Proposition~\ref{LR3onM3monoids}, so condition {\bf
LR3} holds and $(S,r_S)$ is a solution.

If (ii) holds then by Lemma~\ref{S,Smatchedlemma} {\bf lr3}
holds on the quadratic set $(X,r),$ so condition (i) is satisfied
therefore $(S,r_S)$ is a solution.

Condition (iii) implies that $X$ is invariant under the left and
the right actions of the matched pair, therefore the assumptions
of (ii) are satisfied, and $(S,r_S)$ is a solution.

 Case (iv) is already covered in Proposition~\ref{M3cancellative} and included for completeness. 
  \end{proof}
  
At this level the conceptual basis for why $(S,S)$ with
cancellation (or rather with \textbf{LR3}) has a solution of the
YBE on it is then provided by the following. Note that the
monoids $S$ are somewhat analogous to the `FRT bialgebras'
$A(R)$ in the theory of quantum groups and just as there one
has\cite{Ma:book} that  $A(R)\bowtie A(R)\to A(R)$, similarly we
have a monoid homorphism $S\bowtie S\to S$ given by $u.a\mapsto
ua$ (i.e. by the product in $S$). In fact this is the exact
content of the \textbf{M3} condition and is analogous to the
parallel observation in the group case \cite{Takeuchi}. Likewise,
just as one has iterated $\bowtie$s for $A(R)$, we have:

\begin{theorem}
\label{tricrossedproduct} Let $S$ be an {\bf M3}-monoid  with respect
to a matched pair structure $(S,S)$ and associated map $r_S$. Then the following are equivalent
\begin{enumerate}
\item $r_S$ obeys the YBE.
\item $(S,S\bowtie S)$ form a matched pair with actions extending those of  $(S,S)$.
\item $(S\bowtie S,S)$ forms a matched pair with actions extending those of $(S,S)$.
\end{enumerate}
 In this case the respective actions are
\[ {}^{u.a}v={}^{ua}v,\quad (u.a)^v=u^{{}^av}.a^v;\quad v^{u.a}=v^{ua},\quad {}^v(u.a)={}^vu.{}^{v^u}a\]
for all $u.a\in S\bowtie S,\ v\in S$, and $S\bowtie(S\bowtie S)=(S\bowtie S)\bowtie S$.
\end{theorem}
\begin{proof}  We verify the first matched pair, the second is analogous. The left action here is to view a general element $u.a\in S\bowtie S$ as built from $u,a\in S$ where we multiply them and use the given action of $S$ on $S$ (we use the dot to emphasis the product in $S\bowtie S$): We have
\[ {}^a({}^uv)={}^{au}v=^{\bf M3}{}^{{}^au a^u}v={}^{{}^au.a^u}v\]
so the relations in $S\bowtie S$ are represented, and this requirement determines the action uniquely. We similarly have a unique extension of the actions in $(S,S)$ to an action of $S\bowtie S$, as stated. We then check that these form a matched pair:
\[{}^{u.a}(wv)={}^{ua}(wv)={}^{ua}w\,( {}^{(ua)^w}v)={}^{u.a}w\,( {}^{u^{{}^aw}a^w}v)={}^{ua}w\, ({}^{u^{{}^aw}.a^w}v)={}^{u.a}w\, ({}^{(u.a)^w}v)\]
since the action of $S$ on $S\bowtie S$ has the same structural form as the action  of $S$ on a product in $S$. For the other action in the matched pair,
\begin{eqnarray*}((u.a)(v.b))^w&=& (u{}^av. a^vb)^w=(u{}^av)^{({}^{a^vb}w)}.(a^v b)^w=u^{{}^{{}^av a^v b}w}({}^av)^{({}^{a^v b}w)}. a^{v{}^bw}b^w\\
&=&u^{{}^{avb}w}({}^av.a^v)^{{}^bw}b^w= u^{{}^{avb}w}.(a.v)^{{}^bw}.b^w \\
((u.a)^{{}^{(v.b)}w})(v.b)^w&=&(u.a)^{{}^{vb}w} (v^{{}^bw}.b^w)=u^{{}^{avb}w}.a^{{}^{vb}w}.v^{{}^bw}.b^w \end{eqnarray*}
where $(a.v)^w$ is by definition the action of $w$ on $({}^av.a^v)$. To have equality of these expressions we require $(a.v)^w=a^{{}^vw}.v^w$ for all $w$. Using the unique factorisation in $S\bowtie S$, we need
\[ ({}^av.a^v)^w\equiv ({}^av)^{({}^{a^v}w)}.a^{vw}=a^{{}^vw}.v^w\equiv {}^{(a^{{}^vw})}(v^w).(a^{({}^vw\, {}v^w)})\equiv {}^{(a^{{}^vw})}(v^w).a^{vw}\]
which holds if we assume \textbf{LR3} in $S$. Conversely, this condition is also necessary due to the unique factorisation in $S\bowtie S$. As remarked above, given {\bf M3}, the condition {\bf LR3} is equivalent to $r_S$ obeying the YBE. The other matched pair is similar and requires
the same assumption. That the two matched pairs give the same product on $S\times S\times S$
is a matter of direct computation of the products in the two cases, one readily verifies that they give the same on reducing all expressions to $S\times S\times S$ in the obvious way. 
\end{proof}

 Now, the operation $r_S:S\times S\to S\times S$ expresses reordering of two factors: the value in $S\times S$ read on the left and bottom is transported to the value read on the top and right by $\Rightarrow$ in Figure~1(a) (here $S=T$), an equality $a.u={}^au.a^u$ in $S\bowtie S$. Similarly working within the above triple factorisation and  calling the three copies of $S$ as $S,T,U$ to keep them distinct, the `input' of the Yang-Batxer cube in Figure~1(b) is a reverse-ordered expression in $S\bowtie T\bowtie U$. Each $\Rightarrow$ is a reordering and the `output' of the cube is the canonically ordered expression. Since at each stage the same elements in the triple product are involved,
we have the same result going around the front or around  the back of the cube, i.e. the Yang-Baxter equation for $r_S$ holds. This is a geometric reason for $r_S$ to obey the YBE and provides a different point of view than that in Section~2.

 We are also now in position to provide a full characterisation of the exponentiation problem for $r_S$ under some minimal cancellation assumption. This provides  is a ``monoidal" analogue of a theorem for the group
$G(X,r)$, see   \cite{Lu}, see also \cite{Takeuchi}. Note that we
do not know, in general, under which conditions on $(X,r)$ the monoid $S(X,r)$ is
embedded in $G(X,r)$ so one can not deduce our result from the group case.

\begin{theorem}
\label{S,S matchedth} Let $(X,r)$ be a  2-cancellative
quadratic set, and let $S=S(X,r)$ be the associated monoid,
graded by length. Then $(S,r_S)$ is a graded braided monoid
reducing to $(X,r)$ on degree one \textsl{iff }  $(X,r)$ is  a
braided set.
\end{theorem}
\begin{proof} 
If $(X,r)$ is  a 
braided set, then by Theorem
\ref{theoremA} $S$ is a graded {\bf M3}-monoid with actions extending uniquely the canonical left and right actions associated with $r$, next  by Theorem
\ref{theoremAextra} $(S, r_S)$ is a braided monoid. (In this direction we do not need 2-cancellativeness of $r$).
Conversely, assume that $(S,r_S)$ is a graded braided monoid
reducing to $(X,r)$ on degree one, then, clearly, the hypothesis of  
Lemma~\ref{S,Smatchedlemma} is satisfied, so $(X,r)$ is a braided set.
\end{proof}

We close the section with an open question.
\begin{openquestion}
1) Let $(X,r)$ be an involutive  nondegenerate solution of YBE,
$S=S(X,r)$ the associated YB-monoid. Is it true that $S$ is
cancellative?

 2) Is it true that if $(X,r)$ is a 2-cancellative braided set,
 the associated monoid $S=S(X,r)$ is cancellative?
\end{openquestion}
 In the case when $(X,r)$ is a finite
square-free solution, the answer is affirmative, see \cite{T06}.
In this case the monoid $S(X,r)$ is embedded in $G=G(X,r)$ and the
elements of $G$ are of the shape $u^{-1}v,$ where $u,v\in S.$
Proposition~\ref{cancellsolprop} shows that for arbitrary
2-cancellative solution $S$ satisfies cancellative law on
monomials of length 3.

\section{Matched pair approach to extensions of solutions}
\label{secmpe}

In this section we study extensions of solutions and their
relations with matched pairs of monoids. The notions of \emph{a
union} of solutions, extensions  and one-sided extensions  were introduced in
\cite{ESS}, but only for nondegenerate involutive solutions
$(X,r_X),$ $(Y,r_Y).$ We introduce  what we call \emph{regular
extensions} $(Z,r)$ of arbitrary solutions $(X,r_X),$ $(Y,r_Y)$, and provide
necessary and sufficient conditions (in terms of left and
right actions) for a regular extension $(Z,r)$ to satisfy YBE.  Moreover, regular extensions
obeying the YBE have a very natural interpretation in terms of matched pairs. They correspond to  certain types of strong matched pairings between the associated monoids. 

\subsection{Regular extensions and YB-extensions}

In this section we introduce the notion of regular extensions and provide first results
on when a regular extension obeys the YBE, see Theorem~\ref{propositionBZ}.  We also introduce notations to be used throughout the section.

\begin{definition}
\label{extensiongeneraldef} Let $(X,r_X)$ and
$(Y,r_Y)$ be disjoint quadratic sets (i.e. with bijective maps
$r_X: X\times X \longrightarrow X\times X, \; r_Y: Y\times
Y\longrightarrow Y\times Y$). Let $(Z,r)$ be a set with a
bijection $r: Z\times Z\longrightarrow Z\times Z.$ We say that
$(Z,r)$ is \emph{a (general)  extension of} $(X,r_X),(Y,r_Y),$ if
$Z= X\bigcup Y$ as sets, and  $r$ extends the maps $r_X$ and
$r_Y,$ i.e. $r_{\mid X^2}= r_X $, and $r_{\mid Y^2}=r_Y.$ Clearly
in this case $X, Y$ are $r$-invariant subsets of $Z$. $(Z,r)$ is
\emph{a YB-extension of } $(X,r_X)$, $(Y,r_Y)$ if $r$ obeys YBE.
\end{definition}
\begin{remark}
\label{extensionsrem} In the assumption of the above definition,
suppose $(Z,r)$ is a non-degenerate  extension of
$(X,r_X),(Y,r_Y).$ Then the equalities $r(x,y) = ({}^xy,x^y),$
$r(y,x) = ({}^yx,y^x),$ and the non-degeneracy of $r$, $r_X,$
$r_Y$ imply that
\[
{}^yx, x^y \in X, \;\;\text{and }\;\; {}^xy, y^x \in Y,\;\;
\text{for all}\;\; x \in X, y\in Y.
\]
Therefore, $r$ induces bijective maps
\begin{equation}
\label{rhosigma} \rho: Y\times X \longrightarrow X\times Y ,  \;
\text{and} \;\sigma: X\times Y \longrightarrow Y\times X,
\end{equation}
and left and right ``actions"
\begin{equation}
\label{ractions1} {}^{(\;)}{\bullet}: Y\times X \longrightarrow
X,\;\;\; {\bullet}^{(\;)}: Y\times X \longrightarrow Y,\;
\text{projected from}\; \rho
\end{equation}
\begin{equation}
\label{ractions2}
 \la:
X\times Y \longrightarrow Y,\quad  \ra: X\times Y \longrightarrow
X, \ \text{projected\ from}\; \sigma.
\end{equation}
Clearly,  the 4-tuple of maps $(r_X, r_Y, \rho, \sigma)$ uniquely
determine the extension $r.$ The map $r$ is also uniquely
determined by $r_X$, $r_Y$, and the maps (\ref{ractions1}),
(\ref{ractions2}).

However, if we do not assume $(Z,r)$ non-degenerate, there is no
guarantee that $r$ induces maps as (\ref{rhosigma}), neither
actions (\ref{ractions1}), (\ref{ractions2}).
\end{remark}
\begin{remark}
 Clearly, $(X, r)$ with $ r= \id_{X\times X}$ is a solution, which
 is degenerate whenever $X$ is a set with more than one
element, since ${}^xy=x$ for all $x,y \in X$.

 Given two disjoint solutions $(X,r_X)$, and $(Y,r_Y),$
let $Z= X\bigcup Y$. The following two examples yield the
"easiest" extensions we can get:
\begin{enumerate}
\item \label{abelianextension} \cite{ESS} Define $r: Z\times Z
\longrightarrow Z\times Z$ as $r(x_1,x_2):= r_X(x_1,x_2),$
$x_1,x_2\in X$,
 $r(y_1,y_2):= r_Y(y_1,y_2),$ for all  $y_1,y_2\in Y$ and
\[
r(y,x):=(x,y); \; r(x,y):=(y,x),\; \text{for all}\; x\in X, y \in
Y.
\]
Then $r$ is a solution, it is called \emph{the trivial extension}.
 \item \label{identityextension} Define $r: Z\times Z
\longrightarrow Z\times Z$ as $r(x_1,x_2):= r_X(x_1,x_2),$
for all $x_1,x_2\in X$,
 $r(y_1,y_2):= r_Y(y_1,y_2),$ for all  $y_1,y_2\in Y$ and
\[
r(y,x):=(y,x); \; r(x,y):=(x,y),\; \text{for all}\; x\in X, y \in
Y.
\]
Then,  $r$ is an extension, $r$ is bijective,  but  does not
induce maps (\ref{rhosigma}), nor actions (\ref{ractions1}),
(\ref{ractions2}).
Moreover, $r$ obeys YBE if and only if
$r_X= \id_{X\times X}$, and $r_Y=\id_{Y\times Y}$.

As a very particular example of an YB- extension $(Z,r)$ of two
non-degenerate solutions, which does not induces maps
(\ref{rhosigma}), one can consider the extreme case when $X=
\{x\}, Y = \{y\}$ are one element sets, with the trivial solutions
$r_X= \id_{X\times X}, r_Y= \id_{Y\times Y}.$
\end{enumerate}
\end{remark}
If we want to assure the existence of  maps (\ref{rhosigma}), (we
need them if we want to apply the theory of strong matched pairs),
but not assuming necessarily the bijection $r$ to be
non-degenerate, we should consider only \emph{regular} extensions
$(Z,r)$, which are defined below.
\begin{definition}
\label{regularextensionsdef}
 In notation as above, a (general)  extension $(Z,r)$ of $(X,r_X),(Y,r_Y)$
is  \emph{a regular extension} if $r$ is bijective, and the
restrictions $r_{\mid Y\times X}$ and $r_{\mid X\times Y}$ have
the shape
\[
r_{\mid Y\times X}: Y\times X \longrightarrow X\times Y, \quad
r_{\mid X\times Y} = (r_{\mid Y\times X })^{-1}: X\times Y
\longrightarrow Y\times X.
\]
We call 
\[
\label{groundactions} {}^{Y}{\bullet}: Y\times X \longrightarrow X, \;\;
{\bullet}^{X}: Y\times X \longrightarrow Y\]
\[\la:
X\times Y \longrightarrow Y, \;\; \ra:
  X\times Y \longrightarrow X
\]
projected from $r_{|Y\times X}$ and $r_{|X\times Y}$ the associated ground action and accompanying action respectively. 
\end{definition}

It follows from the definition that each regular extension $(Z,r)$
satisfies
 \[(r\circ r)_{\mid Y\times X}= \id_{\mid Y\times X}, \quad (r\circ r)_{\mid
X\times Y}= \id_{\mid X\times Y},\] but $r$ is not necessarily
involutive on $X\times X,$ neither on $Y\times Y.$
\begin{definition} 
\label{associatedregulardef} With respect to solutions $(X,r_X),$
$(Y,r_Y)$, a pair of maps  
\[{}^{Y}{\bullet}: Y\times X \longrightarrow X, \;\;
{\bullet}^{X}: Y\times X \longrightarrow Y\]
 is called  {\em regular} if the map $r_{|Y\times X}(\alpha,x)=_{\rm def}({}^\alpha x,\alpha^x)$ is invertible. 
\end{definition}
There is clearly a 1-1 correspondence between regular extensions $(Z,r)$ and regular pairs of actions $({}^Y\bullet,\bullet^X)$.  Moreover, regularity of a pair of actions is clearly equivalent to the existence of the accompanying actions such that
\begin{equation}
\label{laraground} {}^{\alpha}x\la \alpha^x=\alpha,\quad
{}^{\alpha} x\ra {\alpha}^x=x,\quad {}^{x\la {\alpha}}(x\ra
\alpha)=x,\quad (x\la \alpha)^{x\ra \alpha}=\alpha.
\end{equation}
Sometimes for simplicity we shall write ${}^x{\alpha}$ instead of
$x\la {\alpha},$  or $x^{\alpha},$ instead of $x\ra {\alpha}.$ 

Henceforth we shall assume that $(X,r_X),$ $(Y,r_Y)$ are arbitrary
disjoint braided sets ($r_X, r_Y$  are bijective maps  obeying YBE), not
necessarily involutive, non-degenerate, or finite. Any additional restriction on the solutions will be mentioned explicitly. We shall consider only regular extensions $(Z,r)$ of
$(X,r_X),$ $(Y,r_Y),$ with corresponding regular ground actions  as in Definition~\ref{associatedregulardef}.

Furthermore, assuming the  actions in Definition~\ref{associatedregulardef} are
given, we also deduce automatically a left action of $\;Y$ on
$X^2$ and a right action of $X$ on $Y^2$ defined as
\[{}^{\alpha}{xy}:={}^{\alpha}x.{}^{{\alpha}^x}y,\quad\quad
(\alpha\beta)^x:={\alpha}^{{}^{\beta}x}{\beta}^x, \quad \text{for
all}\quad x,y \in X, \alpha,\beta \in Y.
\]
For convenience we shall use notation $x,y,z,$ for the elements of
$X,$ $\alpha,\beta, \gamma $ for the elements of $Y.$ Then the following lemma is straightforward:

\begin{lemma}\label{easypropertiesextensions}
In notation as above, let $(Z,r)$ be a regular extension of
$(X,r_X),$ $(Y,r_Y).$ Then 
\begin{enumerate}
\item The set of defining relations $\Re(r)$ of $U=S(Z,r)$ is:
\[
\Re(r)= \{\alpha x= {}^{\alpha}x \alpha^x\mid x \in X, \alpha\in Y \}\bigcup \{x\alpha =(x\la \alpha)(x\ra\alpha)\mid x \in X, \alpha\in Y \}\bigcup\Re(r_X)\bigcup \Re(r_Y).
\]
\item $r$ is 2-cancellative \emph{iff} $r_X,$ and $r_Y$ are
2-cancellative; \item $r$ is involutive \emph{iff} $r_X,$ and $r_Y$ are
involutive; \item $r$ is square-free \emph{iff} $r_X,$ and $r_Y$ are
square-free.
\end{enumerate}
\end{lemma}

\begin{notation}\label{mlnotation} In order to study regular extensions further it will be helpful to have a `local' notation for some of our conditions, in which the specific elements for which the condition is being imposed will be explicitly indicated, indicated in lexicographical order of first appearance. Thus for example {\bf l1(x,y,z)} means the condition exactly as written in Lemma~\ref{ybe} for the specific elements  $x,y,z$. Similarly {\bf r2(x,y,z)} means for the elements $x,y,z$ exactly  as appearing as in Definition~\ref{extendedleftaction}.  Then for example the `local' version of Lemma~\ref{leml1r2} proved in the same way as at the end of the proof there, is the result
\[
r^{12}r^{23}r^{12}(x,y,z)=r^{23}r^{12}r^{23}(x,y,z) \Longleftrightarrow {\rm\bf r1(x,y,z), l2(x,y,z)} \Longleftrightarrow {\rm\bf l1(x,y,z), r2(x,y,z)}
\]
for any fixed triple $(x,y,z)$.  In the present section we consider triples in the set $Z^3$ so for example
\[ {\bf l1(\alpha,x,y)}:\quad {}^\alpha({}^xy)={}^{{}^\alpha x}({}^{\alpha^x}y)\]
where $\alpha,x,y\in Z$. Finally, we use this notation to specify the restrictions of any of our conditions to subsets of interest. For example
\[ {\bf l1(Y,X,X)}:=\{{\bf l1(\alpha,x,y)}\quad {for\ all\ } \alpha\in Y,\ x,y\in X\}.\]
Finally, in view of the following key lemma, we concentrate on the particular examples
\[\begin{array}{lclc}
 {\rm\bf ml1}:={\bf l1(Y,Y,X)}:&{}^{\alpha}{({}^{\beta}x)} =
{}^{{}^{\alpha}{\beta}}{({}^{{\alpha}^{\beta}}x)}&{for\ all\ }x\in X, \alpha,\beta\in Y\\
 {\rm\bf
mr1}:={\bf r1(Y,X,X)}:&({\alpha}^x)^y=({\alpha}^{{}^xy})^{x^y}& {for\ all\ }x,y\in X, \alpha\in Y \\
{\rm\bf \mll}:={\bf l2(Y,X,X)}:&{}^{\alpha}{(r_X(xy))}=r_X({}^{\alpha}{(xy)}) & {for\ all\ }x,y\in X, \alpha\in Y\\
 {\rm\bf \mrr}:={\bf r2(Y,Y,X)}:
&(r_Y(\alpha\beta))^x =r_Y({(\alpha\beta)}^x) &{for\ all\ }x\in X, \alpha,\beta\in Y.
\end{array}\]
  \end{notation}

\begin{lemma} \label{VIPlemma} Let $(Z,r)$ be a regular extension of the quadratic sets  $(X,r_X), (Y,r_Y).$ 

{\bf A.} The following are equivalent:
\begin{enumerate}
\item
\label{l2axy}
{\rm\bf r1(Y,X,X)},  {\rm\bf l2(Y,X,X)}
\item
\label{l2xay}
{\rm\bf r1(X,Y,X)}, {\rm\bf l2(X,Y,X)}
\item
\label{l2xya}
{\rm\bf r1(X,X,Y)},  {\rm\bf l2(X,X,Y)}.
\end{enumerate}

{\bf B.} The following are equivalent:
\begin{enumerate}
\item \label{r2abx}
{\rm\bf l1(Y,Y,X)},  {\rm\bf r2(Y,Y,X)}
\item
\label{r2axb}
{\rm\bf l1(Y,X,Y)}, {\rm\bf r2(Y,X,Y)}
\item
\label{r2xab}
{\rm\bf l1(X,Y,Y)},  {\rm\bf r2(X,Y,Y)}.
\end{enumerate}
\end{lemma}
\begin{proof}  These are all restrictions of the YBE  to different parts of $Z^3$. As explained above, the `local' YBE at $(\alpha,x,y)$, say, is equivalent to {\bf r1}$(\alpha,x,y)$, {\bf l2}$(\alpha,x,y)$ which for all $x,y\in X, \alpha\in Y$ is condition {\bf A}(1). That the three parts of {\bf A} are in fact equivalent uses in an essential way that, by definition, a regular extension is involutive on $X\times Y$ and $Y\times X$. Thus, for example, we look more carefully at the  `local' YBE diagram, 
\begin{equation}
\label{ybediagram5}
\begin{CD}
 \alpha xy  @> r^{23} >>\alpha{({}^xy}x^y)\\
@V  r^{12} VV @VV r^{12} V\\
 ({{}^{\alpha}x}{\alpha}^x)y@.[{}^{\alpha}{({}^xy)}{{\alpha}^{{}^xy}}] x^y\\
@V r^{23} VV @VV r^{23} V\\
{{}^{\alpha}x}[{{}^{{\alpha}^x}y}{({\alpha}^x)^y}] @.
[{}^{\alpha}{({}^xy)}][{}^{{\alpha}^{{}^xy}}{(x^y)}][{({\alpha}^{{}^xy})}^{x^y}]=
w_1\\
  @V r^{12} VV  \\
[{}^{{}^{\alpha}x}{({{}^{{\alpha}^x}y})}][{({}^{\alpha}x)}^{{}^{{\alpha}^x}y}][{({\alpha}^x)^y}]=
w_2@.
 \end{CD}
\end{equation} 
for which $w_1=w_2$ is the condition $r_{12}r_{23}r_{12}(\alpha,x,y)=r_{23}r_{12}r_{23}(\alpha,x,y)$. However, all arrows are bijections so inverting the first $r_{12}$ on the left and the last $r_{23}$ on the right we have that $w_1=w_2$ is equivalent to an instance of the equation $r_{23}^{-1}r_{12}r_{23}=r_{12}r_{23}r_{12}^{-1}$ applied to some element of $X\times Y\times X$. But the inverted instances of $r$ are of mixed type and hence involutive as explained, i.e. $r_{12}^{-1}=r_{12}$ acting the other way, $r_{23}^{-1}=r_{23}$ acting the other way here. Hence if condition A(1) holds then so does A(2). Similarly for all the other parts of the lemma. \end{proof}

\begin{theorem}
\label{propositionBZ} Let $(X,r_X),$ $(Y,r_Y)$ be disjoint solutions,
with  a regular pair of ground actions ${}^Y\bullet,\bullet^X$ , and 
let $(Z,r)$ be  the corresponding regular extension.  Then
$(Z,r)$ obeys YBE {\em if and only if}  {\bf ml1,mr1,ml2,mr2}.
\end{theorem}
\begin{proof}
By hypothesis $r$ is an extension of $r_X$ and $r_Y,$ therefore the `local' YBE already holds on all
$(x,y,z)\in X^3$ and $(\alpha,\beta,\gamma)\in Y^3$. All the other cases are covered in Lemma~\ref{VIPlemma} as explained in the proof of that. Part A there covers the cases where exactly one element is from $Y$ and these hold by the lemma {\em iff} {\bf mr1},\mll. Part B covers the cases where exactly two elements are from $Y$ and these hold {\em iff} {\bf ml1},\mrr. \end{proof}

\begin{corollary}\label{dbraidedset} 
If $(X,r)$ is a braided set, it has a canonical `double braided set' $(X\sqcup X,r_D)$ where
\[ r_D(x,y)=r(x,y),\quad r_D(\bar x,\bar y)=r(\bar x,\bar y),\quad r_D(\bar x, y)=r(\bar x,y),\quad r_D(x,\bar y)=r^{-1}(\bar y,x)\]
for all $x,y\in X$. Here the bar denotes that the element is viewed in the second copy of $X$ and the result of $r$ is correspondingly interpreted.
\end{corollary}
\begin{proof} Here $X=Y$ and $r_X=r_Y$. The ground actions are also those of $X$, hence all the conditions {\bf ml1,mr1,ml2,mr2} reduce to {\bf l1,r1,l2,r2} in $X$. \end{proof}

\begin{remark}
\label{explanml1ml2remark} 
Condition \textbf{ml1 } implies that the assignment $\alpha
\longrightarrow {}^{\alpha}\bullet $ 
extends to a left action of the
associated YB-monoid $S(Y,r_Y)$ and YB-group $G(Y, r_Y)$ on $X$. 
\textbf{mr1} 
assures that the assignment $x \longrightarrow {\bullet}^x $ 
extends to a
right action of the associated YB-monoid or YB-group $(X,r_X)$  on $Y$. Hence {\bf ml1,mr1} are degree 1 versions of axioms {\bf ML1, MR1} of a matched pair. Meanwhile,  \mll\ as a restriction of the invariance condition {\bf l2} to parts 
 of $Z$ asserts that for every $\alpha \in Y$ the two maps $r_X$ and  the left
action ${}^{\alpha}{\bullet}$ (extended on $X^2$ as usual) commute
on $X^2$ so
\[({}^{\alpha}{\bullet})\circ r_X = r_X\circ ({}^{\alpha}{\bullet}): X^2\longrightarrow X^2\]
i.e. that the action ${}^\alpha\bullet$ is compatible with the product of the YB-monoid or YB-group of $(X,r_X)$. Similarly, \mrr\ asserts that for every $x \in X$ there is an equality of maps in
$Y^2$:
\[({\bullet}^x)\circ r_Y = r_Y\circ ({\bullet}^x): Y^2\longrightarrow Y^2\]
or compatibility of $\bullet^x$ with the product of the YB-monoid or group of $(Y,r_Y)$. In this way \mll, \mrr\ can be viewed as the degree one analogue of the conditions {\bf ML2, MR2} of a matched pair. This will be developed further in the next section. 
\end{remark}

\begin{remark} Just as {\bf l2} contains {\bf l1} and {\bf r2} contains {\bf r1} in Lemma~\ref{leml1r2}, the local version of the proof there, applied to parts of $Z$, includes the assertions
\[
\mll\quad\Rightarrow\quad {\bf \mla}:={\bf l1(Y,X,X)}: \quad  {}^{\alpha}{({}^xy)}={}^{{}^{\alpha}x}{({}^{{\alpha}^x}y)},\quad {for\ all\ }x,y\in X, \alpha\in Y
\]
\[
\mrr\quad\Rightarrow\quad {\bf m1ra}:={\bf r1(Y,Y,X)}:\quad ({\alpha}^{\beta})^x= ({\alpha}^{{}^{\beta}x})^{{\beta}^x},\quad{for\ all\ } x\in X, \alpha, \beta\in Y.\]
 \end{remark}

We now show that the under minor assumption of 3-cancellation these weaker conditions together are still sufficient.
 
\begin{theorem}
\label{theoremBcancellative}
Suppose  $(X, r_X),$ and $(Y, r_Y)$ are 2-cancellative and
$\U $ has cancellation on monomials of length 3. Then 
$(Z,r)$
obeys YBE {\em iff}   {\rm\bf ml1, mr1},\mla,\mra. \end{theorem}
\begin{proof} 
Under the hypothesis of the theorem, and Lemma~\ref{easypropertiesextensions}, we know that $(Z,r)$ is 2-cancellative and hence by Lemma~\ref{cancellativelemma} applied to $(Z,r)$ we know that it obeys YBE iff {\bf l1,r1} hold for $Z$. The mixed parts of these are {\bf ml1,mr1},\mla,\mra. The parts of {\bf l1,r1} wholly for elements of $X$ already hold as $(X,r_X)$ is a braided set, and similarly for the parts wholly for elements of $Y$. The next lemma then completes the proof. \end{proof}

By looking in detail at the proof of Lemma~\ref{cancellativelemma} applied to $(Z,r)$, one can see that the proof there is again `local' i.e. applies pointwise and hence to restricted versions of {\bf l1,r1,l2,r2}. In this way one can establish using the methods above:

\begin{lemma}\label{2cancellativeXY} Under the hypothesis of Theorem~\ref{theoremBcancellative}:
\begin{equation*} {\rm\bf mr1, \mla }\Longleftrightarrow
{\rm\bf \mll }; \quad\quad {\rm\bf ml1, \mra }  \Longleftrightarrow {\rm\bf
\mrr }.\end{equation*}
\end{lemma}
Hence under the hypotheses of the theorem one can also say that $(Z,r)$ is a braided set {\em iff} \mll,\mrr. Finally, we look at the simplest types of regular extensions.

 \subsection{Construction of matched pair extensions}

In this section we study the regular extensions $(Z,r)$ of braided sets $(X,r_X)$ and $(Y,r_Y)$ in terms of their
YB-monoids  $S= S(X,r_X),$ $T= S(Y,r_Y)$ respectively. We let 
$U=S(Z,r)$ denote the monoid associated to $(Z,r).$ 
Typically $u,v,w $ will denote elements of $S$ or in $\langle X\rangle$ but if there
is ambiguity we shall indicate exactly which monoid is
considered, similarly, $a,b,c$ will denote elements of $T$
or  $\langle Y \rangle.$ Our main Theorem~\ref{theoremB} answers the question: under what conditions do the ground actions corresponding to the 
regular extension extend to a matched pair $(S,T)$. The difference with Theorem~\ref{propositionBZ} for $(Z,r)$ to obey the YBE is that instead of {\bf ml2,mr2} we require:
\[\begin{array}{lclc}
\mlw:&{\rm for\ all }\  x,y \in X,\; \alpha \in Y&\exists k,\ {}^{\alpha}{(r_X(xy))}=r_X^k({}^{\alpha}{(xy)}) 
\\
\mrw:&{\rm for\ all  }\ x \in X,\; \alpha,\beta \in Y
&\exists k,\ {(r_Y(\alpha\beta))}^x =r^k_Y({(\alpha\beta)}^x). 
\end{array}\]
These equalities  are in $X^2$ in the first case and in $Y^2$ in the second. We recall that two words $xy,x'y'\in\langle X\rangle$ of length two are equal in $S=S(X,r_X)$ {\em iff} 
\[xy=r_X^k(x'y')\quad  \text{
in}\; \langle X \rangle\quad{\rm for\ some\ integer}\ k\]
(this follows from standard Groebner basis results for monoids of this type). So \mlw\  is a weaker version of {\rm\bf \mll } in which equality is only required modulo the relations of $S$. Similarly for \mrw .

\begin{theorem}
\label{theoremB} Let $(X,r_X),$ $(Y,r_Y)$ be disjoint solutions,
with associated monoids $S,$ and $T$. Let ${}^Y\bullet,\bullet^X$ be a regular pair of ground actions.
Then
$(S,T)$ is a graded strong matched pair with actions
extending respectively ${}^{Y}{\bullet}$ and ${\bullet}^X,$ 
{\em if and only if} {\rm\bf ml1, mr1,ml2w,mr2w}.
\end{theorem}
\begin{proof} Assume $(S,T)$ is a graded strong matched pair. Then since the
matching actions are graded, conditions \textbf{ML1},  and
\textbf{MR1}, restricted on $X, Y$ give straightforwardly the
identities \textbf{ml1, mr1 }. Since the left action of $T$ on $S$ is graded and agrees with the
relations on $S$, we conclude that $xy=zt$ in $S$ implies
${}^a{(xy)}={}^a{(zt)}$ in $S,$ for all $a$ in $T.$ Clearly
$xy=r(xy)$ in $S$ for all $x,y\in S.$ Hence \mlw\ comes
straightforwardly from the following equalities in $S^2$
\begin{equation}
\label{leftactr=rleftact} {}^{\alpha}{(r(xy))}=
{}^{\alpha}{(xy)}=r({}^{\alpha}{(xy)})
\end{equation}
for all $x,y \in X, \alpha \in Y.$ Similarly we deduce \mrw. 
It follows then that \textbf{ml1, mr1}, \mlw, \mrw\ 
are necessary conditions.

Next we show that these conditions are also sufficient.
 Hence we now assume that the conditions
\textbf{ml1, mr1}, \mlw, \mrw\ are
satisfied. We follow a strategy similar to the one in Section~3 to
extend the actions on the generating sets to actions  of a strong
matched pair $(S,T)$. Under the assumptions of the theorem we
shall prove several statements. The procedure is parallel to the
one in Section~3 and we omit those proofs that are essentially the
same.

Thus, as a first approximation we extend the ground actions
${}^{Y}{\bullet}$ and ${\bullet}^X,$  to  a left action of $T$
onto $\langle X\rangle$, ${}^T{\bullet}: T\times\langle X\rangle
\longrightarrow \langle X\rangle, $ and a right action of $S$ onto
$\langle Y\rangle,$ ${\bullet}^{S}:\langle Y\rangle \times S
\longrightarrow \langle Y\rangle.$ Note that  \textbf{ml1}
implies straightforwardly a left action of $T$ on $X,$ and
\textbf{ mr1 } implies a right action of $S$ on $Y$ defined the
usual way. Clearly ${}^{(ab)}x= {}^a{({}^bx)},$ and
${\alpha}^{(uv)}=({\alpha}^u)^v$, for all $x\in X,\alpha \in Y,
u,v\in S, a,b \in T.$

\textbf{Step 1.} We define recursively  a left ``action"
\begin{equation}
\label{preactionST1} {}^Y\bullet: Y\times \langle X \rangle
\longrightarrow \langle X \rangle \end{equation}  via
${}^{\alpha}xu:={}^{\alpha}x {}^{{\alpha}^x}u$ (assuming that the
actions ${}^{\alpha}v$ are defined for all $\alpha \in Y$, and all
$v\in X^n$, with $n=|v|$. Analogously we define  a right
action \begin{equation} \label{preactionST2}{\bullet}^X: \langle Y
\rangle\times X \longrightarrow \langle Y \rangle \end{equation}
via ${a\beta}^x := a^{{}^{\beta}x}{\beta}^x.$

To be sure that these actions are well defined we need the
following lemma. It is  verified by an argument similar to the
proof of Lemma~\ref{lemML1}. We consider equalities of monomials
in the free monoids $\langle X \rangle,$ and $\langle Y
\rangle.$
\begin{lemma}
\label{ML1ST} The following conditions hold: a) {\bf ML2} holds
for $Y$ ``acting" (on the left) on $\langle X\rangle,$ that is:
\[ {}^{\alpha}{(uv)}=
({}^{\alpha}u).({}^{{\alpha}^u}v)\;\; \text{for all}\;\; \alpha
\in Y, u,v \in \langle X\rangle,\; \; |u|,|v| \geq 1.\]
b) {\bf MR2} holds for $X$ ``acting" (on the right) on $\langle
Y\rangle,$ that is:
\[(ab)^x= (a^{{}^bx})b^x\;\; \text{for all}\;\; x
\in X, \; a,b \in \langle Y\rangle,\; \; |a|,|b| \geq 1.\]
\end{lemma}

So  the left and the right ``actions" (\ref{preactionST1}),
(\ref{preactionST2}) are well defined.

\textbf{Step 2.} We extend the ``actions" (\ref{preactionST1}),
(\ref{preactionST2}) to a left action of $T$ onto $\langle
X\rangle$,  and a right action of $S$ onto $\langle Y\rangle.$
\begin{lemma}
\label{ML1&MR1ST} The actions ${}^{\alpha}{\bullet},$ on $\langle
X\rangle$ and $\bullet^x$ on $\langle Y\rangle$ extend to
${}^a\bullet,$ and $\bullet^u$ for arbitrary monomials $a\in T,
u\in S,$ that is to left and right actions
$
{}^T{\bullet}: T\times \langle X\rangle \longrightarrow \langle
X\rangle$ and ${\bullet}^S: \langle Y\rangle \times
S \longrightarrow \langle Y\rangle$ 
obeying
\[  {}^{ab}u= {}^a{({}^bu)} \;\text{is an equality
in}\;\langle X\rangle,\;\text{for all}\; u \in \langle X\rangle,
a, b\in T\]
\[ a^{uv}= (a^u)^v\;\text{is an equality
in}\;\langle Y\rangle,\; \;\text{for all}\;  a \in \langle
Y\rangle, u,v \in S.\]
\end{lemma}
\begin{proof}
It will be enough to show that the following equalities hold:
\begin{equation}
\label{L11ST} {}^{\alpha}{({}^{\beta}u)} =
{}^{{}^{\alpha}{\beta}}{({}^{{\alpha}^{\beta}}u)}\; \text{for all}
\; \alpha,\beta \in Y, u \in \langle X\rangle.
\end{equation}
\begin{equation}
\label{R11ST}
(a^x)^y = (a^{{}^xy})^{x^y} \; \text{for all} \; x,y \in X, a \in
\langle Y\rangle.
\end{equation}

(\ref{L11ST}) is proven by induction on $|u|=n.$ By
hypothesis (\ref{L11ST}) holds for $|u| =1$, which gives
the base for induction. Assume  (\ref{L11ST}) is true for all $u
\in X^n.$ Let $u \in X^{n+1}.$ Then $u = tv,$ with $t \in X$, $v
\in X^n.$ Consider the equalities:
\[
{}^{\alpha}{({}^{\beta}u)}={}^{\alpha}{({}^{\beta}tv)}
={}^{\alpha}{(({}^{\beta}t).({}^{{\beta}^t}v))}=
({}^{\alpha}{({}^{\beta}t)})({}^{{\alpha}^{{}^{\beta}t}}{({}^{{\beta}^t}v)})
\]
where the second and last are by Lemma~\ref{ML1ST}a), and
\[ =^{\text{inductive
ass.}}\;({}^{\alpha}{({}^{\beta}t)})({}^{(({\alpha}^{{}^{\beta}t}).({\beta}^t))}v)=
({}^{\alpha}{({}^{\beta}t)})({}^{(\alpha\beta)^t}v)\]
using Lemma~\ref{ML1ST}b) for the last equality. 
We have shown
\begin{equation}
\label{eqs4ST} {}^\alpha{({}^{\beta}u)=
({}^{\alpha\beta}t})({}^{(\alpha\beta)^t}v)= w_1.
\end{equation}
 Similarly,  we obtain:
\begin{equation}
\label{eqs5ST} {}^{{}^{\alpha}{\beta}}{({}^{{\alpha}^{\beta}}u)}
 = ({}^{{}^{\alpha}{\beta}}{({}^{{\alpha}^{\beta}}t)})({}^{({}^{\alpha}{\beta}.{\alpha}^{\beta})^t}v)=w_2.
\end{equation}
Now
\[
{}^{{}^{\alpha}{\beta}}{({}^{{\alpha}^{\beta}}t)}=^{\bf ml1}{}^{\alpha\beta}t,\quad 
({}^{\alpha}{\beta}.{\alpha}^{\beta})^t=^\mrw{(\alpha\beta)}^t .
\]
Therefore, the last two equalities and (\ref{eqs4ST})-(\ref{eqs5ST}) imply
${}^{{}^{\alpha}{\beta}}{({}^{{\alpha}^{\beta}}u)}={}^{\alpha}{({}^{\beta}u)}.$
  This verifies (\ref{L11ST}).

  Now  the equality
\[
{}^{(\alpha_1...\alpha_k)}u :=
 {}^{\alpha _1}{(...({}^{\alpha_ {k-1}} {({}^{\alpha _k}u)} )..)}
\]
gives a well defined left action $T\times \langle X \rangle
\longrightarrow \langle X \rangle$ , which implies the first part of the lemma.

The proof of (\ref{R11ST}) is analogous. So the right action
$\langle Y \rangle \times S  \longrightarrow \langle Y \rangle$ is
also well defined, and this implies the second part of the lemma.
\end{proof}

The following lemma is analogous to Proposition~\ref{propML2} and
is proven by similar argument.

\begin{lemma}  \label{propML2ST}
\[ 
{}^a{(uv)}=({}^au).({}^{(a^u)}v), \;\text{is an equality in}\;
\langle X \rangle \;\text{for all}\; a \in T,  u, v\in \langle X
\rangle .
\]
\[
 (ab)^u=(a^{{}^bu}).(b^u) \;\text{is an equality in}\; \langle Y \rangle \;\;\text{for all}\;  a, b\in \langle
Y \rangle, u \in S.
\]
\end{lemma}
So far we have extended the actions ${}^YX,$ and  $Y^X$ to actions
${}^T{\langle X\rangle}$ and ${\langle Y\rangle}^S.$

\textbf{Step 3.} We will show next that these actions respect the
relations in $S,$ $T$, and therefore they induce naturally left
and right actions ${}^{(\;)}{\bullet}: T\times S \longrightarrow S
$, and ${\bullet}^{(\;)}: T\times S \longrightarrow T $. We need
analogues of Lemma~\ref{beautifullemma}, and Proposition
\ref{agreementwithreductions}. For these analogues we use slightly
different arguments.
\begin{lemma}\label{MLST2}
Condition \mlw\ is equivalent to the following
\[ xy= zt \;\text{in }\; S^2 \Longrightarrow {}^{\alpha}{(xy)}= {}^{\alpha}{(zt)}\;\text{in }\; S^2
 \;\text{for all} \; x,y,z,t \in X, \alpha \in
Y.
\]
\end{lemma}
\begin{proof}
Since
clearly $r(xy)=xy$ in $S$, so if we assume the stated condition in the lemma, 
${}^{\alpha}{(r(xy))}={}^{\alpha}{xy}$ in $S,$ which together with
the evident equality   ${}^{\alpha}{xy}= r({}^{\alpha}{xy})$ in
$S$, imply ${}^{\alpha}{(r(xy))}= r({}^{\alpha}{xy}).$ 
Conversely,
assume \mlw\ holds, and let $xy=zt$ in $S^2.$ Remind that
all equalities in $S^2$ come straightforwardly from the defining
relations $\Re(r),$ therefore  there exists positive integers
$k,l$ and a monomial $tu  \in X^2 $ such that $xy = r^l(tu),$ $zt = r^k(tu),$  (We do not
assume $r$ is necessarily of finite order.) 
Suppose now $\alpha \in Y.$ The following
equalities hold in $S$: 
\[{}^{\alpha}{(zt)}=
{}^{\alpha}{(r^k(tu))}=^{\mlw}
\;r({}^{\alpha}{(r^{k-1}(tu))})= \cdots =
r^k{({}^{\alpha}{(tu)})}={}^{\alpha}{(tu)},
\]
Similarly, ${}^{\alpha}{(xy)}={}^{\alpha}{(tu)}$
hence \mlw\   implies the condition stated in the lemma.
\end{proof}
\begin{lemma} \label{bl2}\

\begin{enumerate}
\item [\text{(i)}]  $ xy= zt \;\text{in }\; S^2 \Longrightarrow
{}^a{(xy)}= {}^a{(zt)} \;\text{in }\; S^2, \; \text{for all} \;
  x,y,z,t \in X, a \in T.$
\item [\text{(ii)}] $\alpha\beta= \gamma\delta \;\text{in }\;
T^2 \Longrightarrow (\alpha\beta)^u= (\gamma\delta)^u \;\text{in
}\; T^2, \; \text{for all} \; \alpha,\beta,\gamma,\delta \in Y, u
\in S.$
\end{enumerate}
\end{lemma}
\begin{proof}
We shall prove (i), the proof of (ii) is analogous. We use
induction on $|a|$ with Lemma~\ref{MLST2} as base for the
induction. Assume (i) is true for all $a$, with $|a|=n.$
Consider:
\begin{eqnarray*}
{}^{a\beta}{(xy)}&=&{}^a{({}^{\beta}{(xy)})}\quad :\text{by
{\bf ML1}}\\
&=& {}^a{({}^{\beta}{(zt)})} \quad : \text{by} \;Lemma~\ref{MLST2},\;
\text{and the inductive assumption}\\
&=&{}^{a\beta}{(zt)} \quad :\text{by \textbf{ML1}}.
\end{eqnarray*}
\end{proof}
\begin{corollary}
\label{bcor} The actions ${}^{(\;)}{\bullet}: T\times
S\longrightarrow S$, and ${\bullet}^{(\;)}: T\times S\longrightarrow
T$ are well defined.
\end{corollary}
\begin{proof}
It will be enough to verify:
\[
xy =zt \;\text{is an equality in}\; S^2 \Longrightarrow
{}^a{uxyv}={}^a{uztv} \; \text{is an equality in}\; S
\]
for all $x,y, z,t \in X,u,v \in S, a \in T.$

By Lemma~\ref{propML2ST} the following are equalities in $\langle
X\rangle$:
\begin{equation}
\label{bl3} {}^a{uxyv}=({}^au)({}^{a^u}{xy})({}^{a^{uxy}}v), \quad
{}^a{uztv}=({}^au)({}^{a^u}{zt})({}^{a^{uzt}}v).
\end{equation}
Now the equality $xy=zt$ in $S$ implies ${}^{a^u}{xy} =
{}^{a^u}{zt}$ in $S$ (by Lemma~\ref{ML1&MR1ST}), and
$a^{(uxy)}=a^{(uzt)}$ in $\langle Y\rangle$ (by \textbf{MR1}), so
replacing these in (\ref{bl3}) we obtain ${}^a{uxyv}={}^a{uztv}$
holds in $S.$ It follows then that for any $w_1, w_2\in S,$ and
any $a \in T$, one has
\[
w_1=w_2\; \text{holds in} \;S \Longrightarrow {}^a{w_1}={}^a{w_2}\ 
\text{holds in} \;S .
\]
Hence the left action of $T$ on $S$ is well defined. Similar
argument verifies that the right action of $S$ on $T$ is also well
defined.
\end{proof}
It follows from  Lemmas~\ref{ML1&MR1ST},~\ref{propML2ST} that the
actions obey {\rm\bf ML1, MR1, ML2, MR2}. (Clearly, an equality of
words $u=v$ in $\langle X\rangle$  implies $u=v$ as elements of
$S$, and an equality of words $a=b$ in $\langle Y \rangle$ implies
$a=b$ as elements of $T$). We have proved  Theorem~\ref{theoremB}. 
\end{proof}

\begin{proposition}
\label{MPcancellative}
Suppose  $(X, r_X),$ and $(Y, r_Y)$ are 2-cancellative and
$\U $ has cancellation on monomials of length 3. Then 
$S,T$ is a matched pair {\em iff}  {\rm\bf ml1, mr1} hold for all $x,y \in X, \alpha,\beta \in Y.$ .
\end{proposition}
\begin{proof} 
 Under the hypothesis of the propositon, we first prove that
\begin{equation}\label{ml1ml2w}  {\bf ml1} \Longleftrightarrow \mrw;
\quad\quad {\bf mr1}  \Longleftrightarrow \mlw.\end{equation}
Thus, let $x, y \in X,$ $\alpha \in Y.$ We consider again the diagram
(\ref{ybediagram5}).
As we know each two monomials in the diagram are equal as elements
of $U= S(Z,r).$ Hence $w_1= w_2$ in $U,$ which as before we can write
as \begin{equation}
\label{ar=ra2}
w_1=[{}^{\alpha}{(r_X(xy))}][{({\alpha}^{{}^xy})}^{x^y}]=[r_X({}^{\alpha}{(xy)})][{({\alpha}^x)^y}]=w_2
\end{equation}
as an equality in ${\U }^3.$ By assumption $\U $ has cancellation
on monomials of length 3, so (\ref{ar=ra2}) yields the second of (\ref{ml1ml2w}). 
Analogous argument for a diagram starting with $\alpha\beta x$
gives the first of  (\ref{ml1ml2w}). We then use Theorem~\ref{theoremB}. 
\end{proof}

\begin{proposition}
\label{U=STml1mr1}
Suppose  $(X, r_X),$ and $(Y, r_Y)$ are 2-cancellative, $S,T$ is a matched pair  and 
$\U = S\bowtie T.$ Then $(Z,r)$ is a solution \emph{iff}  {\bf \mla, \mra}.

In this case $S,T$ is a strong matched pair.
\end{proposition}
\begin{proof} 
By Theorem~\ref{propositionBZ} 
$(Z,r)$ obeys YBE {\em if and only if}  {\bf ml1,mr1,ml2,mr2}.
By hypothesis $S,T$ is a matched pair, so {\bf ml1, mr1} hold.
It will be enough then to prove that under our assumptions one has: 
\[
{\bf \mla} \Longleftrightarrow {\bf ml2}; \quad {\bf \mra} \Longleftrightarrow {\bf mr2},
\]
for all $x,y \in X, \alpha,\beta \in Y.$

We will show ${\bf \mla} \Longleftrightarrow {\bf ml2}$. Let 
$x,y \in X, \alpha \in Y.$
Look again  at the monomials $w_1,w_2$ in the the diagram
(\ref{ybediagram5}). We know that 
\begin{equation}
\label{w1=w2}
w_1=[{}^{\alpha}{({}^xy)}][{}^{{\alpha}^{{}^xy}}{(x^y)}][{({\alpha}^{{}^xy})}^{x^y}]=[{}^{{}^{\alpha}x}{({{}^{{\alpha}^x}y})}][{({}^{\alpha}x)}^{{}^{{\alpha}^x}y}][{({\alpha}^x)^y}]=
w_2
\end{equation}
is an equality in ${\U }^3.$ By assumption $\U= S\bowtie T,$ 
hence (\ref{w1=w2}) implies the following equality in $S^2:$
\begin{equation}
\label{w1=w2eq2}
[{}^{\alpha}{({}^xy)}][{}^{{\alpha}^{{}^xy}}{(x^y)}]=[{}^{{}^{\alpha}x}{({{}^{{\alpha}^x}y})}][{({}^{\alpha}x)}^{{}^{{\alpha}^x}y}],
\end{equation}
furthermore 
\begin{equation}
\label{w1=w2eq3}
[{}^{\alpha}{({}^xy)}][{}^{{\alpha}^{{}^xy}}{(x^y)}]= {}^{\alpha}{(r_X(xy))}, \quad [{}^{{}^{\alpha}x}{({{}^{{\alpha}^x}y})}][{({}^{\alpha}x)}^{{}^{{\alpha}^x}y}]= r_X({}^{\alpha}{(xy)}),
\end{equation}
see for example (\ref{rx}).
Then the 2-cancellativeness of $r_X,$ (\ref{w1=w2eq2}), and (\ref{w1=w2eq3}) give the desired implications. 
Analogous argument proves ${\bf \mra} \Longleftrightarrow {\bf mr2}$.
\end{proof} 

It is possible to prove further results along these lines \cite{TGISM}. For example, one can
show  under the same hypotheses as the above proposition, if the left and right ground actions are nondegenerate then  $\U $ has  cancellation on monomials of length 3.

\subsection{Matched pair characterisation of regular YB-extensions}

It is clear comparing Theorem~\ref{propositionBZ} and  Theorem~\ref{theoremB} that regular  extensions of $(X,r_X)$ and $(Y,r_Y)$ obey the YBE require that the associated monoids $S=S(X,r_Y)$ and $T=S(Y,r_Y)$ form a matched pair but that the latter is a weaker assertion. In this section we give a further matched pair requirement that then gives an exact characterisation of when a regular extension obeys the YBE, see Theorem~\ref{bellaT2}. We use the conventions above for $S,T$ and as above we let $U=S(Z,r)$.

First of all, let $(Z,r)$ be an YB-extension of $(X,r_X),$
$(Y,r_Y).$  Then  $(\U ,
r_{\U })$ is a braided monoid, induced from $(Z,r)$, see Theorem
\ref{theoremAextra}.  As we have seen in Section~3, the
"ground" left and right actions ${}^{(\;)}{\bullet}: Z \times Z
\longrightarrow Z, \; {\bullet}^{(\;)}: Z\times Z \longrightarrow
Z$ induced by $r$ extend uniquely to left and right actions
\[
{}^{(\;)}{\bullet}: \U  \times \U  \longrightarrow \U , \;
{\bullet}^{(\;)}: \U  \times \U  \longrightarrow \U ,
\]
which respect the grading and make $(\U ,\U )$ a graded strong
matched pair, with  condition {\bf M3}, see Theorem
\ref{theoremA}.
(Note that its proof in this direction does not need the
assumption that $r$ is 2-cancellative). It is not difficult to see
that these induce left and right actions:
\begin{equation}\label{STactions1}
{}^{T}{\bullet}: T\times S \longrightarrow S, \;\;\text{and}\;\;
{\bullet}^{S}: T\times S \longrightarrow T,
\end{equation}
and accompanying actions
\begin{equation}\label{STactions2}
\la: S\times T \longrightarrow T, \;\;\text{and}\;\;
  S \longleftarrow S\times T:  \ra
\end{equation}
which  make $(S,T)$ a graded strong matched pair. Hence the
following lemma is true.
\begin{lemma}
\label{ST}
 Let $(X,r_X),$ $(Y,r_Y)$ be disjoint solutions, with YB-monoids, respectively $S$ and $T$.
  Let $(Z,r)$ be a regular YB-
   extension of $(X,r_X),$ $(Y,r_Y),$ with a  YB-monoid $\U $. Then $(S,T),$  is a  graded strong
   matched pair with  actions  (\ref{STactions1}) and (\ref{STactions2}) induced from the
   braided monoid $(\U, r_{\U})$.
\end{lemma}
We will show first that in the case when $(Z, r)$ is an YB- extension its associated group $G(Z,r)$ is a double crossed product of 
$G(X, r_X)$ and $G(Y,r_Y).$
\begin{proposition}
\label{BP}
 Let $(X,r_X),$ $(Y,r_Y)$ be disjoint solutions, with 
YB- groups $G_X=G(X,r_X),$ and $G_Y=G(Y,r_Y)$.
Let $(Z,r)$ be a regular YB-
   extension of $(X,r_X),$ $(Y,r_Y),$ with a  YB-group $G_Z=G(Z, r)$. 
Then
\begin{enumerate}
\item \label{BP1} $G_X, G_Y$ is a matched pair of groups with
actions induced
 from  the braided group $(G_Z, r).$
 \item
 \label{BP2} $G_Z$ is isomorphic to the double crossed product $G_X\bowtie
 G_Y$.
\end{enumerate}
\end{proposition}
\begin{proof}
(\ref{BP1}) is straightforward. To prove 
(\ref{BP2}) we first show that $G_Z=G_X.G_Y$  as set, i.e. that every element 
 $w \in
G_Z$ can be presented as $w=ua, u \in G_X, a\in G_Y.$ Indeed,
$G_Z$ is generated by $Z= X\bigcup Y,$
so every element $w \in G_Z, w \neq 1$ has the shape
\begin{equation}
\label{w} w= u_1a_1u_2a_2\cdots u_{s+1}a_{s+1},
\end{equation}
where $u_i \in G_X, a_i \in G_Y, i = 1,... s,$ $u_1=1$, or
$a_s=1,$ is possible, but if $s \geq 0,$ one has $u_i\neq 1$, for $2
\leq i\leq s+1$, and $ a_j \neq 1,$ for $1 \leq j \leq s.$ We
have to "normalize" $w.$ For the purpose we use induction on the number $n(w)= s$
of sub-words $a_i$ of $w$ which occur in a "wrong" place, and
have to be "moved" to the right. There is nothing to prove if $s
=0.$ This is the base for the induction. Suppose every $w$ with
$n(w) \leq s-1$ can be reduced to a presentation $w = ua,$
with $u \in G_X,$ $a \in G_Y.$ Let $w\in G_Z$ has a
presentation (\ref{w}) with  $ n(w) = s, s \geq 1.$  We reduce $w$ using condition {\bf M3}. 
\begin{equation}
\label{w1} w=u_1a_1u_2a_2\cdots u_sa_s u_{s+1}a_{s+1}
=^{{\bf M3}} u_1a_1u_2a_2 \cdots
u_s({}^{a_s}{u_{s+1}})({a_s}^{u_{s+1}}).a_{s+1} =
w_1
\end{equation}
Now (\ref{w1}) gives a new presentation $w_1$ of $w$, $w=w_1$ as elements of $G_Z$, but  $n(w_1)
= s-1,$  so by the inductive
assumption $w_1$ can be reduced to the shape $w_1= ua,$ with
$u \in G_X,$ $a \in G_Y.$ Hence every element of $G_Z$ can be presented in the required normal form.

Moreover, it is easy to see  that for each $w \in G_Z$ the presentation $w
= ua,$ with $u\in G_X, a\in G_Y$ is unique. Indeed, if we assume,
\[ua = vb, u,v \in G_X, a, b \in G_Y\]
then
\[
v^{-1}u = b{a}^{-1} \in G_X\bigcap G_Y.
\]
By hypothesis $X$ and $Y$ are disjoint sets, so $G_X\bigcap
G_Y=1,$ which in view of the previous equality implies $v=u,
b = a.$ This proves part (\ref{BP2}).
\end{proof}

We will show analogous result about monoids. In this case we do not assume cacellation holds in $U.$ 

\begin{lemma} 
\label{mainlemma}
Let $(Z,r)$ be a YB-extension in the setting above. Each
element $w \in \U$ can be presented uniquely as 
$w = ua, u \in S,
a \in T.$ 
\end{lemma}
\begin{proof}
We will make first some general remarks.
Let $W_0 \in \langle Z \rangle.$ Clearly $W_0$ is a product of elemets of  $Z = X\bigcup Y.$ Suppose $|W_0| = N.$
By hypothesis $(Z,r)$ is a solution, therefore the \emph{braid group}
$B_N$ acts on $Z^N.$ Note that the orbit of $W_0$ under this action
consists  exactly of all words $W \in  \langle Z \rangle,$ with the property $W_0 =W$ as elements of $\U$. Indeed,
$W_1=W_2$ in $\U$ \emph{iff } there are finite sequences of
replacements $\rho_i, 1 \leq i \leq p, \sigma_j, 1 \leq j \leq q,$
coming from the relations $\Re(r)$, and such that the following is
an equality of words in the free monoid $\langle Z\rangle:$
\begin{equation}
\label{replacementseq}
\rho_p\circ\cdots \rho_1(W_1)= \sigma_q\circ\cdots
\circ\sigma_1(W_2).
\end{equation}
But each such a replacement is exactly   $r^{ii+1}$, for an appropriate $i, 1 \leq i \leq N-1.$ 
($r^{ii+1} = \id_{Z^{i-1}}\times r\times \id_{Z^{N-i-1}} ).$ Now consider $\deg_{X}$ (respectively $\deg_{Y}$) which count the 
number of symbols from $X,$ (respectively, from $Y$) that occur in an element of $\langle Z\rangle$. Since each relation in $\Re(r)$ has the shape $yz = {}^yz.y^z,$ where $y,z \in Z$ it is immediate that each replacement in (\ref{replacementseq}) leaves the degrees unchanged. Hence both $\deg_X$ and $\deg_Y$ are defined on $U$ itself and in particular are independent of the choice of representative $W_0$, say,  of our given element $w$.  We let $K=\deg_X(w)$ and $L=\deg_Y(w)$.

Since $U$ is an {\bf M3} monoid we can clearly use similar  replacements as in the proof of Proposition~
\ref{BP} to put $w$ in a normal form $w=ua$ with corresponding presentatation in normal form $W_0=U_0A_0$, where $U_0\in X^K$, $A_0\in Y^L$. The following sub-lemma shows that any element of $B_N$ that sends $W_0$ to something again of this form has the same action as an element of  $B_K\times B_L$. Hence any other presentation in normal form  has the shape $\rho_1(U_0)\rho_2(A_0)$ where $\rho_1\in B_K$ and $\rho_2\in L$. But by the remarks above applied to $X$ (respectively $Y$) we see that $U_0=\rho_1(U_0)$ when viewed in $S$ and $A_0=\rho_2(A_0)$ when viewed in $T$. \end{proof}

 \begin{lemma} \label{braidlemma} Let $N=K+L\ge 2$ be positive integers and consider the action of the braid group $B_N$ on $Z^N$ induced by a regular YB-extension  $(Z,r)$ of  $(X,r_X)$, $(Y,r_Y)$.   The action of any braid group element sending $X^K\times Y^L\to X^K\times X^L$ can be presented as the action of an element of $B_K\times B_L$. 
 \end{lemma}
 \begin{proof} A regular extension means that the restrictions $r_{X,X}=r_X$, $r_{Y,Y}=r_Y$, $r_{Y,X}$ and $r_{X,Y}=r_{Y,X}^{-1}$ obey the mixed braided relations where strands are labelled by $X$ or $Y$ (i.e. there is a braided category for which objects are arbitrary products of $X,Y$ and the braiding is given by the various restrictions of $r$). We present an element of $b\in B_N$ as a series of elementary crossings and we represent it as a map $X^K\times Y^L\to X^K\times Y^L$ defined as the corresponding composition of $r_{X,X},r_{Y,Y},r_{X,Y},r_{Y,X}$ or their inverses according to the orientation and  labelling of strands as at each crossing. A standard crossing say with a left strand passing down and over a right strand is represented by one of the $r$, while a reversed crossing where the left strand passes down and under is represented by one of the $r^{-1}$. The most important fact for us is that in this representation an  $X-Y$ crossing is represented by the same map for either orientation of crossing. With these preliminaries, the proof of the lemma proceeds by induction on the total number $M$ of $X-X$ or $Y-Y$ crossings in a presentation of $b$. 
 
 \begin{figure}
\[ \includegraphics{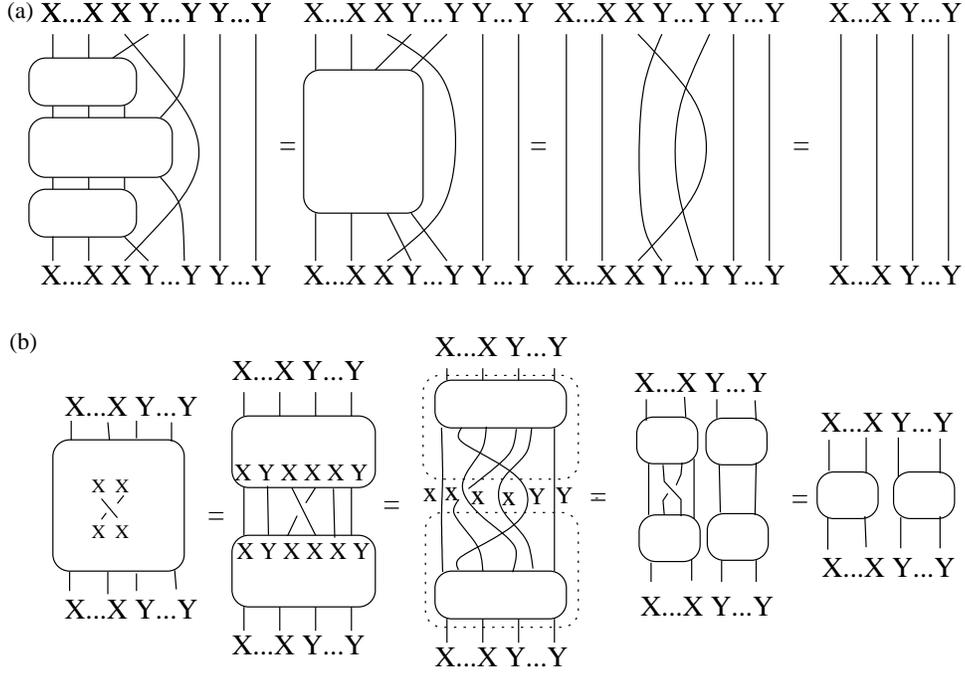}\]
\caption{Diagrams in proof of Lemma~\ref{braidlemma} needed for  $U=S\bowtie T$, (a) case $M=0$ and (b) case $M\ge 1$. }
\end{figure}

 {\bf Case $M=0$}. This means that the only crossings are of $X-Y$ or $Y-X$ type. We are going to do an induction in the even number $P$ of such crossings to prove that the representation of such a braid is the identity. If $P=0$ we have the identity braid and are done. By the remarks above the representation does not depend on whether the crossings are `under' or `over'.   Our notation in Figure~2 therefore does not distinguish these. Moreover, if a strand crosses another and then immediately crosses back as the braid presentation is read from top down, the representation is the same as if the two crossings were replaced by identity maps. This is then a presentation of a braid of the same type but with $P-2$ crossings, and we assume the result is true for all such braid presentations as induction hypothesis. Therefore it suffices to prove the result for braid presentations with $M=0$ and no `double crosses' in which a strand is crossed and then crossed back. In this case consider the right-most $X$ strand in $X^K\times Y^L$. It cannot cross to the left as this would need an $X-X$ crossing but it can cross to the right. However  there are only a finite number of $Y$' to the right so it cannot keep crossing to the right. At some point the $X$ strand must cross over some  $Y$ strand in maximal position $L'$ to the right followed at some point later by crossing back to the left. By the same argument of no double-crossings, it must then keep crossing $Y$'s to the left until it eventually crosses back over all the $Y$ strands. This is depicted in Figure~2(a) where the boxes depict unkown braid of $X-Y$ and $Y-X$ type. Since operations in disjoint strands commute we can group all these boxes as one large box which is then a presentation of a braid on $X^{K-1}\times Y^{L'}$ and less than $P$ crossings. Hence by our induction hypothesis the representation of this unknown braid is trivial. Hence our original  braid has the same representation as one of the form shown in the middle in Figure~2(a) which is then trivial for reasons as above. An informal explanation of the idea behind the above argument is that $r$ is involutive for $X-Y$ and $Y-X$ crossings and since only these are involved, the representation is ``effectively''  as for a representation of the symmetric group $S_N$, which would be trivial here.
 
{\bf Case $M\ge 1$}. Mark an $X-X$ or $Y-Y$ crossing according to what is available, for example we illustrate the former case. We divide the braid presentation into two boxes, those that came before and those that came after the marked crossing. In each case we use only $X-Y$ and $Y-X$ crossings to move all the $Y$ strands to the right after the first box, and before the second box. This is shown in Figure~2(b). In each case the composites shown dashed are braid presentations in our standard form on $X^K\times Y^L$ but with smaller total number of $X-X$ and $Y-Y$ crossings. By our induction hypothesis, these composites have the same representation as the product of some braid on $X^K$ and some braid on $Y^L$. This is the middle equality in Figure~2(b). Hence the original braid has representation in the same required  form. Analogous arguments apply if the marked crossing is a $Y-Y$ one. \end{proof}

We have done most of the work in these lemmas for one direction of the following theorem.

\begin{theorem}
\label{bellaT2}
Let $(X,r_X),$ $(Y,r_Y)$ be disjoint solutions,
with 
associated monoids $S,$ and $T$. Let ${}^Y\bullet,\bullet^X$ be
a regular pair of ground actions, $(Z,r)$   the corresponding 
regular extension and $\U$ its associated monoid.  Then
$(Z,r)$ obeys YBE {\em if and only if}  the following condition hold:
\begin{enumerate}
\item 
\label{it1} 
$U =  S\bowtie T$, where $(S,T)$ is a strong matched pair, with actions extending the ground actions; and
\item
\label{it2} 
$(S,U)$ and $(U,T)$ are matched pairs with canonical actions induced from the actions of $S,T$.
\end{enumerate}

In this case $G_Z= G_X\bowtie G_Y$.
\end{theorem}
\begin{proof}
Assume first that  $(Z,r)$ is  a YB-regular extension
of $(X,r_X),$ $(Y,r_Y).$ 
By Lemma~\ref{ST}
$(S,T),$  is a  graded strong
   matched pair with  actions  (\ref{STactions1}) and (\ref{STactions2}) induced from the
braided monoid $(\U, r_{\U})$. The last lemma shows that $U =  S\bowtie T$,  which is part (1). For (2) we use 
again the fact that 
the "ground" 
left and right actions extend uniquely to left and right actions which make $(\U ,\U )$ a graded strong
matched pair, see Theorem
\ref{theoremA}.
These actions immediately restrict to actions
\[
{}^{(\;)}{\bullet}: U  \times S  \longrightarrow \U , \quad 
{\bullet}^{(\;)}: U  \times S  \longrightarrow S ,
\] which make $(S,U )$ a graded 
matched pair. Similarly, $(U,T)$ is a matched pair with actions restricted from those of from $(U,U).$

Conversely, assume conditions (\ref{it1}), (\ref{it2}) are satisfied. We will show
that $(Z,r)$ obeys YBE. By Theorem~\ref{propositionBZ} it will be enough to prove conditions
{\bf mr1, \mll , ml1, \mrr }.  By assumption $(S,\U)$ is a matched pair, so the action of $S$ respects the relations of $U$. Using {\bf MR2} from the matched pair axioms we have in particular that
\begin{equation}\label{SUproof}({}^\alpha y)^{({}^{\alpha^y}x)}.(\alpha^y)^x =({}^\alpha y.\alpha^y)^x=(\alpha y)^x=\alpha^{{}^yx}.y^x=
{}^{\alpha^{({}^yx)}}(y^x).(\alpha^{{}^yx})^{y^x}\end{equation}
in $U$ for all $x,y\in X, \alpha\in Y$ (this can also be written as {\bf mr2w} with $Z$ in the role of $Y$ for construction of the matched pair $(S,U)$).  By assumption part  (1) holds, so $\U$ decomposes into $ST$ and each element $w \in \U$ has unique presentation as $w=ua, u\in S, a\in T.$ Hence equality of the two factors of (\ref{SUproof}) holds separately, which we  identify as ${\bf lr3(\alpha,y,x)}$ and ${\bf r1(\alpha,y,x)}$. Also, since $U$ acts on $S$, its relations are respected which means ${}^{\alpha}({}^yx)={}^{{}^\alpha y}({}^{\alpha^y}x)$ or ${\bf l1(\alpha,y,x)}$. We group this with ${\bf lr3(\alpha,y,x)}$ as ${\bf l2(\alpha,y,x)}$, as in the local version of the proof of Lemma~\ref{leml1r2} applied in a part of $Z^2$. Since these remarks hold for all $x,y\in X$ and $\alpha\in Y$, we have proven {\bf mr1, ml2} as required. The other two follow similarly from a matched pair $(U,T)$.  
\end{proof}

We have  found explicit minimal conditions on the ground actions
(\ref{groundactions})
which are necessary and sufficient to have extensions with nice properties.  
Now we give a global description of the nature of
extensions of this type. Theorem
\ref{S,S matchedth}  tells us that $r_Z$ obeys the YBE when
constructed in this form of a matched pair extension $(U,U)$, and
we shall see that essentially this form is forced on us for
any extension with nice properties.

\begin{definition} A regular extension of {\bf M3}-monoids $S,T$ is an {\bf M3}-monoid $U$ such that $U=S\bowtie T$ where $(S,T)$ is a strong matched pair and the actions of $(U,U)$ extend the actions of $(S,S),(T,T),(S,T),(T,S)$. We denote the last of these by $\la,\ra$. 
\end{definition}

\begin{remark}  If the actions in the initial matched
pairs extend, then 
\begin{equation}\label{STaction} {}^v(u.a)={}^vu(v^u\la a),\quad  {}^b(u.a)={}^bu.{}^{b^u}a \end{equation}
are the only possible definitions for the actions of $S,T$. Hence the extended actions necessarily take the form
 \begin{equation}\label{extaction} {}^{(v.b)}(u.a)={}^v({}^bu).(v^{{}^bu}\la {}^{b^u}a),\quad (v.b)^{(u.a)}=(v^{{}^bu}\ra {}^{b^u}a) . (b^u)^a.\end{equation}
Equivalently, 
the associated map $r_U$  in terms of the associated maps for each matched pair  takes the form 
\begin{equation}\label{extrU}r_U=r^{23}_{T,S}{}^{-1}\circ r^{34}_{T}\circ r^{12}_{S}\circ r^{23}_{T,S}:S\times T\times S\times T\to S\times T\times S\times T,\end{equation}
where the numerical indices denote the positions in which the map is applied.
\end{remark}

\begin{theorem}\label{M3exttheorem} Let $U=S\bowtie T$ for $(S,T)$ a strong matched pair of {\bf M3}-monoids. The following are equivalent:

(1) $U$ is a regular extension of {\bf M3}-monoids.

(2) $(U,T), (S, U)$ are matched pairs extending the given actions.

(3) {\bf ml2,mr2} defined analogously to Notation~\ref{mlnotation} but for $S, T$, hold.
\end{theorem}
\begin{proof} The conditions in part (3) are the monoid versions of {\bf ml1, mr2} in the notation used previously, except that we apply them now to the sets of the monoids $S,T$ in the role of $X,Y$. Let us also decompose {\bf ml2=l2(T,S,S)} into pieces  {\bf ml1a=l1(T,S,S)} and {\bf lr3a}  analogously to Section~2, and similarly for {\bf mr2= r2(T,T,S)}. Explicitly the conditions of part (3) therefore read:
\[ {\bf ml1a}:\quad {}^{{}^au}({}^{a^u}v)={}^a({}^uv),\quad 
 {\rm\bf lr3a:}\quad ({}^au)^{({}^{a^u}v)}={}^{(a^{{}^uv})}(u^v)\]
 \[ {\bf mr1a}:\quad (a^b)^u=(a^{{}^bu})^{b^u},\quad  {\rm\bf lr3b:}\quad
  ({}^ab)^{({}^{a^b}u)}={}^{(a^{{}^bu})}(b^u)\]
  for all $u,v\in S$ and $a,b\in T$. The {\bf ml1a, mr1a} are the same conditions as previously but applied here to $S, T$ in the role of $X,Y$. We first show that these four are equivalent to part (2).
 
Under the hypothesis of the theorem we have matched pairs $(S,S), (T,S), (S,T), (T,S)$ and we
define left actions of $T,S$ on $U$ by the necessary formula (\ref{STaction}).  That these are
separately well-defined actions of $S,T$ follows from
\[ {}^{b}({}^c(ua)={}^b({}^cu.{}^{c^u}a)={}^b({}^cu).{}^{b^{{}^cu}c^u}a={}^{bc}u.{}^{(bc)^u}a={}^{bc}(u.a)\]
and a similar computation for $\la$, using only that the initial
matched pair data. Similarly for the right actions of $S,T$ on $(u.a)$.

In particular for the $(U,T)$ matched pair we look at the above left action of $T$ on $U$ and a right action of $U$ on $T$ given necessarily by  $b^{u.a}=(b^u)^a$. Since the cross-relations of $U=S\bowtie T$ are $au={}^au.a^u$, the latter is an action exactly when {\bf mr1a} holds. In this case
\[ (bc)^{u.a}=((bc)^u)^a=(b^{{}^cu} c^u)^a=(b^{{}^cu})^{{}^{(c^u)}a}.(c^u)^a=b^{({}^cu.{}^{c^ua})}(c^u)^a=b^{{}^c(u.a)} c^{u.a}\]
again from just the initial matched pair data and the definitions. So the matched pair condition on this side holds automatically. For the condition on the other side 
\begin{eqnarray*} {}^c((u.a)(v.b))&=& {}^c(u.{}^av. a^v. b)={}^c(u.{}^av){}^{c^{(u.{}^av)}}(a^v.b)={}^cu.{}^{(c^u.a)}v. {}^{(c^{(u.{}^av)})}(a^v).{}^{(c^{(u.{}^av)})^{a^v}}b\\
{}^c(u.a).{}^{c^{u.a}}(v.b)&=&({}^cu.{}^{c^u}a).({}^{c^{u.a}}v.{}^{(c^{u.a})v}b)={}^cu.{}^{(({}^{c^u}a)(c^u)^a)}v.({}^{c^u}a)^{({}^{c^{u.a}}v)}.{}^{(c^{u.{}^av})^{a^v}}b
\end{eqnarray*}
using only the definitions. The second factors agree by the {\bf M3} property in $T$. The outer factors agree, and the third factors agree if {\bf lr3b} holds. This is also necessary for agreement if one sets $u=b=1$ and uses unique factorisation in $U=S\bowtie T$ (i.e. an equality means equality of each factor in normal form). Similarly, we have the actions for a matched pair $(S,U)$ precisely when {\bf ml1a} holds and they form a matched pair precisely when {\bf lr3a} holds. Hence (2) $\Leftrightarrow$ (3).

Clearly (1) implies (2) by restricting the actions in the $(U,U)$ matched pair. We now show the converse, i.e. that the $(S,U)$ and $(U,T)$ matched pairs automatically extend to a $(U,U)$ one with actions of the form (\ref{extaction}). Note that  if our actions do form a matched pair, then
\begin{eqnarray*} (vb)(ua)&=&v.{}^bu.b^u.a={}^v({}^bu)(v^{{}^bu})({}^{b^u}a)(b^u)^a={}^v({}^bu).(v^{{}^bu}\la {}^{b^u}a).(v^{{}^bu}\ra {}^{b^u}a) . (b^v)^a\\
&=&{}^{(vb)}(ua) (vb)^{(ua)}.\end{eqnarray*} i.e. \textbf{M3}
necessarily holds for $(U,U)$.

\begin{figure}
\[ \includegraphics{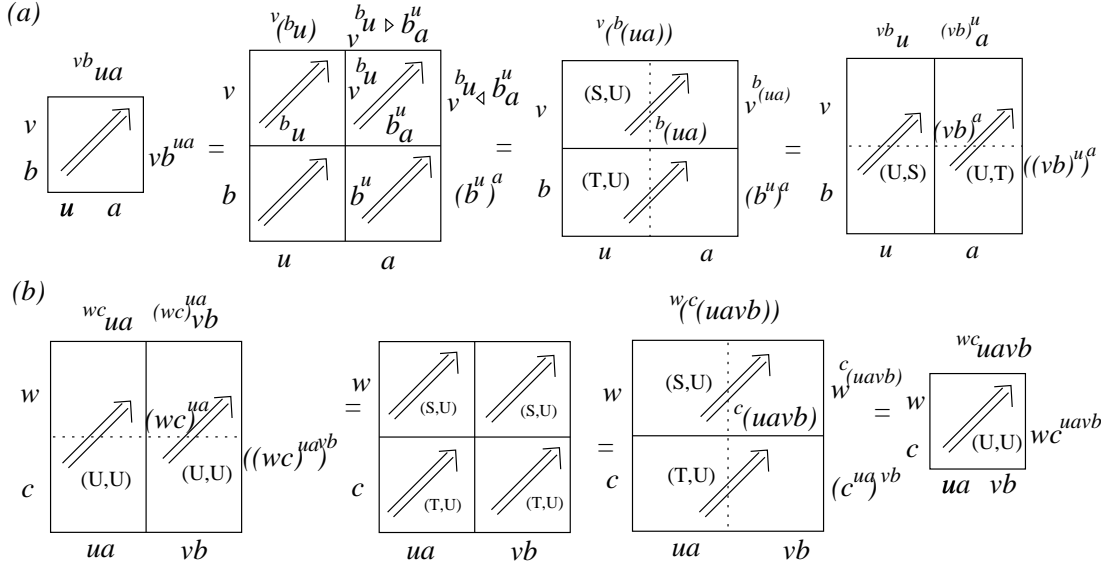}\]
\caption{Proof that extension forms a matched pair. (a) definition
of extended actions, (b) proof of vertical
subdivision property.}
\end{figure}

In order to proceed, we first prove an equivalent version of part (2); we can equally well construct matched pairs $(U,S),(T,U)$ in the same way. To do this, we use Lemma~\ref{VIPlemma} applied now to the sets $S,T$. The conditions {\bf mr1=r1(T,S,S)} which is part of the $(S,T)$ matched pair hypothesis  and {\bf l2(T,S,S)} together are equivalent to {\bf r1(S,S,T), l2(S,S,T)}. Similarly {\bf ml1=l1(T,T,S)} and  {\bf r2(T,T,S)} together are equivalent by the lemma to {\bf r2(S,T,T), l1(S,T,T)}. Note that we do not actually need $r_S,r_T$ to be invertible in the proof of the lemma and we do not have this here, otherwise the proof is the same (alternatively one can understand {\bf L2} as a covariance property of $r_{T,S}$ and this holds if an only if its inverse is covariant; one can eventually deduce the required result from this). We write  {\bf mr1a'=r1(S,S,T)}  and {\bf ml1', lr3a'}  for the two parts of {\bf l2(S,S,T)}. Explicitly:
\[ {\bf mr1a'}:\quad u^v\ra a=(u\ra{(v\la a)})\ra (v\la a),\quad {\bf lr3a'}:\quad {}^uv\ra {({u^v}\la a)}={(u\ra({v\la a}))}\la (v\ra a).\]
Here {\bf ml1'} is  part of the $(T,S)$ matched pair data that $S$ acts on $T$ from the left. Meanwhile {\bf mr1a'}  says that $U=T\bowtie S$ acts on $S$ from the right as an extension of the given actions. That $S$ acts on $U$ from the matched pair hypothesis and one of the matching conditons requires precisely the  {\bf lr3a'} condition. The details are strictly analogous to our proof above so we omit them. In the same way, that the relevant actions extend to matched pairs $(T,U)$ uses {\bf ml1a'=l1(S,T,T)} and the two parts {\bf mr1'} and {\bf lr3b'} of {\bf l2(S,T,T)}. The middle one is part of the $(T,S)$ matched pair data and the the other two read explicitly
\[ {\bf ml1a'}:\quad u\la {}^ab={}^{u\la a}((u\ra a)\la b),\quad {\bf lr3b'}:\quad 
  (u\la a)\ra ((u\ra a)\la b)={(u\ra {{}^ab})}\la a^b.\]

In principle one can now proceed to use these variant conditions as well as the original conditions (3) to verify the matched pair structure for $(U,U)$. We have given them explicitly in case the reader wishes to verify any parts of this (for example that  
  that  $v\in S,b\in T$ acting on $U$ as above indeed form a representation of
$S\bowtie T$ is a reasonable computation).  Here we provide a diagrammatic proof  using the  subdivision form of the matched pair axioms explained in Figure~1(a) and the results already proven.  In our case the proposed 
action of $S\bowtie T$ on itself is given in terms of its
composite square in Figure~3(a). By definition it
consists of composing the four actions in our initial data
according to the same 'transport' rules, namely that the action
takes place as the element is taken through the box (to the top or
to the right). That this composite could also be viewed as a
vertical composite of two horizontal rectangular boxes or a horizontal composite of two vertical boxes, is nothing other than a restatement of the definitions. In the first case the $(U,U)$ actions are represented as the $(S,U)$ actions placed above the $(T,U)$ actions. In the second case it is represented as the $(U,S)$ actions placed to the left of the $(U,T)$ actions. 

We can now prove the subdivision property as follows. The horizontal one is shown in Figure~3(b). On the left is the horizontal composition of two instances of the $(U,U)$ actions. We wish to show that this can be merged into a single $(U,U)$ box as on the right. To do this we break the boxes on the dotted line to give four boxes. The lower two are two horizontally composed instances of the $(T,U)$ actions and since this is a matched pair by the variant of part (2) it composes to a single horizontal box. Similarly, the two upper boxes are horizontally composed $(S,U)$ actions and hence by part (2) they compose to a single horizontal box. This is shown in the middle of Figure~3(b). Finally, the placement of one horizontal box vertically over the other is nothing other than one definition of the $(U,U)$ actions as explained above. 

The vertical subdivision property is entirely analogous. We break two vertically composed $(U,U)$ boxes into four boxes. The two on the left are vertically composed instances of $(U,S)$ boxs, the two on the right are vertically composed instances of $(U,T)$ boxes. Since we have proven above that these are matched pairs as part (2) and its variant, we can merge these boxes and recognise a single $(U,U)$ box as a single $(U,S)$ box followed horizontally by a single $(U,T)$ box, which is our alternate definition of the $(U,U)$ box as explained above. \end{proof}

\begin{corollary} \label{braextcor} Let $U=S\bowtie T$ be a regular extension of {\bf M3}-monoids.

(1) $U$ is strong {\em iff} $S,T$ are.

(2) $U$ is  braided {\em iff} $S,T$ are braided.
\end{corollary}
\begin{proof} Part (1) is clear as $r_U$ is invertible precisely when both $r_S,r_T$ are. For part (2) 
we recall that {\bf mr1=r1(T,S,S), l2(T,S,S)} have the interpretation explained in  the proof of Lemma~\ref{VIPlemma} of the YBE restricted to $T\times S\times S$, the difference is that we now work on the sets $S,T$ with the restrictions of $r_U$ but the diagram involved is analogous. Similarly for {\bf ml1, r2(T,T,S)}. By Lemma~\ref{VIPlemma} we  deduce that the YBE holds for  all mixed triples (those involving elements of both $S,T$) using the restrictions $r_S,r_T,r_{S,T},r_{T,S}$. We have seen in the theorem that that these are precisely the conditions for a regular extension of {\bf M3}-monoids. Hence, given this, the YBE holds for all combinations precisely when it holds for $r_S$ and $r_T$ separately. \end{proof}

\begin{corollary}\label{doublebramon} If $(S,r_S)$ is a braided monoid 
then it has a canonical `quantum double' braided monoid $(S\bowtie S,r_{S\bowtie
S})$. 
\end{corollary}
\begin{proof} The conditions in  part (3) of the theorem hold because all the actions involved are based on the actions of $(S,S)$ and $(S,r_S)$ is assumed to be a braided monoid.  In more explicit terms    {\bf ml1a, mr1a}  for $S, T$ hold automatically as the actions
are given by the actions of $S$ on itself, which is a strong matched pair $(S,S)$ by assumption, while
conditions {\bf lr3a, lr3b} hold as part of the assumption that $r_S$ obeys YBE. We then use the theorem and the preceding corollary. \end{proof}

We remark that this is exactly
analogous to the construction of a quasitriangular structure or
`doubled braiding' on $A(R)\bowtie A(R)$ as a version of
Drinfeld's quantum double, see \cite{Ma:book}.   In the case of
$S$ given by a braided set $(X,r)$, it corresponds to the braided
set $(X\sqcup X,r_D)$ in Corollary~\ref{dbraidedset}. 

\begin{proposition} Suppose that $S,T$ are cancellative {\bf M3}-monoids forming a strong matched pair. Then $U=S\bowtie T$ is a regular extension {\em iff} {\bf ml1a, mr1a} hold for $S,T$. In  this case $U$ is cancellative.
\end{proposition}
\begin{proof} Since $(S,T)$ are a matched pair we can argue that \begin{eqnarray*}{}^a(uv)&=&{}^a({}^uv.u^v)={}^a({}^uv).{}^{(a^{{}^uv})}(u^v)\\
&=&{}^au.{}^{a^u}v={}^{{}^au}({}^{a^u}v).
({}^au)^{({}^{a^u}v)}\end{eqnarray*} using the {\bf M3} property for the first and last equalities. So if we assume the
\textbf{ml1a} condition stated and left cancellation in $S$ we
will have \textbf{lr3a} above.  Similarly for \textbf{lr3b} using
\textbf{mr1a} and right cancellation in $T$. If we assume cancellation on both sides then
suppose $(u.a)(v.b)=(u'.a').(v.b)$ in $U$. Using the cross-relations in $U=S\bowtie T$ and the unique
factorisation there, this means
\[ u. {}^{a}v=u'.{}^av,\quad a^v.b=a'{}^v.b\]
From the second and right cancellation in $T$ we have $a'{}^v=a^v$. Using (\ref{lara}) we deduce
\[ u\la a=(u {}^av)\la a^v=(u'{}^{a'}v)\la a{}^v=u'\la({}^{a'}v\la a'{}^v)=u'\la a'\]
\[(u\ra a)v=u\ra({}^av\la a^v).({}^av\ra a^v)= (u{}^av)\ra a^v= (u'{}^{a'}v)\ra a^v= (u'{}^{a'}v)\ra a'{}^v=(u'\ra a').v\]
hence if we have right cancellation in $S$ we deduce that $u'\ra a'=u\ra a$. It follows from (\ref{lara}) that $u'=u, a'=a$. From the other side, $(u.a)(v.b)=(u.a).(v'.b')$ means 
\[ u.{}^av=u. {}^av',\quad a^v.b=a^{v'}.b'\]
so by left cancellation in $S$ we have ${}^av={}^{a}v'$. This and the remainder of the proof is analogous to the proof on the other side already given.  \end{proof}

\section{Application to symmetric sets and strong twisted unions}

In this section we apply the matched pair approach developed in the paper to construct extensions in some important special cases such as nondegenerate involutive solutions (or as they are often called \emph{symmetric sets}) $(X,r)$,    strong twisted unions of solutions, see Definition~\ref{STUdef},  finite square-free solutions, etc. Note that every finite non-degenerate  involutive
square-free solution $(Z,r)$ is an extension of some nonempty
disjoint square-free involutive solutions $(X,r_X),$ and
$(Y,r_Y),$ \cite{Rump}. Furthermore,
in this case the monoid $S(Z,r)$ and the YB-algebra $\Acal(k,Z,r)$
(over arbitrary  field $k$)  have remarkable algebraic and homological properties, see \cite{TM, T2, T06, T07, JO}.
In \cite{TGISM} our matched pair approach is applied to various special cases of extensions of solutions, in particular to strong twisted union of solutions, finite symmetric sets, finite square-free involutive solutions, etc.

\begin{definition}
\label{STUdef}
We call a regular extension $(Z,r)$ a \emph{strong twisted union} of 
the quadratic sets $(X,r_X)$ and $(Y, r_Y)$ if
\begin{enumerate}
\item 
\label{a)}
The assignment $\alpha
\longrightarrow {}^{\alpha}\bullet $ 
extends to a left action of the
associated  group $G(Y, r_Y)$ (and the associated monoid $S(Y,r_Y)$) on $X$,  
and the assignment $x \longrightarrow {\bullet}^x $ 
extends to a
right action of the associated  
group of $G(X,r_X)$  (and the associated monoid $S(X,r_X)$)
  on $Y$;
\item
\label{b} The regular pair of ground actions satisfy 
 \[\begin{array}{lclc}
 {\rm\bf stu :}\quad &{}^{{\alpha}^y}x = {}^{\alpha}x;
\quad  &{\alpha}^{{}^{\beta}x}={\alpha}^x,\quad  \text{for all}\quad x, y \in X, \alpha,\beta \in Y 
\end{array}\]
\end{enumerate}
 \end{definition}
\begin{remark}
In \cite{ESS}, Definition 3.3. the notion of \emph{a generalized twisted union} $(Z,r)$ of the solutions 
$(X,r_X)$ and $(Y, r_Y),$ is introduced in the class of symmetric sets. More precisely,  a symmetric set $(Z,r)$ is a generalized twisted union of the disjoint symmetric sets  $(X,r_X)$, and $(Y, r_Y)$
if it is an extension, and for every $x\in X, \alpha\in Y$ the ground action 
${}^{{\alpha}^x}{\bullet} : Y\times X \longrightarrow X $ does not depend on $x$, and the ground action 
${\bullet}^{{}^{\alpha}{x}} : Y\times X \longrightarrow Y$ does not
depend on  $\alpha.$
 
Note that in contrast with a generalized twisted union,  a strong twisted union $(Z,r)$ of 
$(X,r_X)$ and $(Y, r_Y)$ does not necessarily obey YBE and is not limited to symmetric sets. It is easy to see that a strong twisted union $(Z,r)$ of symmetric sets which obeys YBE, is a generalized twisted union.  Furthermore, it will be shown in  \cite{TGISM} that 
a generalized twisted union $(Z,r)$ of two involutive square-free solutions is a strong twisted union. 
\end{remark} 

\begin{remark} Let $(X,r)$ be a quadratic set.  A permutation $\tau \in \Sym(X)$ is called \emph{an automorphism
of $(X, r)$} (or shortly an \emph{r-automorphism}) if
$(\tau \times \tau) \circ  r = r \circ (\tau \times \tau)$. The
group of $r$-automorphisms of $(X, r)$ will be denoted by $\Aut(X,
r)$. In \cite{TGISM} it will be shown that  a strong twisted union $(Z,r)$ of 
solutions $(X,r_X)$ and $(Y, r_Y)$ obeys the YBE {\em if and only if} the assignment 
$\alpha
\longrightarrow {}^{\alpha}\bullet $
extends to a a group homomorphism 
\[G(Y,r_Y) \longrightarrow Aut(X,r) \]
and the assignment $x \longrightarrow {\bullet}^x$
extends to a a group homomorphism 
\[ G(X,r_X) \longrightarrow Aut(Y,r).\]
\end{remark}

Here we limit ourselves to some first results connecting with those of Section~2.

\begin{proposition}
\label{stuprop}
Let  $(Z,r)$ be a strong twisted union of the  2-cancellative disjoint solutions, 
$(X,r_X)$ and $(Y, r_Y)$. Suppose {\bf lri} holds for $(Z,r)$ and  $U=S(Z,r)$ has cancellation on monomials of length 3. Then $(Z,r)$ obeys YBE \emph{iff} 
\[\begin{array}{lclc}
 {\rm\bf csla:}\quad 
&{}^{{}^y{\alpha}}{({}^yx)}={}^{{}^{\alpha}y}{({}^{\alpha}x)}
\quad  &{\rm\bf csra:}\quad 
({\alpha}^{\beta})^{x^{\beta}}=({\alpha}^x)^{{\beta}^x}, 
\end{array}\]
for all $x, y \in X, \alpha,\beta \in Y.$
\end{proposition}
\begin{proof}
Note first that the hypothesis of Theorem~\ref{theoremBcancellative}
is satisfied, so 
$(Z,r)$
obeys YBE {\em iff}   {\rm\bf ml1, mr1},\mla,\mra, where we use our original notions for $Z=X\sqcup Y$. 
As a strong twisted union $(Z,r)$ satisfies {\rm\bf ml1, mr1}.
Now we interprete condition \mla:
\begin{equation}
\label{mlaeq1}
{}^{\alpha}{({}^yx)}={}^{{}^{\alpha}y}{({}^{{\alpha}^y}x)}=^{ {\rm\bf stu}} {}^{{}^{\alpha}y}{({}^{\alpha}x)},
\end{equation}
Apply  {\rm\bf stu } again to yield
\[
{}^{\alpha}{({}^yx)}={}^{{}^y{\alpha}}{({}^yx)},
\]
which together with (\ref{mlaeq1}) gives $\mla \Longleftrightarrow {\rm\bf csla}$. 
One similarly finds $\mra \Longleftrightarrow {\rm\bf csra}$ under our hypotheses.
\end{proof}
\begin{remark}
\label{extcyclesetscor}
Under the hypothesis of Proposition \ref{stuprop}
suppose  $(X,r_X)$ and  $(Y, r_Y)$ are symmetric sets. Then \begin{enumerate} 
\item $(Z,r)$  is nondegenerate, involutive  and the
cyclic conditions hold.
\item $(Z,r)$ is a solution of YBE \emph{iff} $(Z,r)$ 
is a cycle set, i.e. conditions {\bf csl, csr} hold.
In this case $(Z,r)$ is a symmetric set. 
\end{enumerate}
This follows immediately from Proposition~\ref{lri&INVOLUTIVE} on noting that the involutiveness of  $(X,r_X)$ and  $(Y, r_Y)$ imply $(Z,r)$ involutive while {\bf lri} holds  by assumption. \end{remark}

\begin{proposition}
\label{extensiontrivialsolutions}
Let $(Z,r)$ be a regular extension  of two trivial solutions  $(X,r_X),$ $(Y, r_Y)$. Suppose the monoid $U=S(Z,r)$ is with cancellation on monomials of length 3. Then  $(Z,r)$  obeys YBE {\em iff} it is a strong twisted union
of $(X,r_X), (Y, r_Y)$. 
\end{proposition}
\begin{proof}
More generally, under the hypothesis of 
Theorem
\ref{theoremBcancellative}, suppose $(Z,r)$ is a regular extension of the involutive solutions  $(X,r_X), (Y, r_Y)$. Suppose furthermore that $(X,r_X)$ is a trivial solution, i.e. $r_X(xy)=yx$ for all $x,y
\in X.$ This gives ${}^xy=y,$ for all $x,y \in X,$ 
thus an easy computation shows
 \begin{equation}
\label{extoftrivialsol1}
{\rm\bf \mla } \Longleftrightarrow {}^{{\alpha}^x}y =
{}^{\alpha}y,\quad \quad {\rm\bf mr1} \Longleftrightarrow (\alpha^x)^y
= (\alpha^y)^x 
 \end{equation}
 for all $x,y \in X, \alpha \in Y$.
Analogously,
in the case when $(Y, r_Y)$ is a
trivial solution
one has 
\begin{equation}
\label{extoftrivialsol2}
{\rm\bf \mra } \Longleftrightarrow {\alpha}^{{}^{\beta}x} =
{\alpha}^x 
\quad \quad {\rm\bf ml1} \Longleftrightarrow {}^{\alpha}{({}^{\beta}x)}= {}^{\beta}{({}^{\alpha}x)}.
 \end{equation}
Assume now both $(X,r_X), (Y,r_Y)$ are trivial solutions.
 Clearly the  hypothesis of Theorem
\ref{theoremBcancellative} is satisfied so 
(\ref{extoftrivialsol1}),  (\ref{extoftrivialsol2}) are in force.
This yields
\[
{\rm\bf \mla, \mra } \Longleftrightarrow {\bf stu},
\] and  
\[{\rm\bf ml1, mr1 } \Longleftrightarrow \text{part\ (1)\ of\ Definition}\ \ref{STUdef},
\] 
which completes the proof.  
\end{proof}

Given a quadratic set $(X,r_X),$ the set of 
$r_X$-\emph{fixed pairs} is $\{x, y \in X \mid r_X(xy)=xy\}$.  We have seen in section 2 that in the case of involutive nondegenerate solutions there is a close relation between the set of $r_X$-fixed pairs and  2-cancellativeness, see Lemma~\ref{Zrlemma3} and Corollary~\ref{involcancelcorol}. 
We will show now that for involutive nondegenerate solutions 
$(X,r_X),$ and
$(Y,r_Y),$ the 
conditions {\bf \mll , \mrr } can be interpreted in terms of the $r_X$ and $r_Y$-fixed pairs. This is especially useful in the case 
of square free solutions, where the set of $r_X$-fixed pairs is exactly the diagonal, \emph{diag}$X\times X.$

\begin{lemma} 
\label{involutiveextensionscor}
In the notation of Theorem~\ref{theoremB},
suppose $(X, r_X)$ and $(Y,r_Y)$ are involutive solutions, and   
$S,T$ is a matched pair of monoids. Let $(Z,r)$ be the associated regular extension. Then 
\[
(i)\quad {\bf \mll } \Longleftrightarrow 
[r_X(xy)=xy\Longrightarrow r_X({}^{\alpha}{xy})={}^{\alpha}{xy}, \;\text{for all}\; x, y \in X, \alpha\in Y]. 
\]
\[
(ii)\quad {\bf \mrr } \Longleftrightarrow [
r_Y(\alpha\beta)=\alpha\beta\Longrightarrow r_Y((\alpha\beta)^x)=(\alpha\beta)^x, \;\text{for all}\; x \in X, \alpha, \beta\in Y]. 
\]
\end{lemma} 
\begin{proof}
It is clear that when 
$r_X$ is involutive condition \mlw\ implies that for every
$x,y \in X, \alpha\in Y$ one of the following is an equality of
words, either 
\begin{equation}
\label{eqfixed1}
{}^{\alpha}{(r_X(xy))}=r_X({}^{\alpha}{(xy)})\quad \text{holds in}
\;\;  X^2 ; 
\end{equation}
or
\begin{equation}
\label{eqfixed2} {}^{\alpha}{(r_X(xy))}={}^{\alpha}{(xy)}\quad \text{holds in}
\;\;  X^2.
\end{equation}
Assume now (\ref{eqfixed2}) holds. Then the equalities 
\[ 
{}^{\alpha}{(r_X(xy))}={}^{\alpha}{({}^xy.x^y)}={}^{\alpha}{({}^xy)}.
{}^{{\alpha}^{({}^xy)}}{(x^y)}
\]
and (\ref{eqfixed2}) imply 
\[
{}^{\alpha}{(xy)}= {}^{\alpha}{x}.{}^{{\alpha}^x}{y)}={}^{\alpha}
{({}^xy)}.
{}^{{\alpha}^{({}^xy)}}{(x^y)}\quad\text{in}
\;\;  X^2 . 
\]
Therefore
\[
{}^{\alpha}{x}={}^{\alpha}
{({}^xy)},
\]
which by the nondegeneracy implies $x ={}^xy.$
This by Lemma~\ref{Zrlemma3}, is equivalent  to $r(x,y) = (x,y),$ i.e. 
$(x,y)$ is an $r_X$ fixed pair. Clearly for a $r_X$-fixed pair
$(x,y)$ (\ref{eqfixed1}) holds \emph{iff} $r_X({}^{\alpha}{(xy)})={}^{\alpha}{(xy)}$. This proves (i). Analogous argument proves (ii). 

\end{proof}
\begin{proposition} 
\label{theoremD} Let $(X,r_X), (Y,r_Y)$ be
non-degenerate involutive disjoint solutions  with YB-monoids
$S=S(X,r_X)$, $T=S(Y,r_Y).$  Suppose $(S,T)$ is a strong matched
pair, and  let $(Z,r),$ be the associated regular extension. Then
\begin{enumerate}
 \item
\label{i1theoremD}
 $r$ obeys YBE \emph{iff} 
\[
\begin{array}{lclcr}
&{}^xy=x\Longrightarrow
{}^{{}^{\alpha}x}{({}^{{\alpha}^x}y)}={}^{\alpha}x,\quad&
{\alpha}^{\beta}={\beta}\Longrightarrow
(\alpha{}^{{}^{\beta}x})^{{\beta}^x}={\beta}^x,& {\rm for\ all}\ x, y \in X, \alpha,\beta \in Y.
\end{array}
\]
\item \label{i3theoremD} Suppose 
 that both $(X,r_X), (Y,r_Y)$  are square
free. Then  $r$ obeys YBE {\em iff} the following ``mixed'' weak cyclic conditions
hold:
\[
 {}^{\alpha^x}x = {}^{\alpha}x, \;\;\text{and} \;\;
 {\alpha}^{{}^{\alpha}x} = {\alpha}^x, \;\; \text{for all}
\;\; x \in X, \alpha \in Y.
\]
In this case $(Z,r)$ is also square-free and satisfies {\rm\bf lri}.
\end{enumerate}
\end{proposition}
\begin{proof}
 We give sketch of the proof. To verify (1) 
one uses Lemma~\ref{Zrlemma3}  again and  interpretes the implications in Lemma~\ref{involutiveextensionscor} in terms of the left and the right actions.   

Assume now  $(X,r_X)$ and $(Y, r_Y)$ are square-free. Then,
$(x,y)$ is an $r_X$ fixed point iff $y=x,$ similarly, $(\alpha\beta)$
is an $r_Y$ fixed point iff $\alpha=\beta$, (see Lemma 
\ref{someproperties}\ref{SFnondeg}).
 
Next we replace $y=x$, $\alpha = \beta $ in (\ref{i1theoremD}) and obtain 
the implications in (\ref{i3theoremD}).
\end{proof}

 \begin{figure}
 \[ \includegraphics{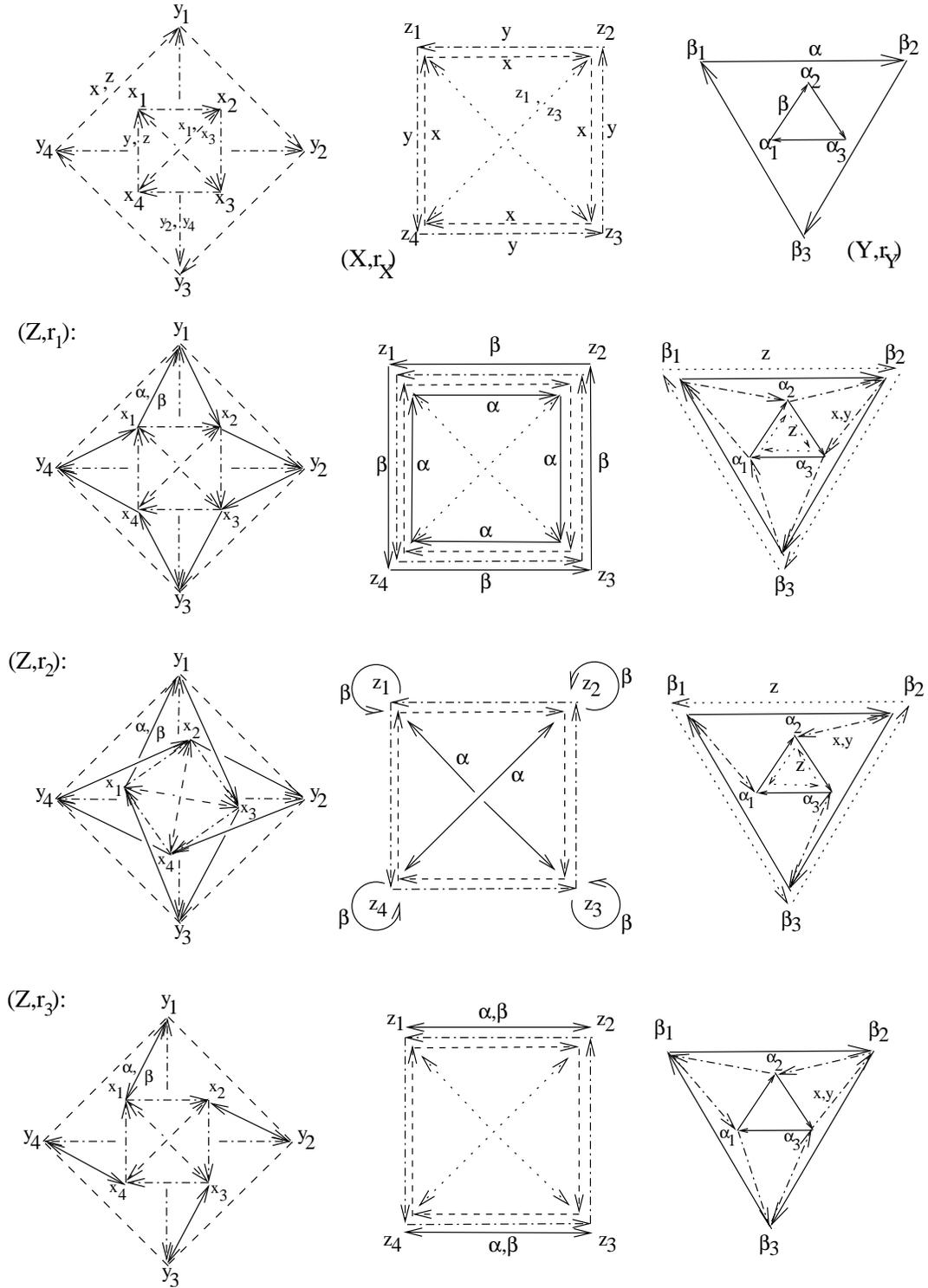}\]
 \caption{Example~\ref{ext1}  of extensions $(Z,r)$}
 \end{figure}

We now give various examples of extensions $(Z, r)$ of a fixed
pair of solutions $(X,r_X), (Y,r_Y).$ All solutions  are
involutive, non-degenerate, square-free,  with \textbf{lri}. In
this case $\Rcal_z = \Lcal_z^{-1}$ for all $z \in Z.$  We find the
extensions applying effectively Theorem~\ref{theoremB} and our
theory about the behaviour of finite square-free solutions (rather
than by computer).

\begin{example}\label{ext1}

Let $X= \{x_1,x_2,x_3,x_4; y_1,y_2,y_3,y_4; z_1,z_2,z_3,z_4 \}$
Let $\rho,\sigma, \tau$ be the following cycles of length 4 in
$Sym(X)$
\[\rho=(x_1,x_2,x_3,x_4), \quad \sigma =(y_1, y_2,
y_3,y_4), \quad \tau= (z_1,z_2,z_3,z_4) \] Define the involutive
map $r_X: X^2 \longrightarrow X^2$ as:
\[\begin{array}{lclc}
  r_X(x_iy_j)= \sigma(y_j)\rho^{-1}(x_i),&
 r_X(x_iz_j)=\tau(z_j)\rho^{-1}(x_i),\quad
 &r_X(y_ix_j)= \rho(x_j)\sigma^{-1}(y_i),\\
 r_X(y_iz_j)= \tau^{-1}(z_j)\rho^{-1}(y_i), &
r_X(z_ix_j)= \rho(x_j)\tau^{-1}(z_i), \quad &
r_X(z_iy_j)= \rho(y_j)\tau(z_i) \\
r_X(x_1x_2)= x_4x_3,\;&r_X(x_1x_4)= x_2x_3,\quad
r_X(x_3x_2)= x_4x_1,\; & r_X(x_3x_4)= x_2x_1,\\
r_X(y_1y_2)= y_4y_3,\; &r_X(y_1y_4)= y_2y_3,\quad
r_X(y_3y_2)= y_4y_1,\; &r_X(y_3y_4)= y_2y_1,\\
r_X(z_1z_2)= z_4z_3,\; &r_X(z_1z_4)= z_2z_3,\quad
r_X(z_3z_2)= z_4z_1,\; &r_X(z_3z_4)= z_2z_1,\\
r_X(x_1x_3)= x_3x_1, \; &r_X(x_2x_4)= x_4x_2,\quad r_X(y_1y_3)=
y_3y_1, &r_X(y_2y_4)= y_4y_2,\\
 \quad \quad &r_X(z_1z_3)= z_3z_1, \quad\quad
r_X(z_2z_4)= z_4z_2.\quad&
\end{array}\]
For the left action on $X$  we have:
\[\begin{array}{lclc}
\Lcal_{x_1} = \Lcal_{x_3}=\sigma \circ \tau \circ (x_2 x_4)\quad&
\Lcal_{x_2}
= \Lcal_{x_4}=\sigma \circ \tau \circ (x_1 x_3);\\
\Lcal_{y_1} = \Lcal_{y_3}=\rho \circ\tau^{-1} \circ(y_2 y_4)\quad&
\Lcal_{y_2} = \Lcal_{y_4}=\rho \circ \tau^{-1} \circ
(y_1 y_3)\\
\Lcal_{z_1} = \Lcal_{z_3}=\rho \circ \sigma \circ (z_2 z_4)\quad&
\Lcal_{z_2} = \Lcal_{z_4}=\rho \circ \sigma \circ (z_1
z_3).\\
\end{array}\]

The second solution $(Y,r_Y)$ is simpler. Let $Y= \{\alpha_1,
\alpha_2, \alpha_3; \beta_1, \beta_2,\beta_3\}.$ Take the cycles
$f= (\alpha_1, \alpha_2, \alpha_3), g= (\beta_1, \beta_2,\beta_3)$
in $Sym(Y)$ and define $r_Y: Y^2\longrightarrow Y^2$ as
\[\begin{array}{lclc}
  r_Y(\alpha_i\beta_j)=
g(\beta_j)f^{-1}(\alpha_i),\quad &
 r_Y(\beta_j \alpha_i)= f(\alpha_i)g^{-1}(\beta_j)\\
r_Y(\alpha_i\alpha_j)= \alpha_j\alpha_i, \quad &
r_Y(\beta_i\beta_j)= \beta_j\beta_i,
\end{array}\]
for all $i,j, 1 \leq i,j\leq 3.$  The left actions on $Y$ are:
$\Lcal_{\beta_i}= f, \Lcal_{\alpha_i}= g,$ for all $1 \leq i \leq
3.$

These initial data are shown in Figure~4. The definition of the
graph is given after the examples. We omit most of the labels on
the arrows in order not to clutter the diagram.  We present
several examples of YB extension $(Z,r)$ of $X,Y$ with $Z=X\sqcup
Y.$

Note that all $(Z,r)$ are square free, as extensions of
square-free solutions, therefore we can apply the combinatorics
developed in \cite{T06}. In particular \textbf{cc } is in force
for $(Z, r).$ Clearly, ${\bf x}= \{ x_1, x_2, x_3 ,x_4 \}$ and
${\bf y}=\{ y_1, y_2, y_3 ,y_4 \}, {\bf z}= \{ z_1, z_2, z_3, z_4
\} $ are the orbits of $X$ under the action of $G(X, r_X),$ and
${\bf \alpha} = \{ \alpha _1, \alpha _2, \alpha _3 \}, {\bf
\beta}= \{\beta _1, \beta _2, \beta _3 \}$ are the orbits of $Y$,
under the action of $G(Y, r_Y),$. Each orbit in $X$ is $r_X$
invariant subsets, so we have to bear in mind also
"mini-extensions" of pairs like $(\textbf{x}, Y),$ $(\textbf{y},
Y),$ $(\textbf{z}, Y),$ $(\textbf{x}\sqcup \textbf{y}, Y)$ etc.
One can show that with this initial data, it is impossible to have
extension $(Z,r)$ in which some $x_i, y_j, z_k$ belong to the same
orbit.

We are interested in cases when the left action of $Y$ onto $X$
provides "links" between the two orbits $\textbf{x}, \textbf{y}.$
So we need extensions for the pairs $\textbf{x}\sqcup \textbf{y},
Y$ and for $\textbf{z}, Y,$ which are compatible.  We shall write
$r({\bf \alpha x})= {\bf y \beta}$,(respectively  $r({\bf \alpha
x})= {\bf y \alpha}$ to indicate that
\begin{equation}
\label{A} r(\alpha_i x_p)=  y_q \beta_j \; (\text{resp.} \;
r(\alpha_i x_p)=  y_q \alpha_j\;\text{for some}\; 1 \leq i,j \leq
3, 1\leq p,q\leq 4.
\end{equation}
More detailed study shows what kind of pairs  $(p,q), 1 \leq p, q
\leq 4$ are admissible with the structure of $(X,r_X).$ For this
particular  $Y$ every pair $(i,j), 1 \leq i,j \leq 3$ is
admissible.  Note that no restriction on compatibility between the
two pairs $(i,j)$ and $(p,q)$ are necessary. There are two types
of extensions which connect $\textbf{x}, \textbf{y},$ in one
orbit, these are \textbf{A:} $r({\bf \alpha x})=  {\bf y\beta}$; and
\textbf{B:} $r({\bf \alpha x})= {\bf y\alpha}$. The admissible
actions of $Y$ on $\textbf{z}$ depend only on the general type
\textbf{A}, or \textbf{B}, and do not depend on the particular
pairs of indices $p, q, \alpha, \beta,$ occurring in (\ref{A}).
Furthermore, for simplicity we consider the case $\Lcal_{\beta
\mid {\bf x} \sqcup {\bf y}} = \Lcal_{\alpha \mid {\bf x} \sqcup
{\bf y}},$ for all $\alpha$ and $\beta.$ We shall discuss only
case \textbf{A}. Clearly, under these assumptions $Y$ does not act
as automorphisms on $X$ iff $\Lcal_{\alpha\mid \textbf{z}}\neq
\Lcal_{\beta\mid \textbf{z}}.$ We start with a list of the left
actions of $Y$ onto ${\bf x} \sqcup {\bf y},$ satisfying
\textbf{A} and producing solutions $r.$ Assume \[\Lcal_{\beta \mid
{\bf x} \sqcup {\bf y}} = \Lcal_{\alpha \mid {\bf x} \sqcup {\bf
y}}.\] Three subcases are possible.

\textbf{A1}. The set ${\bf x} \sqcup {\bf y}$ becomes a cycle of
length $8$. Denote  $\theta=(x_1 y_1 x_2 y_2 x_3 y_3 x_4 y_4) \in
Sym({\bf x}\sqcup {\bf y})$. The admissible actions in this
subcase are:
\begin{equation}
\label{A1} {\bf (i)}\quad \Lcal_{\alpha \mid {\bf x} \sqcup {\bf
y}}= \theta; \quad {\bf (ii)} \quad\Lcal_{\alpha \mid {\bf x}
\sqcup {\bf y}} ={\theta} ^3; \quad {\bf (iii)} \quad
\Lcal_{\alpha \mid {\bf x} \sqcup {\bf y}} = {\theta} ^5; \quad
{\bf (iv)}\quad \Lcal_{\alpha \mid {\bf x }\sqcup {\bf
y}}={\theta}^7 .
\end{equation}
\textbf{A2.} $({\bf x} \sqcup {\bf y})$ splits into two disjoint
cycles of length $4$. Only two cases are admissible.
\begin{equation}
\label{A2}   {\bf (i)} \Lcal_{\alpha \mid {\bf x }\sqcup {\bf y}}=
(x_1 y_1 x_3 y_3)(x_2 y_2 x_4 y_4) =\vartheta; \quad {\bf (ii)}
\Lcal_{\alpha \mid{\bf x\sqcup y}}= ( y_3 x_3 y_1 x_1)( y_4 x_4
y_2 x_2 ) =\vartheta^{-1}.
\end{equation}
\textbf{A3.} $X_0$ splits into four disjoint cycles of length $2$.
There are three admissible actions.
\begin{equation}
\label{A3}   \Lcal_{\alpha \mid{\bf x\sqcup y}}
 = (x_1{\sigma}^i(y_1)) (x_2 {\sigma}^i(y_2))
(x_3{\sigma}^i(y_3)) (x_4{\sigma}^i(y_4)) , \;1 \leq i \leq 3.
\end{equation}
To determine the left actions of $Y$ on $X$ completely,  we need
to know admissible  actions $\Lcal_{\alpha \mid \textbf{z}}$ and
$\Lcal_{\beta \mid\textbf{z}}.$ One can verify that $\Lcal_{\alpha
\mid \textbf{z}}$  determines uniquely $\Lcal_{\beta
\mid\textbf{z}}.$ In four cases $\Lcal_{\alpha} \neq
\Lcal_{\beta},$  see ${\rm\bf (a)},{\rm\bf (b)}, {\rm\bf (g)},
{\rm\bf (h)}$, each of which produces solutions $(Z,r)$ with
$G(Z,r)$ acting on $Z$ not as automorphisms.
 We give now the list of admissible actions of $Y$ on $z,$ which agree
 with the initial data, and the assumption \textbf{A}. (Note that
 in case B, the list of admissible actions of $Y$ on $z,$ is
 different.)
\begin{equation} \label{actiononz}
\begin{array}{lclc}
 {\rm\bf (a)}\quad & \Lcal_{\alpha \mid \textbf{z}}=\tau , \Lcal_{\beta
\mid\textbf{z}}=\tau ^{-1}; \quad &{\rm\bf (b)}\quad \Lcal_{\alpha
\mid\textbf{z}}=\tau ^{-1},
\Lcal_{\beta \mid\textbf{z}}=\tau; \\
 {\rm\bf (c)}\quad& \Lcal_{\alpha \mid\textbf{z}}=\Lcal_{\beta \mid\textbf{z}}=(z_1z_2)(z_3z_4)
; \quad &{\rm\bf (d)}\quad \Lcal_{\alpha
\mid\textbf{z}}=\Lcal_{\beta \mid\textbf{z}}=(z_1z_4)(z_2z_3);\\
{\rm\bf (e)}\quad&\Lcal_{\alpha \mid\textbf{z}}= \Lcal_{\beta
\mid\textbf{z}}=(z_1z_3); \quad &{\rm\bf (f)}\quad \Lcal_{\alpha
\mid\textbf{z}}=
\Lcal_{\beta \mid\textbf{z}}=(z_2z_4);\\
{\rm\bf (g)}\quad& \Lcal_{\alpha \mid\textbf{z}}=(z_1z_3)(z_2z_4);
\Lcal_{\beta \mid\textbf{z}} = \id_{\bf z};   \quad &{\rm\bf
(h)}\quad \Lcal_{\alpha \mid\textbf{z}}=id_{\textbf{z}};
\Lcal_{\beta \mid\textbf{z}}=(z_1z_3)(z_2z_4).
\end{array}\end{equation}
We now have a list of admissible actions of $Y$ on $X$ that are
compatible with $r$ obeying the YBE; it remains to present
similarly admissible actions of $X$ on $Y.$ We give next two types
of actions of ${\bf x}\sqcup{\bf y}$ upon $Y$ which are admissible
with the  list of actions already specified above. (This is not a
complete list of admissible choices). Each of the actions below
glues $Y$ in one orbit:
\begin{equation}
\label{XactsonY1} \Lcal_{x_j \mid Y} = \Lcal_{y_j \mid Y}={\pi}^q,
\quad\text{where},\quad \pi = (\alpha_1\beta_1
\alpha_2\beta_2\alpha_3\beta_3), \quad q =  1, 3, 5,
\end{equation}
and
\begin{equation}
\label{XactsonY2}
  \Lcal_{x_j \mid Y} = \Lcal_{y_j \mid Y}=
 (\alpha_1g^k(\beta_1))\circ(\alpha_2g^k(\beta_2))\circ(\alpha_3g^k(\beta_3)),\quad 0\leq k\leq 2.
\end{equation}
 As a final step we determine the admissible actions of ${\bf
z}$ on $Y$. The choice is  limited:
\begin{equation}
\label{XactsonY3} \Lcal_{z_j \mid Y }= (f\circ g)^{k},\quad 0\leq
k\leq 2.
\end{equation}
Note that to make a list of what we call \emph{admissible actions}
is enough to verify conditions {\bf ml1, \mla }. We leave as an
exercise to the reader to check that  any 4tiple of actions
\[\Lcal_{\alpha \mid {\bf x}\sqcup {\bf y}}\circ \Lcal_{\alpha \mid {\bf z}}, \quad
\Lcal_ {\beta \mid {\bf x}\sqcup {\bf y}}\circ \Lcal_{\beta \mid
{\bf z}}, \quad  \Lcal_{x \mid Y} = \Lcal_{y \mid Y}, \quad
\Lcal_{z \mid Y}\] chosen from the lists (\ref{A1})--(\ref{A3}),
(\ref{actiononz}),  (\ref{XactsonY1})--(\ref{XactsonY2}), and
(\ref{XactsonY3}), respectively,  satisfies {\bf ml1, \mla }, and
therefore defines a solution $(Z, r).$  Note that different
triples can produce isomorphic extensions.

In all cases the left action of the group $G(Z,r)$ splits $Z$ into
three orbits, i.e. 3 invariant subsets, namely $O_1= {\bf x}\sqcup
{\bf y}, O_2={\bf z}, O_3=Y$. We now give three concrete
extensions, with graphs presented in Figure~4.

\textbf{1)} $(Z,r_1)$ is determined by the actions
\[
\Lcal_{\alpha_i \mid X}=(x_1 y_1 x_2 y_2 x_3 y_3 x_4 y_4) \circ
\tau, \quad \Lcal_{\beta _i \mid X} = (x_1 y_1 x_2 y_2 x_3 y_3 x_4
y_4) \circ \tau^{-1},
\]
\[ \Lcal_{x_j
\mid Y}=\Lcal_{y_j \mid Y}= (\alpha _1 \beta _1 \alpha _2 \beta _2
\alpha _3 \beta _3), \; \Lcal_{z_j \mid Y}= f\circ g \quad ,
\]
for all  $1 \leq i \leq 3, 1\leq j \leq 4.$

 \textbf{2)} $(Z,r_2)$
is determined by the actions
\[
\Lcal_{\alpha_i \mid X}=(x_1 y_1 x_3 y_3) \circ (x_2 y_2 x_4 y_4)
\circ (z_1 z_3)(z_2 z_4), \quad \Lcal_{\beta _i \mid X} =(x_1 y_1
x_3 y_3) \circ (x_2 y_2  x_4 y_4)\circ \id_{\bf z},
\]
\[ \Lcal_{x_j
\mid Y}=\Lcal_{y_j \mid Y}= (\alpha _1\beta _1)\circ(\alpha _2
\beta _2)\circ(\alpha _3 \beta _3), \quad \Lcal_{z_j \mid Y}=
(f\circ g)^2,
\]
for all  $1 \leq i \leq 3, 1\leq j \leq 4.$

\textbf{3)} $(Z,r_3)$ is determined by the actions
\[
\Lcal_{\alpha_i \mid X}= \Lcal_{\beta _i \mid X} =(x_1 y_1)(x_2
y_2)(x_3 y_3)(x_4 y_4)(z_1z_2)(z_3z_4),
\]
\[ \Lcal_{x_j
\mid Y}=\Lcal_{y_j \mid Y}=(\beta _3 \alpha _3 \beta _2 \alpha _2
\beta _1 \alpha _1),\quad
 \Lcal_{z_j \mid Y}= \id_Y,
\]
for all  $1 \leq i \leq 3, 1\leq j \leq 4.$

We make some comments on the solutions and their graphs.
 For arbitrary solution $(X,r)$ with \textbf{lri}, we define the
graph $\Gamma=\Gamma (X,r)$ as follows. It is an oriented graph,
which reflects the left action of $G(X,r)$ on $X$. The set of
vertices of $\Gamma$ is exactly $X.$ There is a labelled arrow  $x
{\buildrel z \over \longrightarrow}y,$ if $x,y,z \in X, x\neq y$
and ${}^zx = y.$ Clearly $x{\buildrel z\over \longleftrightarrow}
y$ indicates that ${}^zx = y$ and ${}^zy = x.$ (One can make such
graph for arbitrary solutions but then it should be indicated
which action is considered). The graphs $\Gamma(X,r_X),
\Gamma(Y,r_Y)$, and $\Gamma(Z,r)$ for the above three extensions
are presented in Figure~4. To avoid clutter we typically omit
self-loops unless needed for clarity or contrast (for example
$\Gamma(Z,r_2)$  shows these explicitly to indicate
${}^{\beta}z_i=z_i$). Also for the same reason, we use the line
type to indicate when the same type of element acts, rather than
labelling every arrow.

Moreover, these  extensions are non-isomorphic. This can be  read
directly from the choice of the actions, but also from the graphs
$\Gamma(Y,r_i)$. Note that two solutions
 are isomorphic if and only if their oriented
graphs    are isomorphic.

In the cases \textbf{1)}, and \textbf{2)}  $G(Y, r_Y)$ does not
act as automorphisms on $X$. For $(Z, r_1)$ this follows from the
\[
{}^{\alpha_k}{(x_iz_j)} = {}^{\alpha_k}{x_i}
{}^{{\alpha_k}^{x_i}}{z_j} ={}^{\alpha_k}{x_i} {}^{\beta
_{k-1}}{z_j} = y_j z_{j-1} \neq {}^{\alpha_k}{x_i}
{}^{\alpha_k}{z_j}= y_j z_{j+1}.
\]
In the second case it is also easy to verify that
${}^{\alpha_k}{(x_iz_j)} \neq {}^{\alpha_k}{x_i}
{}^{\alpha_k}{z_j}.$ In all cases the group $G(X, r_x)$ acts as
automorphisms on $Y.$ In case \textbf{3)} $G(Z, r_3)$ acts as
automorphisms on $Z.$
\end{example}

 {\bf Acknowledgments}. The first author thanks The Abdus Salam
International Centre for Theoretical Physics, where she worked on
the paper in the summers of  2004 and 2005. It is her pleasant
duty to thank Professor Le Dung Trang and  the Mathematics group
of ICTP for the inspiring and creative atmosphere during her
visits. Thanks also to
 QMUL for partial support during her two visits there.

 \ifx\undefined\bysame
\newcommand{\bysame}{\leavevmode\hbox to3em{\hrulefill}\,}
\fi

\end{document}